\documentclass[12pt,a4paper,titlepage]{book}
\usepackage{amsmath,amsfonts,amssymb,bbm,mathrsfs,theorem}
\usepackage{t1enc}
\usepackage[latin1]{inputenc}
\usepackage[english]{babel}
\usepackage{epsfig}
\usepackage{float}
\usepackage{times}
\usepackage{color}
\usepackage{exscale}

\newcommand{\CT}{\ensuremath{\mathscr{C}}}
\newcommand{\C}{\ensuremath{\mathbbm{C}}}

\newcommand{\D}{\ensuremath{\mathscr{D}}}

\newcommand{\N}{\ensuremath{\mathbbm{N}}}
\newcommand{\R}{\ensuremath{\mathbbm{R}}}
\newcommand{\Z}{\ensuremath{\mathbbm{Z}}}
\newcommand{\Rs}{\ensuremath{R^\ast}}
\newcommand{\laplace}{\ensuremath{\triangle}}

\newcommand{\gS}{\ensuremath{\mathscr{S}}}
\renewcommand{\(}{\left(}
\renewcommand{\)}{\right)}
\renewcommand{\exp}{\mbox{e}}
\newcommand{\uint}{\int\limits}
\newcommand{\usum}{\sum\limits}
\newcommand{\ph}{\ensuremath{\varphi}}
\newcommand{\pa}{\partial}

\newcommand{\geg}[2]{\raisebox{-2.7mm}{ $\stackrel{\longrightarrow}{{\scriptstyle #1 \to #2}}$ }}
\newcommand{\bn}{\ensuremath{\mbox{J}_0}}

\newcommand{\dkbild}[3]{
\begin{figure}[!hbp]
\begin{picture}(390,183)
\put(5,10){\includegraphics[width=6cm]{BilderDiplomkorrektur/#1.eps}}
\put(200,10){\includegraphics[width=6cm]{BilderDiplomkorrektur/#2.eps}}
\put(5,0){\small $0$}
\put(81,0){\small $128$}
\put(86,-10){\small $y$}
\put(160,0){\small $256$}
\put(-14,11){\small $128$}
\put(-25,91){\small $x$}
\put(-7,91){\small $0$}
\put(-21,172){\small $-128$}
\put(371,92){\small $y$}
\put(292,178){\textcolor{white}{\circle*{10}}}
\put(295,185){\small $f(0,y)$}
\end{picture}
\caption{#3}
\end{figure}
}

\newcommand{\dkdbild}[4]{
\begin{figure}[!htbp]
\begin{picture}(390,183)
\put(5,10){\includegraphics[width=6cm]{BilderDiplomkorrektur/#1.eps}}
\put(200,10){\includegraphics[width=6cm]{BilderDiplomkorrektur/#1Skala.eps}}
\put(5,0){\small $0$}
\put(81,0){\small $128$}
\put(86,-10){\small $y$}
\put(160,0){\small $256$}
\put(-14,11){\small $128$}
\put(-25,91){\small $x$}
\put(-7,91){\small $0$}
\put(-21,172){\small $-128$}
\put(371,92){\small $y$}
\put(295,185){\small $f(0,y)$}
\end{picture}
\caption{#2}
\vspace{\floatsep}
\vspace{\floatsep}
\begin{picture}(390,183)
\put(5,10){\includegraphics[width=6cm]{BilderDiplomkorrektur/#3.eps}}
\put(200,10){\includegraphics[width=6cm]{BilderDiplomkorrektur/#3Skala.eps}}
\put(5,0){\small $0$}
\put(81,0){\small $128$}
\put(86,-10){\small $y$}
\put(160,0){\small $256$}
\put(-14,11){\small $128$}
\put(-25,91){\small $x$}
\put(-7,91){\small $0$}
\put(-21,172){\small $-128$}
\put(371,92){\small $y$}
\put(292,178){\textcolor{white}{\circle*{9}}}
\put(295,185){\small $f(0,y)$}
\end{picture}
\caption{#4}
\end{figure}
}

\newcommand{\dkdbildp}[4]{
\begin{figure}[!p]
\begin{picture}(390,183)
\put(5,10){\includegraphics[width=6cm]{BilderDiplomkorrektur/#1.eps}}
\put(200,10){\includegraphics[width=6cm]{BilderDiplomkorrektur/#1Skala.eps}}
\put(5,0){\small $0$}
\put(81,0){\small $128$}
\put(86,-10){\small $y$}
\put(160,0){\small $256$}
\put(-14,11){\small $128$}
\put(-25,91){\small $x$}
\put(-7,91){\small $0$}
\put(-21,172){\small $-128$}
\put(371,92){\small $y$}
\put(295,185){\small $f(0,y)$}
\end{picture}
\caption{#2}
\vspace{\floatsep}
\vspace{\floatsep}
\begin{picture}(390,183)
\put(5,10){\includegraphics[width=6cm]{BilderDiplomkorrektur/#3.eps}}
\put(200,10){\includegraphics[width=6cm]{BilderDiplomkorrektur/#3Skala.eps}}
\put(5,0){\small $0$}
\put(81,0){\small $128$}
\put(86,-10){\small $y$}
\put(160,0){\small $256$}
\put(-14,11){\small $128$}
\put(-25,91){\small $x$}
\put(-7,91){\small $0$}
\put(-21,172){\small $-128$}
\put(371,92){\small $y$}
\put(292,178){\textcolor{white}{\circle*{9}}}
\put(295,185){\small $f(0,y)$}
\end{picture}
\caption{#4}
\end{figure}
}

\newcommand{\dkdbildt}[4]{
\begin{figure}[!th]
\begin{picture}(390,183)
\put(5,10){\includegraphics[width=6cm]{BilderDiplomkorrektur/#1.eps}}
\put(200,10){\includegraphics[width=6cm]{BilderDiplomkorrektur/#1Skala.eps}}
\put(5,0){\small $0$}
\put(81,0){\small $128$}
\put(86,-10){\small $y$}
\put(160,0){\small $256$}
\put(-14,11){\small $128$}
\put(-25,91){\small $x$}
\put(-7,91){\small $0$}
\put(-21,172){\small $-128$}
\put(371,92){\small $y$}
\put(295,185){\small $f(0,y)$}
\end{picture}
\caption{#2}
\vspace{\floatsep}
\vspace{\floatsep}
\begin{picture}(390,183)
\put(5,10){\includegraphics[width=6cm]{BilderDiplomkorrektur/#3.eps}}
\put(200,10){\includegraphics[width=6cm]{BilderDiplomkorrektur/#3Skala.eps}}
\put(5,0){\small $0$}
\put(81,0){\small $128$}
\put(86,-10){\small $y$}
\put(160,0){\small $256$}
\put(-14,11){\small $128$}
\put(-25,91){\small $x$}
\put(-7,91){\small $0$}
\put(-21,172){\small $-128$}
\put(371,92){\small $y$}
\put(292,178){\textcolor{white}{\circle*{9}}}
\put(295,185){\small $f(0,y)$}
\end{picture}
\caption{#4}
\end{figure}
}

\newcommand{\asbildg}[3]{
\begin{figure}[!htbp]
\begin{picture}(390,470)
\put(30,10){\includegraphics[width=5cm]{BilderAntisym/#1.eps}}
\put(195,10){\includegraphics[width=5cm]{BilderAntisym/#2.eps}}
\put(30,152){\textcolor{white}{\line(1,0){160}}}
\put(30,294){\textcolor{white}{\line(1,0){160}}}
\put(10,75){\small $x$}
\put(30,0){\small $0$}
\put(90,0){\small $128$}
\put(95,-10){\small $r$}
\put(155,0){\small $256$}
\put(4,142){\small $-128$}
\put(23,75){\small $0$}
\put(13,11){\small $128$}
\put(4,284){\small $-128$}
\put(23,217){\small $0$}
\put(13,154){\small $128$}
\put(4,427){\small $-128$}
\put(23,361){\small $0$}
\put(13,296){\small $128$}
\put(10,217){\small $x$}
\put(10,361){\small $x$}
\put(35,426){\small \textcolor{white}{$g_0(x,r)$}}
\put(35,283){\small \textcolor{white}{$g_2(x,r)$}}
\put(35,140){\small \textcolor{white}{$g_5(x,r)$}}
\put(340,220){\small $r$}
\put(230,440){\small $g_0(0,r), g_2(0,r), g_5(0,r)$}
\end{picture}
\caption{#3}
\end{figure}
}

\newcommand{\asbildhg}[3]{
\begin{figure}[!htbp]
\begin{picture}(390,290)
\put(15,10){\includegraphics[width=6cm]{BilderAntisym/#1.eps}}
\put(195,10){\includegraphics[width=6cm]{BilderAntisym/#2.eps}}
\put(15,93){\textcolor{white}{\line(1,0){180}}}
\put(15,181){\textcolor{white}{\line(1,0){180}}}
\put(5,0){\small $-256$}
\put(100,0){\small $0$}
\put(170,0){\small $256$}
\put(100,-11){\small $y$}
\put(-2,11){\small $128$}
\put(8,48){\small $0$}
\put(-11,84){\small $-128$}
\put(-7,48){\small $x$}
\put(-2,95){\small $128$}
\put(8,134){\small $0$}
\put(-11,171){\small $-128$}
\put(-7,134){\small $x$}
\put(-2,182){\small $128$}
\put(8,219){\small $0$}
\put(-11,257){\small $-128$}
\put(-7,219){\small $x$}
\put(20,81){\small \textcolor{white}{$f_5^e(x,y)$}}
\put(20,168){\small \textcolor{white}{$f_2^e(x,y)$}}
\put(20,255){\small \textcolor{white}{$f_0^e(x,y)$}}
\put(366,135){\small $y$}
\put(255,270){\small $f_0^e(0,y), f_2^e(0,y), f_5^e(0,y)$}
\end{picture}
\caption{#3}
\end{figure}
}

\newcommand{\asbilds}[3]{
\begin{figure}[!htbp]
\begin{picture}(390,100)
\put(0,10){\includegraphics[width=6cm]{BilderAntisym/#1.eps}}
\put(195,10){\includegraphics[width=6cm]{BilderAntisym/#2.eps}}
\put(-10,0){\small $-256$}
\put(85,0){\small $0$}
\put(155,0){\small $256$}
\put(85,-11){\small $y$}
\put(-17,11){\small $128$}
\put(-7,48){\small $0$}
\put(-25,86){\small $-128$}
\put(-22,48){\small $x$}
\put(366,50){\small $y$}
\put(282,95){\textcolor{white}{\circle*{7}}}
\put(290,96){\small $f(0,y)$}
\end{picture}
\caption{#3}
\end{figure}
}

\newcommand{\asbildz}[4]{
\begin{figure}[!htbp]
\begin{picture}(390,100)
\put(0,10){\includegraphics[width=6cm]{BilderAntisym/#1.eps}}
\put(195,10){\includegraphics[width=6cm]{BilderAntisym/#1Skala.eps}}
\put(-10,0){\small $-256$}
\put(85,0){\small $0$}
\put(155,0){\small $256$}
\put(85,-11){\small $y$}
\put(-17,11){\small $128$}
\put(-7,48){\small $0$}
\put(-25,86){\small $-128$}
\put(-22,48){\small $x$}
\put(290,96){\small $f(0,y)$}
\put(366,50){\small $y$}
\end{picture}
\caption{#2}
\vspace{\floatsep}
\begin{picture}(390,100)
\put(0,10){\includegraphics[width=6cm]{BilderAntisym/#3.eps}}
\put(195,10){\includegraphics[width=6cm]{BilderAntisym/#3Skala.eps}}
\put(-8,0){\small $+0$}
\put(32,0){\small $+\frac{\pi}{2}$}
\put(75,0){\small $\pm \pi$}
\put(156,0){\small $-0$}
\put(115,0){\small $-\frac{\pi}{2}$}
\put(85,-11){\small $\eta$}
\put(-17,11){\small $128$}
\put(-7,48){\small $0$}
\put(-25,86){\small $-128$}
\put(-22,48){\small $x$}
\put(290,96){\small $f^{I,F}(0,\eta)$}
\put(366,50){\small $\eta$}
\end{picture}
\caption{#4}
\end{figure}
}

\newcommand{\asbildd}[6]{
\begin{figure}[!htbp]
\begin{picture}(390,100)
\put(0,10){\includegraphics[width=6cm]{BilderAntisym/#1.eps}}
\put(195,10){\includegraphics[width=6cm]{BilderAntisym/#1Skala.eps}}
\put(-10,0){\small $-256$}
\put(85,0){\small $0$}
\put(155,0){\small $256$}
\put(85,-11){\small $y$}
\put(-17,11){\small $128$}
\put(-7,48){\small $0$}
\put(-25,86){\small $-128$}
\put(-22,48){\small $x$}
\put(290,96){\small $f(0,y)$}
\put(366,50){\small $y$}
\end{picture}
\caption{#2}
\vspace{\floatsep}
\begin{picture}(390,100)
\put(0,10){\includegraphics[width=6cm]{BilderAntisym/#3.eps}}
\put(195,10){\includegraphics[width=6cm]{BilderAntisym/#3Skala.eps}}
\put(-10,0){\small $-256$}
\put(85,0){\small $0$}
\put(155,0){\small $256$}
\put(85,-11){\small $y$}
\put(-17,11){\small $128$}
\put(-7,48){\small $0$}
\put(-25,86){\small $-128$}
\put(-22,48){\small $x$}
\put(290,96){\small $f(0,y)$}
\put(366,50){\small $y$}
\end{picture}
\caption{#4}
\vspace{\floatsep}
\begin{picture}(390,100)
\put(0,10){\includegraphics[width=6cm]{BilderAntisym/#5.eps}}
\put(195,10){\includegraphics[width=6cm]{BilderAntisym/#5Skala.eps}}
\put(-10,0){\small $-256$}
\put(85,0){\small $0$}
\put(155,0){\small $256$}
\put(85,-11){\small $y$}
\put(-17,11){\small $128$}
\put(-7,48){\small $0$}
\put(-25,86){\small $-128$}
\put(-22,48){\small $x$}
\put(290,96){\small $f(0,y)$}
\put(366,50){\small $y$}
\end{picture}
\caption{#6}
\end{figure}
}

\newcommand{\fbildr}[2]{
\begin{figure}[!htbp]
\begin{picture}(390,195)
\put(5,15){\includegraphics[width=6cm,angle=90]{BilderFehler/fehlerrange-#1.eps}}
\put(200,15){\includegraphics[width=6cm]{BilderFehler/fehlerrange-#1Skala.eps}}
\put(5,6){\small $0$}
\put(83,4){\small $\frac{R}{2}$}
\put(83,-10){\small $y$}
\put(169,5){\small $R$}
\put(-4,16){\small $0$}
\put(-25,98){\small $x$}
\put(-7,98){\small $\frac{L}{2}$}
\put(-5,177){\small $L$}
\put(371,97){\small $y$}
\put(295,190){\small $f(\frac{1}{2},y)$}
\end{picture}
\caption{#2}
\end{figure}
}

\newcommand{\fbildt}[2]{
\begin{figure}[!htbp]
\begin{picture}(390,195)
\put(5,15){\includegraphics[width=6cm,angle=90]{BilderFehler/fehlertrackneu-#1.eps}}
\put(200,15){\includegraphics[width=6cm]{BilderFehler/fehlertrackneu-#1Skala.eps}}
\put(5,6){\small $0$}
\put(83,4){\small $\frac{R}{2}$}
\put(83,-10){\small $y$}
\put(169,5){\small $R$}
\put(-4,16){\small $0$}
\put(-25,98){\small $x$}
\put(-7,98){\small $\frac{L}{2}$}
\put(-5,177){\small $L$}
\put(371,97){\small $y$}
\put(295,190){\small $f(\frac{1}{3},y),f(\frac{2}{3},y)$}
\end{picture}
\caption{#2}
\end{figure}
}

\theoremstyle{changebreak}

\theoremheaderfont{\bfseries\scshape}

\newtheorem{defi}{Definition}[section]
\newtheorem{lem}[defi]{Lemma}
\newtheorem{prop}[defi]{Proposition}
\newtheorem{rem}[defi]{Remark}
\newtheorem{cor}[defi]{Corollary}
\theorembodyfont{\itshape}
\newtheorem{theo}[defi]{Theorem}
\newenvironment{proof}{\noindent{\bfseries Proof:}\par\nopagebreak}{\nopagebreak\hspace*{\fill} $
\Box$}

\begin{document}


\begin{sloppypar}






















\pagestyle{empty}

\begin{center}

\vspace*{2cm}

{\bf {\LARGE Jens Klein}}\\

\vspace*{3cm}

{\bf {\LARGE {Mathematical Problems in Synthetic Aperture Radar}}}\\

\vspace*{3cm}

{\bf {\LARGE 2004}}\\

\vspace*{5cm}


\vspace*{1cm}

{\bf {\LARGE Universit\"at M\"unster}}\\

\end{center}

\cleardoublepage

\textcolor{white}{o}

\cleardoublepage

\begin{center}


{\large Angewandte Mathematik}\\

\vspace*{3cm}

{\bf {\LARGE {Mathematical problems in Synthetic Aperture Radar}}}\\

\vspace*{3cm}

{\large Inaugural-Dissertation zur Erlangung des Doktorgrades der Naturwissenschaften im Fachbereich Mathematik und Informatik der Mathematisch-Naturwissenschaftlichen Fakult\"at der Westf\"alischen Wilhelms Universit\"at M\"unster}\\

\vspace*{3cm}

vorgelegt von\\

\vspace*{1cm}

{\Large Jens Klein}\\

\vspace*{1cm}

{\large aus Aachen}\\

\vspace*{3cm}

{\Large 2004}

\end{center}

\pagebreak

\textcolor{white}{o}

\vspace*{15cm}

\begin{tabbing}
Dekan:\hspace{7cm}\=Prof. Dr. K. Hinrichs\\
\\
Erster Gutachter:\>Prof. Dr. Dr. h. c. F. Natterer\\
\\
Zweiter Gutachter:\>Prof. Dr. C. W. Cryer\\
\\
Tag der m\"undlichen Pr\"ufung:\\
\\
Tag der Promotion:
\end{tabbing}

\cleardoublepage

\vspace*{8cm}

\begin{center}

For Birgit

\end{center}

\cleardoublepage

\textbf{\Huge Abstract}

\vspace*{2cm}

This thesis is concerned with problems related to Synthetic Aperture Radar (SAR), a technique of making images of the surfaces of planets using electromagnetic waves. Reconstructing images of the surfaces from the gathered data is an inverse problem as are other young and thriving imaging techniques used for example in optical tomography and transient elastography. In optical tomography one tries to create images of the human body employing light, whereas in transient elastography ultrasound is used to measure the propagation of shear waves and thus reconstruct the stiffness of human tissue. But the field of inverse problems covers older and established topics as well like computerized tomography that creates images of the human body by means of x-rays and magnetic resonance imaging using electromagnetic fields to measure the distribution of atoms. All these very different applications have in common that the gathered data is difficult to interpret. Therefore mathematical processing is necessary in order to create an intelligible image of the measured object. Some of the problems related to this mathematical processing necessary in creating SAR-images are analyzed in this thesis.

The thesis is structured as follows:
The first chapter explains what SAR is, and the physical and mathematical background is illuminated.

The following chapter points out a problem with a divergent integral in a common approach and proposes an improvement. Some numerical comparisons are shown that indicate that the improvements allow for a superior image quality.

Thereafter an important problem is analyzed - the problem of limited data. In a realistic SAR-measurement the data gathered from the electromagnetic waves reflected from the surface can only be collected from a limited area. However the reconstruction formula requires data from an infinite distance. The chapter gives a comprehensive analysis of the artifacts which can obscure the reconstructed images due to this problem. Additionally, some numerical examples are shown that point to the severity of the problem.

In chapter \ref{Fast} the fact that data is available only from a limited area is used to propose a new inversion formula. This inversion formula has the potential to make it easier to suppress artifacts due to limited data and, depending on the application, can be refined to a fast reconstruction formula.

In the penultimate chapter a solution to the problem of left-right ambiguity is presented. This problem exists since the invention of SAR and is caused by the geometry of the measurements. This leads to the fact that only symmetric images can be obtained. With the solution from this chapter it is possible to reconstruct not only the even part of the reflectivity function, but also the odd part, thus making it possible to reconstruct asymmetric images. Numerical simulations are shown to demonstrate that this solution is not affected by stability problems as other approaches have been.

The final chapter lists some conclusions drawn from the preceding chapters and develops some continuative ideas that could be pursued in the future.\\

\textbf{Acknowledgements}

First, I would like to thank Prof. Dr. Dr. h. c. F. Natterer for stimulating the work on this thesis and for offering helpful advice.\\
My thank also goes to Dr. F. W\"ubbeling who contributed to this thesis in many fruitful discussions and to all members of the institute for the pleasant atmosphere.\\
Finally I would like to especially thank my parents and my significant other, \mbox{Birgit}, for their ceaseless support that made this thesis possible.

\clearpage

\pagestyle{plain}

\pagenumbering{Roman}

\tableofcontents

\clearpage

\pagenumbering{arabic}

\chapter{Introduction}

This chapter will give a short introduction into the mathematical and physical concepts necessary to deal with SAR. After a closer look at the mathematical terms of well-posed versus ill-posed and direct versus inverse problem, the second section will describe the physical and technical background related to SAR. The last section will then give an introduction to the mathematical model of the inverse problem involved in SAR.

\section{Inverse problems}

From a mathematician's point of view the problems associated with SAR belong to the field of inverse problems. It is difficult to exactly define the term inverse problem. Therefore to begin with, the closely related ideas of  well- and ill-posedness and the concept of stability will be illuminated. To this end a trip into history might prove helpful. The scientific community was aware of the problem of instability as early as 1873 when Maxwell wrote in an essay from February 11th, 1873 \cite[p. 434]{Maxwell}:\\
{\it "There are certain classes of phenomena, as I have said, in which a small error in the data only introduces a small error in the result. ... The course of events in these cases is stable. There are other classes of phenomena which are more complicated, and in which cases of instability may occur, the number of such cases increasing, in an exceedingly rapid manner, as the number of variables increases."}\\
 The concept of a well-posed problem was first formulated by Hadamard in 1902 \cite{Hadamard}. A well-posed problem according to Hadamard requires the existence, uniqueness and stability of a solution, originating from the philosophy that the mathematical model of a physical problem has to have these properties. If one of the properties fails to hold, a problem is called ill-posed. Due to this background only well-posed problems were studied extensively for quite some time and ill-posed problems were neglected.\\
This led to a vast knowledge of and familiarity with direct problems which are usually well-posed. These are problems that are common in physics with an input $x$ and an operator $H$ modelling the physical system. In these cases it is the goal to find the system's response $y=Hx$. Therefore they are called direct problems as it is the classical physical problem to have a system state and let it evolve according to a certain formula. The term inverse problem is derived from this notion of direct problems. The corresponding inverse problems consist in determining the cause $x$ for a known system $H$ and a known response $y$ or in determining parameters of the system $H$ for a known input $x$ and a known response $y$. Unlike direct problems that are usually well-posed, inverse problems are often ill-posed. Inverse problems are relevant for many different physical applications, for example problems related to imaging. One of these applications is SAR which will be described in detail in the following.

\section{Physical background}

\begin{figure}[!hbt]
\begin{picture}(390,170)
\put(0,-10){\includegraphics[width=6.5cm]{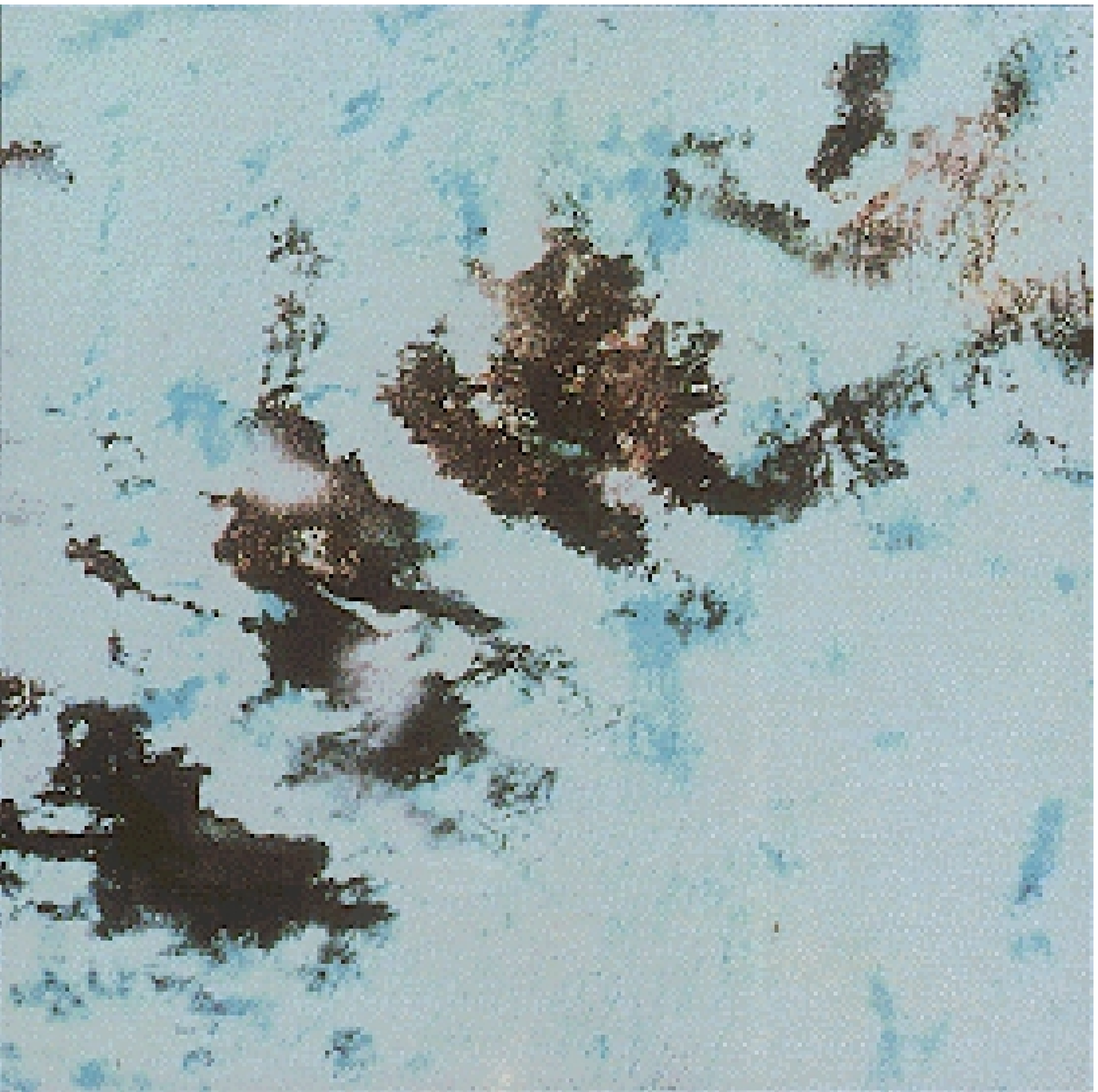}}
\put(195,-10){\includegraphics[width=6.5cm]{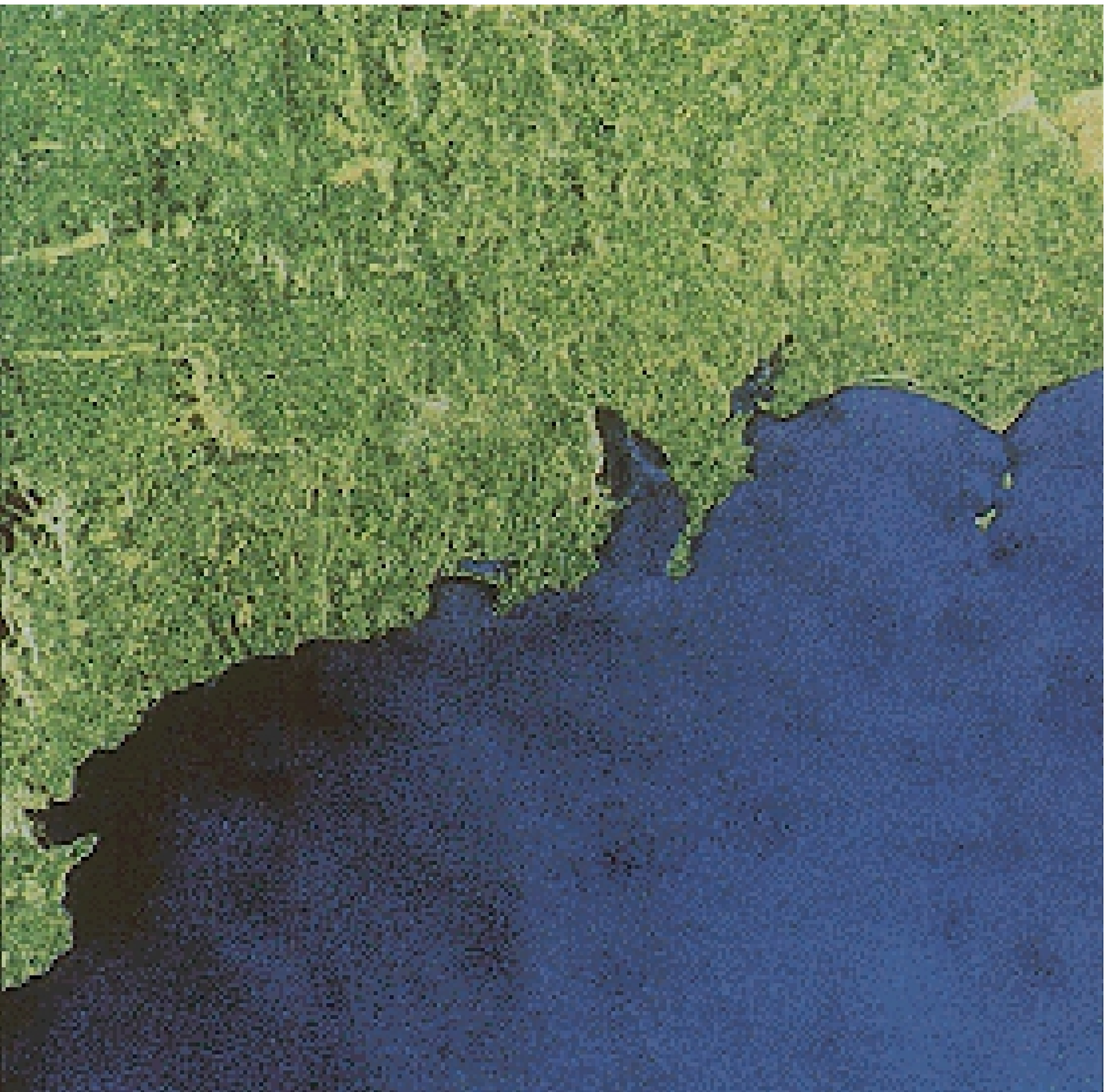}}
\end{picture}
\caption{Waterford (Ireland) on the morning of August 9th, 1991. Left: Optical Landsat satellite image. Right: SAR ERS-1 satellite image. Source: \cite{ESA}\label{fog}}
\end{figure}

Synthetic Aperture Radar is a technique of taking pictures of the surface of planets from an airplane or satellite. As the name indicates, SAR utilizes Radar, electromagnetic waves with a much longer wavelength than used in optical imaging. These waves are emitted by antennas mounted on the airplane or satellite. They are then reflected from the surface and detected by the same antennas. The use of such large wavelength leads to a great advantage in comparison with photography. Since waves of these wavelengths can penetrate clouds and even foliage, SAR images can be taken in foggy or cloudy weather. This capability is demonstrated in figure \ref{fog} that shows in both images the region of Waterford in Ireland on the morning of August 9th, 1991. The left picture shows an optical image that is almost completely obscured by clouds whereas the SAR image on the right is not at all affected by them. This is an important feature since in Europe only one out of ten optical images is free of clouds \cite{ESA}. However, this comes at the cost of a reduced resolution compared with optical imaging.

The typical frequency used in SAR is largely dependent on the technical implementation. Commonly frequencies ranging from 20 MHz to 10 GHz are used \cite{Gustavsson}, \cite{Mirkin}. This corresponds to wavelengths of 15 m to 3 cm. They even permit to detect concealed object, e. g. covered by trees, and to measure the biomass of a region \cite{Gustavsson}.\\
An example of a SAR system is CARABAS. This is an airborne VHF SAR system developed by the FOA (National Defence Research Establishment) in Sweden. CARABAS consists of two $5.5$ m long antennas mounted parallel on a Rockwell Sabreliner aircraft. The distance of the antennas corresponds to the shortest emitted wavelength and translates to three meters since the used frequency spectrum ranges from $20$ to $90$ MHz. The airplane usually travels at a height of $1500$ to $10000$~m at a speed of $100$ to $130 \frac{\mbox{m}}{\mbox{s}}$ \cite{Gustavsson}. For a frequency of $70$ MHz theoretically the best achievable resolution is $1$ m parallel to the flight track and $2$ m perpendicular to the flight track \cite{Hellsten3}. However it turned out that for real measurements the resolution obtainable is only half as good \cite{Hellsten}.

\section{Mathematical model}

The correct mathematical model for radar emission and scattering is given by Maxwell's equations of electromagnetism. Therefore Maxwell's equations are used to understand certain effects that are common in SAR and to enhance the reconstruction \cite{URSI}. But due to the complexity of deriving a reconstruction formula based directly on these vector equations, usually only a scalar wave model is used in mathematical models \cite{Cheney}, \cite{Kostinski}. In today's real life applications not even algorithms based on the wave equation are commonly used, but a simple backprojection algorithm is applied that is now explained in more detail.\\

\begin{figure}[!htb]
\begin{picture}(390,170)
\put(0,-10){\includegraphics[width=14cm]{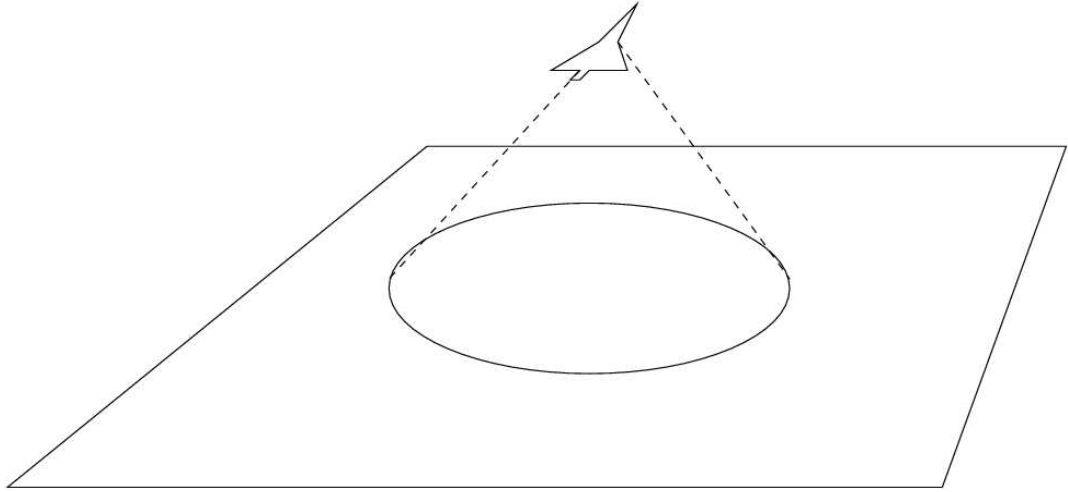}}
\end{picture}
\caption{The geometry of SAR. Source: \cite{Cheney2}\label{geometry}}
\end{figure}

For the simplest case assume a point source that is flying along a straight line above a flat plain and emitting infinitely short pulses in constant time intervals. For a single pulse this results in an expanding sphere centered at the emitting position. When the expanding sphere is large enough, it hits the ground and is reflected. Since the speed of light in air is much larger than the speed of the airplane, the airplane is considered static during each send-receive cycle (start-stop approximation). As the scattered wave is therefore received at exactly the same position at which it was emitted, the received signal comprises an integral of the ground reflectivity over the circle that is obtained intersecting the plane and the expanding sphere at a certain point in time as depicted in figure \ref{geometry}. The equation that describes this measurement is
$$g(x,r)=Rf(x,r)=\frac{1}{|S^1|}\uint_{S^1}f(x+r\xi,r\eta) \, d S_1(\xi,\eta)$$
where $g$ is the measured data and $f$ is the function to be reconstructed (this notation will also be adopted in the following chapters). $S^1$ denotes the unit sphere in $\R^2$, $d S_1$ stands for the canonical surface measure on $S^1$, $r$ denotes the radius of the circle of the intersection of the plane and the expanding sphere, and $x$ is the position of the airplane on the track projected onto the plane. This equation is called the spherical Radon transform. It is obvious that the collected data contains only the even part of the reflectivity function with respect to the flight track as the whole setup is left-right symmetric. The inverse problem now consists of recovering $f$, either analytically by finding a reconstruction formula or numerically. The easiest way to do so is a simple backprojection which numerically distributes each point of collected data equally over all points of the reconstructed image that lie on the circle where the data comes from.\\
The first thorough mathematical treatment based on this model was developed by Andersson \cite{Andersson}. Chapter \ref{Diplomkorrektur} is based solely on these ideas and analyzes them more thoroughly. This model is also the foundation of chapters \ref{Fehler} and \ref{Fast}, but chapter \ref{Fehler} should also yield insight into a more general problem associated with SAR and the ideas of chapter \ref{Fast} may be transferred onto other approaches to SAR. In chapter \ref{Antisym} the algorithms based on this simple model are only used for numerical examples, whereas the underlying theory developed in this chapter is much more general. Therefore the ideas of this chapter should be applicable to all SAR systems.




\chapter{An improved inversion formula for the spherical Radon transform}

\label{Diplomkorrektur}

The determination of a function from spherical averages is a problem often encountered in physical applications such as SAR and SONAR (SOund Navigation And Ranging is a technique that uses sound propagation under water to navigate or to detect other watercraft). The work related to this topic, which has lead to a great amount of insight and refinement today, began with the proposal of a reconstruction formula by Fawcett \cite{Fawcett}. The mathematical analysis of the problem was later improved by Andersson \cite{Andersson}, and two refined reconstruction formulas were derived. This sparked a host of activity \cite{Ulander1}, \cite{Hellsten}, \cite{Nolan}, \cite{Borden}, \cite{Louis} so that much of today's research is based on Anderssons's ideas.\\

In \cite{Andersson} two reconstruction formulas were derived from the Fourier inversion formula, but it was not checked whether they are properly defined. In the following it is shown that the first contains an integral that diverges under physically sensible conditions. An alternative is presented. Additionally it is shown that the other reconstruction formula might be difficult to compute numerically.

\section{Introduction}

At first, for the benefit of the reader, some results from \cite{Andersson} will be recalled. Note that besides the aforementioned problems, a few minor errors occurred, which do not essentially obscure the results in \cite{Andersson}. For a detailed analysis of these errors see \cite{Klein}. Additionally it will be shown that an integral in the first reconstruction formula in \cite{Andersson} does not converge.

\begin{defi}

Let $n \in \N$. Then

\begin{enumerate}

\item

$\N_0:=\N\cup\{0\}$

\item

$\gS(\R^{n}):=\{ \ph \in C^\infty(\R^{n}) : \|x^\beta \pa^\alpha \ph\|_\infty < \infty \, \forall \, \beta, \alpha \in \N_0^{n}\}$ is the Schwartz Space.

\item

A sequence $(\ph_l)_{l \in \N} \subseteq \gS(\R^{n})$ is called zero sequence if $\|x^\alpha \pa^\beta \ph_l\|_\infty \geg{l}{\infty} 0 \forall \beta, \alpha \in \N_0^n$.

\item

A linear functional $f: \gS(\R^n) \rightarrow \C$, $\ph \mapsto <f,\ph>$ is continuous if $<f,\ph_l> \geg{l}{\infty} 0$ for all zero sequences $(\ph_l)_{l \in \N} \subseteq \gS(\R^n)$.

\item

$\gS_e(\R^{n+1}):=\{\varphi \in \gS(\R^{n+1}): \varphi(x,-y)=\varphi(x,y), \forall x \in \R^n, y \in \R \}$.

\item

$\gS_r(\R^n \times \R^{n+1}):=\{\varphi \in \gS(\R^{2n+1}) : \forall \mbox{ orthonormal transformations}$
$$ U:\R^{n+1} \rightarrow \R^{n+1} \: \forall x \in \R^n , z \in \R^{n+1}: \varphi (x,z) = \varphi (x,Uz)\} \mbox{.}$$

\item

$\gS'_e(\R^{n+1})$ and $\gS'_r(\R^n \times \R^{n+1})$ are the dual spaces of $\gS_e(\R^{n+1})$ and $\gS_r(\R^n \times \R^{n+1})$, respectively. The weak-$\star$ topology is used.

\item

$(\gS'_r)_{cone}(\R^n \times \R^{n+1}):=\{g \in \gS'_r(\R^n \times \R^{n+1}): \mbox{supp }\hat{g} \subseteq \{(\xi,\eta):\|\eta\| \geq \|\xi\|\}\} \mbox{.}$

\item

$S^n$ denotes the unit sphere in $\R^{n+1}$, and $S_n$ denotes the canonical measure on the surface of $S^n$.

\item

Let $f \in \gS(\R^n)$. Then $\hat{f}$ denotes the Fourier transform in $\R^n$:
$$\hat{f}(\xi)=\frac{1}{(2\pi)^\frac{n}{2}}\uint_{\R^n}\exp^{-i <x,\xi>} f(x) \, dx \mbox{.}$$

\item

Let $f \in \gS(\R^n)$ and $g \in \gS'(\R^n)$. Then the Fourier transform of $g$ is defined by
$$<\hat{g},f>=<g,\hat{f}>\mbox{.}$$

\item

Let $f \in \gS_e(\R^{n+1})$. Then the operator $R$ is defined by
$$Rf(x,r):=\frac{1}{|S^n|}\uint_{S^n}f(x+r\xi,r\eta) \, d S_n(\xi,\eta) \mbox{.}$$

\end{enumerate}

\end{defi}

In this definition the operator $R$ describes, in the case of SAR and SONAR, the measurement of the reflectivity function $f$ that represents the ground reflectivity. The measurement is modeled as a $\delta$-impulse wave that propagates as concentric spheres. The ground is approximated as a plane. The single scatter approximation for the wave hitting the ground results in integrals over circles.\\
For simplicity, in the following sometimes only $\gS$, $\gS_e$, etc. is written instead of $\gS(\R^{n+1})$, $\gS_e(\R^{n+1})$, etc.

\begin{rem}

\begin{enumerate}

\item

Let $f \in \gS_e(\R^{n+1})$. Then for $g=Rf \in \gS'_r(\R^n \times \R^{n+1})$, $x \in \R^n$, and $r \geq 0$ in the following $g(x,r)$ is sometimes written with abuse of notation. This identification of $g(x,y)$ with $g(x,|y|)$ for $y \in \R^{n+1}$ is justified because $g$ depends only radially on the last $n+1$ variables.

\item

Since $\gS_r(\R^n \times \R^{n+1})$ and $\gS_e(\R^{n+1})$ are subsets of $\gS(\R^n \times \R^{n+1})$ and $\gS(\R^{n+1})$, the definition of the Fourier transform can be extended to $\gS_r$, $\gS'_r$, $\gS_e$, and $\gS'_e$.

\item

It is easily seen that $\hat{f} \in \gS_e$ and $\hat{g} \in \gS_r$ if $f \in \gS_e$ and $g \in \gS_r$, respectively.

\end{enumerate}

\end{rem}

The essential result in \cite{Andersson} is the (Fourier-) inversion formula:

\begin{theo}

\label{inversionformula}

If $\gS_e(\R^{n+1})$ is given the topology of $\gS'(\R^{n+1})$, the mapping
$$R:\gS_e(\R^{n+1}) \rightarrow \gS'_r(\R^n \times \R^{n+1})$$
is continuous and can, by continuity, be extended to a mapping
$$R:\gS'_e(\R^{n+1}) \rightarrow \gS'_r(\R^n \times \R^{n+1}) \mbox{.}$$
The range of this extended mapping $R$ is the closed subspace
$$(\gS'_r)_{cone}(\R^n \times \R^{n+1}) \subseteq \gS'_r(\R^n \times \R^{n+1}) \mbox{.}$$
$R$ is one-to-one and the inverse mapping
$$R^{-1}:(\gS'_r)_{cone}(\R^n \times \R^{n+1}) \rightarrow \gS'_e(\R^{n+1})$$
is continuous. Moreover, if $g=Rf$ and $\hat{f}(\xi,\eta)$ or $\hat{g}(\xi,\eta)$ are integrable functions for $\xi, \eta$ in some open set, then in that set
$$\hat{g}(\xi,\eta)=\left\{ \begin{array}{lr}
\displaystyle{(2\pi)^n \frac{2}{|S^n|} \frac{\hat{f}\(\xi,\sqrt{\|\eta\|^2-\|\xi\|^2}\)}{\|\eta\|^{n-1} \sqrt{\|\eta\|^2-\|\xi\|^2}}} &
\mbox{for }\|\eta\|>\|\xi\|\\
0 & \mbox{for }0 < \|\eta\| \leq \|\xi\|\\
\end{array} \right.$$
or
$$\hat{f}(\xi,\eta)=\frac{1}{(2\pi)^n} \frac{|S^n|}{2}|\eta|\(\|\xi\|^2+\eta^2\)^{\frac{n-1}{2}}\hat{g}\(\xi,\sqrt{\|\xi\|^2+\eta^2}\) \mbox{.}$$

\end{theo}

\begin{proof}

\cite[Theorem 2.1]{Andersson}

\end{proof}

This means that the Fourier transform $\hat{g}$ of the data $g$ can under certain conditions be used to extract the Fourier-transform $\hat{f}$ of the reflectivity function $f$.

\begin{defi}

\begin{enumerate}

\item

Let $g \in \gS_r(\R^n \times \R^{n+1})$. Then
$$\Rs g(x,y)=\uint_{\R^n} g\(z,\sqrt{\|z-x\|^2+y^2}\) \, dz \mbox{.}$$

\item

Let $f \in L_2(\R)$. Then the Hilbert transform of $f$ is defined by the principal value integral
$$\frac{1}{\pi} \uint_\R \frac{f(y)}{x-y} \, dy \mbox{.}$$

\end{enumerate}

\end{defi}

In \cite{Andersson} two reconstruction formulas were derived from theorem \ref{inversionformula}, the first of which was already essentially given in \cite{Fawcett}. They are given in the following corollary.

\begin{cor}

\label{reformulation}

With $c_n=\frac{1}{(2\pi)^n}\frac{|S^n|}{2}$ two reformulations of the inversion formula are possible:

\begin{enumerate}

\item

For $g \in \gS_r$
$$f=c_n H_y \frac{\pa}{\pa y} \laplace^\frac{n-1}{2} \Rs g$$
with the Hilbert transform in $y$, $H_y$, and the Laplace-Operator $\laplace=\laplace_x+\frac{\pa^2}{\pa y^2}$. This formula is essentially also given by Fawcett \cite{Fawcett}.

\item

For $g \in (\gS_r)_{cone}:=\{g \in \gS_r: \mbox{supp }\hat{g} \subseteq \{(\xi,\eta):\|\eta\| \geq \|\xi\|\}\}$
$$f=c_n \Rs K g$$
with the operator $K$ defined by $\widehat{Kg}(\xi,\eta)=\sqrt{\|\eta\|^2-\|\xi\|^2}\,\|\eta\|^{n-1}\hat{g}(\xi,\eta)$.

\end{enumerate}

\end{cor}

\begin{proof}

\cite[Section 3]{Andersson}

\end{proof}

The corollary states that under certain restrictions it is possible to reconstruct the reflectivity function $f$ directly from the data $g$ without taking the detour through the Fourier space.

\begin{rem}

Note, that the essential restriction of this formulation of the corollary in comparison to \cite{Andersson} is that the data $g$ has to be in $\gS_r$ or $(\gS_r)_{cone}$ respectively. This is necessary because otherwise the application of $\Rs$ to $g$ in the first case is undefined or the application of $K$ to $g$ in the second case.\\
Unfortunately $g=Rf$ is usually not in $\gS_r$. Therefore the two reformulations of the inversion formula are only valid in the distributional sense with an appropriately defined $\Rs$.\\
For example a physically reasonable $f \in \gS_e$, $f:\R^2 \rightarrow [0, \infty)$ with $f(x_0,y_0)>0$ for some $x_0,y_0 \in \R$
yields $f(x,y)>c>0$ for all $(x,y) \in K_\epsilon (x_0,y_0)$ with appropriate $c,\epsilon >0$.\\
Moreover
$$2\pi g\(z,\sqrt{(z-x)^2+y^2}\)=\frac{1}{\sqrt{(z-x)^2+y^2}} \uint_{\|r\|=\sqrt{(z-x)^2+y^2}} f(z+r_1,r_2) \, d\sigma(r)$$
with $r=(r_1,r_2)$ and $\sigma$ the canonical measure on the sphere with radius $r$ in $\R^2$.
For $\sqrt{(z-x_0)^2+y_0^2}>\epsilon$ we obtain with a simple geometrical consideration:
$$2\pi g\(z,\sqrt{(z-x_0)^2+y_0^2}\) \geq \frac{\epsilon c}{2 \sqrt{(z-x_0)^2+y_0^2}} \mbox{.}$$
Thus
\begin{align*}
2\pi \Rs g(x_0,y_0)&=\uint_\R g \( z,\sqrt{(z-x_0)^2+y_0^2} \) \, dz\\
&\geq \frac{\epsilon c}{2} \uint_{\sqrt{(z-x_0)^2+y_0^2}\,>\,\epsilon} \frac{1}{\sqrt{(z-x_0)^2+y_0^2}} \, dz =\infty \mbox{.}
\end{align*}
Since this integral diverges, the reconstruction of a non-negative function $f$ with $f(x_0,y_0)>0$ for some $(x_0,y_0) \in \R^2$ is impossible using the first reconstruction formula.

\end{rem}

This result is in accordance with a result from Nessibi, Rachdi, and Trimeche \cite{Nessibi}. They gave reconstruction formulas for functions $g=Rf$ with
$$\uint_0^\infty P(y) f(x,y) \, dy=0$$
for all $x \in \R^n$ and for all one-variable polynomials $P$.

\section{Properties of the function $g=R f$}

Before showing an important property of the data function $g=Rf$, which will be necessary for the derivation of the new reconstruction formula representing an alternative to Andersson's, some definitions will be needed.

\begin{defi}

Let $n \in \N$. Then

\begin{enumerate}

\item

$\D(\R^{n+1}):=C_0^\infty(\R^{n+1})$.

\item

$\D_e(\R^{n+1}):=\{\varphi \in \D(\R^{n+1}): \varphi(x,-y)=\varphi(x,y), \forall x \in \R^n, y \in \R \}$.

\end{enumerate}

\end{defi}








\begin{lem}

\label{gincinf}

If $f \in \gS_e(\R^{n+1})$, then $g=Rf \in C^\infty$.

\end{lem}

\begin{proof}

The interchange of differentiation and integration is justified because $f \in \gS_e(\R^{n+1})$.

\end{proof}

Now an important property of $g$ can be shown.

\begin{theo}

\label{ghatinl1}

If $f \in \gS_e(\R^{n+1})$, then $\hat{g}=\widehat{Rf} \in L^1(\R^n \times \R^{n+1})$.

\end{theo}

\begin{proof}

With $f \in \gS$, also $\hat{f} \in \gS \subseteq L^1$. Theorem \ref{inversionformula} implies that
$$\hat{g}(\xi,\eta)=\left\{ \begin{array}{lr}
\displaystyle{(2\pi)^n \frac{2}{|S^n|} \frac{\hat{f}\(\xi,\sqrt{\|\eta\|^2-\|\xi\|^2}\)}{\|\eta\|^{n-1} \sqrt{\|\eta\|^2-\|\xi\|^2}}} &
\mbox{for }\|\eta\|>\|\xi\|\\
0 & \mbox{otherwise.}\\
\end{array} \right.$$
Therefore
$$\uint_{\R^n \times \R^{n+1}} |\hat{g}(\xi,\eta)|\,d\xi \, d\eta=(2\pi)^n\frac{2}{|S^n|}\uint_{\|\eta\| \geq \|\xi\|}\left|\frac{\hat{f}\(\xi,\sqrt{\|\eta\|^2-\|\xi\|^2}\)}{\|\eta\|^{n-1} \sqrt{\|\eta\|^2-\|\xi\|^2}}\right|\,d\xi \, d\eta \mbox{.}$$
The substitution $\|\eta\|=\rho'$ results in
$$\uint_{\R^n \times \R^{n+1}} |\hat{g}(\xi,\eta)|\,d\xi \, d\eta=2(2\pi)^n\uint_{\R^{n}}\uint_{\rho' \geq \|\xi\|}\left|\frac{\rho' \hat{f}\(\xi,\sqrt{\rho'^2-\|\xi\|^2}\)}{\sqrt{\rho'^2-\|\xi\|^2}}\right|\,d\xi \, d\rho' \mbox{.}$$
The substitution $\rho'=\rho + \|\xi\|$ leads to
$$\uint_{\R^n \times \R^{n+1}} |\hat{g}(\xi,\eta)|\,d\xi \, d\eta=2(2\pi)^n\uint_{\R^n}\uint_{\rho \geq 0}\left|\frac{(\rho + \|\xi\|) \hat{f}\(\xi,\sqrt{\rho^2+2\rho \|\xi\|}\)}{\sqrt{\rho^2+2\rho \|\xi\|}}\right|\,d\xi \, d\rho \mbox{.}$$
$f$ is in $\gS$, therefore
\begin{align*}\uint_{\R^n \times \R^{n+1}}& |\hat{g}(\xi,\eta)|\,d\xi \, d\eta\\
&\leq 2(2\pi)^n\uint_{\R^n}\uint_{\rho \geq 0}\frac{C (\rho + \|\xi\|)}{\sqrt{\rho^2+2\rho \|\xi\|} \(1+\sqrt{\|\xi\|^2+\rho^2+2\rho \|\xi\|}\)^{n+3}}\,d\xi \, d\rho\\
&=2(2\pi)^n\uint_{\R^n}\uint_{\rho \geq 0}\frac{C (\rho+\|\xi\|)}{\sqrt{\rho^2 + 2\rho \|\xi\|}(1+\|\xi\|+\rho)^{n+3}} \,d\xi \, d\rho\\
&\leq 2(2\pi)^n\uint_{\R^n}\uint_{\rho \geq 0} \frac{C}{\sqrt{\rho^2 + 2\rho \|\xi\|}(1+\|\xi\|+\rho)^{n+2}} \,d\xi \, d\rho\end{align*}
and the integral
$$\uint_{\R^n}\uint_{\rho \geq 1} \frac{C}{\sqrt{\rho^2 + 2\rho \|\xi\|}(1+\|\xi\|+\rho)^{n+2}} \,d\xi \, d\rho$$
converges.
So only the integral
$$\uint_{\R^n}\uint_{0 \leq \rho \leq 1} \frac{C}{\sqrt{\rho^2 + 2\rho \|\xi\|}(1+\|\xi\|+\rho)^{n+2}} \,d\xi \, d\rho$$
remains to be examined.
$$\uint_{\R^n}\uint_{0 \leq \rho \leq 1} \frac{C}{\sqrt{\rho^2 + 2\rho \|\xi\|}(1+\|\xi\|+\rho)^{n+2}} \,d\xi \, d\rho$$
$$\leq \uint_{\R^n}\uint_{0 \leq \rho \leq 1} \frac{C}{\sqrt{\rho \|\xi\|}(1+\|\xi\|)^{n+2}} \,d\xi \, d\rho$$
$$=\uint_{\R^n}\frac{2C}{\sqrt{\|\xi\|}(1+\|\xi\|)^{n+2}} \,d\xi < \infty \mbox{.}$$

\end{proof}

\begin{rem}

An analogous proof shows $\eta \hat{g} \in L^1$.

\end{rem}

\section{Modification of Andersson's first inversion formula}

To derive an inversion formula that overcomes the problem of the diverging integral, a modified operator $\Rs$ is defined. This enables a convenient formulation of a new inversion formula.

\subsection{Definition and properties of a modified $\Rs$}

Now a modified version of the operator $\Rs$ is introduced, and its properties are discussed.

\begin{defi}

\label{rsternpa}

For $f \in \D_e$ and $g=Rf \in C^\infty(\R^n \times \R^{n+1})$ we define
$$(\Rs_{\pa} g)(x,y):=\uint_{\R^n} \frac{\pa}{\pa y} g\(z,\sqrt{\|x-z\|^2 + y^2}\)\, dz \mbox{.}$$

\end{defi}

This slight modification by an additional derivation turns out to ensure the convergence of the integral applied by the operator $\Rs_{\pa}$ under the minor and physically feasible constraint that $f$ is in $\D_e$. Thereby the formulation of a mathematically exact reconstruction formula is possible.

\begin{prop}

\label{rsternpawohldef}

$\Rs_{\pa} g$ is well defined for $f \in \D_e$ and $g=Rf$.

\end{prop}

\begin{proof}

Let $x \in \R^n$, $y \in \R$. Then
\begin{align*}
&\left|(\Rs_{\pa} g)(x,y)\right|=\left|\,\uint_{\R^n} \frac{\pa}{\pa y} g\(z,\sqrt{\|x-z\|^2 + y^2}\)\, dz\right|\\
&=\left|\,\uint_{\R^n} \frac{\pa}{\pa y} \frac{1}{|S^n|} \uint_{S^n} f\(z+\xi\sqrt{\|x-z\|^2 + y^2},\eta\sqrt{\|x-z\|^2 + y^2}\)\, dS_n(\xi,\eta) \, dz\right|\mbox{.}
\end{align*}
$f \in \D_e$ and therefore it is possible to interchange differentiation and integration.
\begin{align*}
&\left|(\Rs_{\pa} g)(x,y)\right|\\
&=\left|\,\uint_{\R^n} \frac{1}{|S^n|} \uint_{S^n} \frac{\pa}{\pa y} f\(z+\xi\sqrt{\|x-z\|^2 + y^2},\eta\sqrt{\|x-z\|^2 + y^2}\)\, dS_n(\xi,\eta) \, dz\right|\\
&=\left|\,\uint_{\R^n} \frac{1}{|S^n|} \uint_{S^n} \frac{y}{\sqrt{\|x-z\|^2 + y^2}} \right.\\
&\quad \left. \times {\xi \choose \eta} \cdot (\nabla f)\(z+\xi\sqrt{\|x-z\|^2 + y^2},\eta\sqrt{\|x-z\|^2 + y^2}\)\, dS_n(\xi,\eta) \, dz\right|\\
&\leq \uint_{\R^n} \frac{1}{|S^n|} \uint_{S^n} \frac{|y|}{\sqrt{\|x-z\|^2 + y^2}}\\
&\quad\times \left\|(\nabla f)\(z+\xi\sqrt{\|x-z\|^2 + y^2},\eta\sqrt{\|x-z\|^2 + y^2}\)\right\|\, dS_n(\xi,\eta) \, dz \mbox{.}
\end{align*}
The substitution $r=(\xi,\eta)\sqrt{\|x-z\|^2 + y^2}$ yields
\begin{align*}
&\left|(\Rs_{\pa} g)(x,y)\right|\\
& \leq \uint_{\R^n} \frac{1}{|S^n|} \uint_{\|r\|=\sqrt{\|x-z\|^2 + y^2}} \frac{|y|}{\(\sqrt{\|x-z\|^2 + y^2}\)^{n+1}} \|(\nabla f)((z,0)+r)\|\, d\sigma(r) \, dz \mbox{.}
\end{align*}
Since $f \in \D$,\\
$|(\Rs_{\pa} g)(x,y)|$
$$\leq \uint_{\R^n} \frac{|y|}{\(\sqrt{\|x-z\|^2 + y^2}\)^{n+1}} \mbox{ max }(\|\nabla f\|) \mbox{ diam}(\mbox{supp } f) \, dz < \infty \mbox{.}$$

\end{proof}

\begin{cor}

\label{rsternpainsstrich}

For $f \in \D_e$ and $g=Rf$, $\Rs_{\pa}g \in C^\infty \cap L^\infty$.

\end{cor}

\begin{proof}

This is guaranteed by the estimate in proposition \ref{rsternpawohldef}.

\end{proof}

\begin{defi}

Let $f \in \gS'(\R^n)$ and $\ph \in \gS(\R^n)$. Then
$$<f,\ph>_{\gS(\R^n)}:=f(\ph)$$
is the functional $f$ applied to the test function $\ph$.
Here the subscript $\gS(\R^n)$ is a reminder of the space of the test function $\ph$.

\end{defi}

Now in analogy to \cite{Andersson} an expression for $\widehat{\Rs_{\pa} g}$ is derived.

\begin{theo}

\label{rsternpagdachgleichetagdach}

Let $f \in \D_e$ and $g=Rf$. Then $\widehat{\Rs_{\pa}g}(\xi,\eta)=i \eta \hat{g}\(\xi,\sqrt{\|\xi\|^2+\eta^2}\)$.

\end{theo}

\begin{proof}

Let $\ph\in \gS$, $f \in \D_e$, and $g=Rf$. Then
$$<\widehat{\Rs_{\pa}g},\ph>_{\gS(\R^n \times \R)}=<\Rs_{\pa}g,\hat{\ph}>_{\gS(\R^n \times \R)}=\uint_{\R^n \times \R} (\Rs_{\pa}g)(x',y)\hat{\ph}(x',y)\,dx' \, dy$$
$$=\uint_{\R^n \times \R} \uint_{\R^n} \frac{\pa}{\pa y} g\(z,\sqrt{\|x'-z\|^2 + y^2}\)\, dz \, \hat{\ph}(x',y)\,dx' \, dy \mbox{.}$$
As $\hat{\ph} \in \gS$ and with corollary \ref{rsternpainsstrich}, Fubini's theorem implies
$$<\widehat{\Rs_{\pa}g},\ph>_{\gS(\R^n \times \R)}=\uint_{\R^n} \uint_{\R^n \times \R} \frac{\pa}{\pa y} g\(z,\sqrt{\|x'-z\|^2 + y^2}\) \hat{\ph}(x',y)\,dx' \, dy \, dz \mbox{.}$$
Continuing in the distributional sense
\begin{align*}
<\widehat{\Rs_{\pa}g},\ph&>_{\gS(\R^n \times \R)}=-\uint_{\R^n} \uint_{\R^n \times \R} g\(z,\sqrt{\|x'-z\|^2 + y^2}\) \frac{\pa}{\pa y} \hat{\ph}(x',y)\,dx' \, dy \, dz\\
&=-\uint_{\R^n} \uint_{\R^n \times \R} g\(z,\sqrt{\|x'-z\|^2 + y^2}\)\\
&\quad\times \frac{\pa}{\pa y} \( \frac{1}{(2\pi)^\frac{n+1}{2}} \uint_{\R^n \times \R} \exp^{-i (<x',\xi>+y\eta)} \ph(\xi,\eta)\,d\xi \, d\eta \) dx' \, dy \, dz\\
&=i \uint_{\R^n} \uint_{\R^n \times \R} \uint_{\R^n \times \R} \frac{1}{(2\pi)^\frac{n+1}{2}} \exp^{-i <z,\xi>} \exp^{-i (<x'-z,\xi>+y\eta)}\\
&\quad\times \eta g\(z,\sqrt{\|x'-z\|^2 + y^2}\) \ph(\xi,\eta)\,d\xi \, d\eta \, dx' \, dy \, dz\mbox{.}
\end{align*}
The substitution $x'=x+z$ yields
\begin{align*}
<\widehat{\Rs_{\pa}g},\ph>_{\gS(\R^n \times \R)}&=i \uint_{\R^n} \uint_{\R^n \times \R} \uint_{\R^n \times \R} \frac{1}{(2\pi)^\frac{n+1}{2}} \exp^{-i <z,\xi>} \exp^{-i (<x,\xi>+y\eta)}\\
&\quad\times \eta g\(z,\sqrt{\|x\|^2 + y^2}\) \ph(\xi,\eta)\,d\xi \, d\eta \, dx \, dy \, dz
\end{align*}
\begin{align*}&=\uint_{\R^n \times \R} i \eta \hat{g}\(\xi,\sqrt{\|\xi\|^2 + \eta^2}\) \ph(\xi,\eta)\,d\xi \, d\eta\\
&=<i \eta \hat{g_1},\ph>_{\gS(\R^n \times \R)}\mbox{.}
\end{align*}
Here $\hat{g_1}(\xi,\eta):=\hat{g}\(\xi,\sqrt{\|\xi\|^2+\eta^2}\)$ and with this definition $\hat{g_1} \in \gS'(\R^n \times \R)$.

\end{proof}

\subsection{A modified inversion formula}

With the results of the preceding sections a well defined reconstruction formula is also attainable for physically meaningful reflectivity functions.

\begin{theo}

\label{modinvform}

Let $f \in \D_e$ and $g=Rf$. Then $f=c_n H_y \laplace^{\frac{n-1}{2}} \Rs_{\pa} g$ with the constant \mbox{$c_n:=\frac{1}{(2\pi)^n}\frac{|S^n|}{2}$} and the Hilbert transform $H_y$.

\end{theo}

\begin{proof}

Let $f \in \D_e$. It follows from theorem \ref{ghatinl1} and theorem \ref{inversionformula} that
\begin{align*}
\hat{f}(\xi,\eta)&=c_n|\eta|\(\|\xi\|^2+\eta^2\)^{\frac{n-1}{2}}\hat{g}\(\xi,\sqrt{\|\xi\|^2+\eta^2}\)\\
&=c_n (-i) \, sgn(\eta)\(\|\xi\|^2+\eta^2\)^{\frac{n-1}{2}} i \eta \hat{g}\(\xi,\sqrt{\|\xi\|^2+\eta^2}\)\mbox{.}
\end{align*}
Theorem \ref{rsternpagdachgleichetagdach} yields
\begin{align*}
\hat{f}(\xi,\eta)&=c_n (-i) \, sgn(\eta)\(\|\xi\|^2+\eta^2\)^{\frac{n-1}{2}} \widehat{\Rs_{\pa} g}(\xi,\eta)\\
&=c_n (-i) \, sgn(\eta) \(\laplace^{\frac{n-1}{2}} \Rs_{\pa} g\)\,\hat{}\,(\xi,\eta)\\
&=c_n \(H_y \laplace^{\frac{n-1}{2}} \Rs_{\pa} g\)\,\hat{}\,(\xi,\eta)\mbox{.}
\end{align*}
This completes the proof, because $f$ is in $\gS_e$.

\end{proof}

\section{Problems with Andersson's second inversion formula}

In the following it is shown that the second reformulation of the inversion formula in \cite{Andersson} is only valid in the distributional sense.

\begin{prop}

\label{kgdachinlp}

Let $f \in \gS_e$ and $g=Rf$. Then $\widehat{Kg} \in L^p(\R^n \times \R^{n+1})$ for all $p>0$.

\end{prop}

\begin{proof}

Theorem \ref{inversionformula} implies that
$$\hat{g}(\xi,\eta)=\left\{ \begin{array}{lr}
\displaystyle{(2\pi)^n \frac{2}{|S^n|} \frac{\hat{f}\(\xi,\sqrt{\|\eta\|^2-\|\xi\|^2}\)}{\|\eta\|^{n-1} \sqrt{\|\eta\|^2-\|\xi\|^2}}} &
\mbox{for }\|\eta\|>\|\xi\| \\
0 & \mbox{otherwise.} \\
\end{array} \right.$$
Hence
\begin{align*}
&\left| \widehat{Kg}(\xi,\eta)\right|=\left\{ \begin{array}{lr}
\displaystyle{(2\pi)^n \frac{2}{|S^n|} \hat{f}\(\xi,\sqrt{\|\eta\|^2-\|\xi\|^2}\)} & \mbox{for }\|\eta\|>\|\xi\| \\
0 & \mbox{otherwise} \\
\end{array} \right.\\
&\leq \left\{ \begin{array}{lr}
\displaystyle{(2\pi)^n \frac{C}{|S^n|} \(1+\sqrt{\|\xi\|^2+|\|\eta\|^2-\|\xi\|^2|}\)^{-\frac{2n+2}{p}}} & \mbox{for }\|\eta\|>\|\xi\| \\
0 & \mbox{otherwise} \\
\end{array} \right.\\
&\leq \left\{ \begin{array}{lr}
\displaystyle{(2\pi)^n \frac{C}{|S^n|} (1+\|\eta\|)^{-\frac{2n+2}{p}}} & \mbox{for }\|\eta\|>\|\xi\| \\
0 & \mbox{otherwise} \\
\end{array} \right.\\
&\leq (2\pi)^n \frac{C}{|S^n|} (1+\mbox{max}(\|\xi\|,\|\eta\|))^{-\frac{2n+2}{p}}
\end{align*}
with an appropriate $C > 0$.

\end{proof}

This proposition is sufficient for the validity of the second reconstruction formula in the distributional sense.However it is noteworthy that in general $g$ is only in $C^\infty$ and not for example in $L^2$. Therefore $\hat{g}$ has to be computed in the distributional sense and no further improvement for this inversion formula is achievable.

\begin{cor}

Let $f \in \D_e$ with $\hat{f}(0,0) \not= 0$ and $g=Rf$. Then $g,\hat{g} \not\in L^2(\R^n \times \R^{n+1})$.

\end{cor}

\begin{proof}

Theorem \ref{inversionformula} yields
$$\hat{g}(\xi,\eta)=\left\{ \begin{array}{lr}
\displaystyle{(2\pi)^n \frac{2}{|S^n|} \frac{\hat{f}\(\xi,\sqrt{\|\eta\|^2-\|\xi\|^2}\)}{\|\eta\|^{n-1} \sqrt{\|\eta\|^2-\|\xi\|^2}}} &
\mbox{for }\|\eta\|>\|\xi\| \\
0 & \mbox{otherwise.} \\
\end{array} \right.$$
Therefore
$$\uint_{\R^n \times \R^{n+1}} |\hat{g}(\xi,\eta)|^2\,d\xi \, d\eta=(2\pi)^{2n}\frac{4}{|S^n|^2}\uint_{\|\eta\| \geq \|\xi\|}\frac{\left|\hat{f}\(\xi,\sqrt{\|\eta\|^2-\|\xi\|^2}\)\right|^2}{\|\eta\|^{2n-2} \(\|\eta\|^2-\|\xi\|^2\)}\,d\xi \, d\eta \mbox{.}$$
The substitution $\|\eta\|=\rho'$ results in
\begin{align*}
\uint_{\R^n \times \R^{n+1}}& |\hat{g}(\xi,\eta)|^2\,d\xi \, d\eta\\
&=(2\pi)^{2n}\frac{4}{|S^n|}\uint_{\R^{n}}\uint_{\rho' \geq \|\xi\|}\frac{\rho'^{2-n}\left|\hat{f}\(\xi,\sqrt{\rho'^2-\|\xi\|^2}\)\right|^2}{\rho'^2-\|\xi\|^2}\,d\xi \, d\rho' \mbox{.}
\end{align*}
The substitution $\rho'=\rho + \|\xi\|$ leads to
\begin{align*}
\uint_{\R^n \times \R^{n+1}}& |\hat{g}(\xi,\eta)|^2\,d\xi \, d\eta\\
&=(2\pi)^{2n}\frac{4}{|S^n|}\uint_{\R^n}\uint_{\rho \geq 0}\frac{(\rho + \|\xi\|)^{2-n} \left|\hat{f}\(\xi,\sqrt{\rho^2+2\rho \|\xi\|}\)\right|^2}{\rho^2+2\rho \|\xi\|}\,d\xi \, d\rho \mbox{.}
\end{align*}
We consider
$$\uint_{\|\xi\| \leq \frac{1}{2}}\uint_{0 \leq \rho \leq \frac{1}{2}}\frac{(\rho + \|\xi\|)^{2-n} \left|\hat{f}\(\xi,\sqrt{\rho^2+2\rho \|\xi\|}\)\right|^2}{\rho^2+2\rho \|\xi\|}\,d\xi \, d\rho$$
and assume without loss of generality $\left|\hat{f}\(\xi,\sqrt{\rho^2+2\rho\|\xi\|}\)\right| \geq 1$ for $\|\xi\| \leq \frac{1}{2}$, \mbox{$0 \leq \rho \leq \frac{1}{2}$:}
$$\uint_{\|\xi\| \leq \frac{1}{2}}\uint_{0 \leq \rho \leq \frac{1}{2}}\frac{(\rho + \|\xi\|)^{2-n} \left|\hat{f}\(\xi,\sqrt{\rho^2+2\rho \|\xi\|}\)\right|^2}{\rho^2+2\rho \|\xi\|}\,d\xi \, d\rho$$
$$\geq \uint_{\|\xi\| \leq \frac{1}{2}}\uint_{0 \leq \rho \leq \frac{1}{2}}\frac{\(\rho + \|\xi\|\)^{2-n} }{\rho^2+2\rho \|\xi\|}\,d\xi \, d\rho \mbox{.}$$
With the substitution $r=\|\xi\|$ we obtain
$$\uint_{\|\xi\| \leq \frac{1}{2}}\uint_{0 \leq \rho \leq \frac{1}{2}}\frac{(\rho + \|\xi\|)^{2-n} \left|\hat{f}\(\xi,\sqrt{\rho^2+2\rho \|\xi\|}\)\right|^2}{\rho^2+2\rho \|\xi\|}\,d\xi \, d\rho$$
$$\geq |S^{n-1}|\uint_{0 \leq r \leq \frac{1}{2}}\uint_{0 \leq \rho \leq \frac{1}{2}}\frac{(\rho + r)^{2-n} r^{n-1}}{\rho^2+2\rho r}\,dr \, d\rho \mbox{.}$$
For $n=1$ this is
\begin{align*}
\uint_{\|\xi\| \leq \frac{1}{2}}&\uint_{0 \leq \rho \leq \frac{1}{2}}\frac{(\rho + \|\xi\|) \left|\hat{f}\(\xi,\sqrt{\rho^2+2\rho \|\xi\|}\)\right|^2}{\rho^2+2\rho \|\xi\|}\,d\xi \, d\rho\\
&\geq \uint_{0 \leq r \leq \frac{1}{2}}\uint_{0 \leq \rho \leq \frac{1}{2}}\frac{\rho + r}{\rho^2+2\rho r}\,dr \, d\rho
\geq \uint_{0 \leq r \leq \frac{1}{2}}\uint_{0 \leq \rho \leq \frac{1}{2}}\frac{1}{2\rho}\,dr \, d\rho=\infty \mbox{.}
\end{align*}
For $n \geq 2$ this results in
$$\uint_{\|\xi\| \leq \frac{1}{2}}\uint_{0 \leq \rho \leq \frac{1}{2}}\frac{(\rho + \|\xi\|)^{2-n} \left|\hat{f}\(\xi,\sqrt{\rho^2+2\rho \|\xi\|}\)\right|^2}{\rho^2+2\rho \|\xi\|}\,d\xi \, d\rho$$
$$\geq \uint_{0 \leq r \leq \frac{1}{2}}\uint_{0 \leq \rho \leq \frac{1}{2}}\frac{r^{n-1}}{\rho^2+2\rho r}\,dr \, d\rho=\infty \mbox{.}$$

\end{proof}

\begin{rem}

Unfortunately usually $\hat{f}(0,0) \not= 0$ for a reflectivity function $f \not\equiv 0$ with physically realistic properties. Moreover $\hat{g}$ is not continuous, and therefore $g \notin L^1$, so it might be difficult to compute $\hat{g}$ numerically with sufficient accuracy. Therefore the other inversion formula seems to be a better approach numerically.

\end{rem}

\section{Numerical simulations}

As can be seen in theorem \ref{modinvform}, the exact reconstruction of the reflectivity function $f$ requires data from the whole half plane. This is impossible in practice, therefore several ways to handle this problem were developed, ranging from a simple cutoff to applying an exponential decay towards the edges of the data. In the following a new approach is proposed to lessen the artifacts caused by this limitation of the data. Then the results of computer simulations using this new approach are compared to the results using a simple cutoff. The problem of limited data in the spherical Radon transform is thoroughly examined in chapter~\ref{Fehler}.

\subsection{New approach\label{approx}}

The proof of proposition \ref{rsternpawohldef} yields
\begin{align*}
&(\Rs_{\pa} g)(x,y)=\uint_{\R^n} \frac{1}{|S^n|} \uint_{S^n} \frac{y}{\sqrt{\|x-z\|^2 + y^2}}\\
&\quad\quad\times {\xi \choose \eta} \cdot (\nabla f)\(z+\xi\sqrt{\|x-z\|^2 + y^2},\eta\sqrt{\|x-z\|^2 + y^2}\)\, dS_n(\xi,\eta) \, dz\\
&=\uint_{\R^n} \frac{1}{|S^n|} \frac{y}{\(\sqrt{\|x-z\|^2 + y^2}\)^{n+1}}\\
&\quad\quad\times \uint_{\|r\|=\sqrt{\|x-z\|^2 + y^2}} \frac{r}{\|r\|} \cdot (\nabla f)((z,0)+r)\, d\sigma(r) \, dz\mbox{.}
\end{align*}
The new approach presented here is based upon the idea that the integral
$$\uint_{\|r\|=\sqrt{\|x-z\|^2 + y^2}} \frac{r}{\|r\|} \cdot (\nabla f)((z,0)+r)\, d\sigma(r)$$
does not vary much for large values of $z$ because for large $z$ the integral describes a circle with a large radius that runs through the support of $f$. Therefore the curvature only changes slightly and since $f \in \D$ the same should hold for the integral. The missing data is replaced by an approximation that uses the first and last known data with regard to the variable of integration $z$ as an approximation for the interval of integration where the data is unknown. Since it is easily seen from the formula above that $\frac{\pa}{\pa y} g(x,y)=0$ for $y=0$, it is sufficient that the following approximation is well defined for $y \neq 0$, where $z_{min}$ and $z_{max}$ denote the smallest and the largest values of $z$ with data available:
\begin{align*}
(&\Rs_{\pa} g)(x,y)  \approx \uint_{z_{min}}^{z_{max}}  \frac{\pa}{\pa y} g\(z,\sqrt{\|x-z\|^2 + y^2}\) \,dz\\
& \quad\quad +\uint_{-\infty}^{z_{min}} \frac{1}{|S^n|} \frac{y}{\(\sqrt{\|x-z\|^2 + y^2}\)^{n+1}}\\
& \quad\quad \times \frac{\(\sqrt{\|x-z_{min}\|^2 + y^2}\)^{n+1}}{y} \frac{y}{\(\sqrt{\|x-z_{min}\|^2 + y^2}\)^{n+1}}\\
& \quad\quad \times \uint_{\|r\|=\sqrt{\|x-z_{min}\|^2 + y^2}} \frac{r}{\|r\|} \cdot (\nabla f)((z_{min},0)+r)\, d\sigma(r) \, dz\\
& \quad\quad +\uint_{z_{max}}^\infty \frac{1}{|S^n|} \frac{y}{\(\sqrt{\|x-z\|^2 + y^2}\)^{n+1}}\\
& \quad\quad \times \frac{\(\sqrt{\|x-z_{max}\|^2 + y^2}\)^{n+1}}{y} \frac{y}{\(\sqrt{\|x-z_{max}\|^2 + y^2}\)^{n+1}}\\
& \quad\quad \times \uint_{\|r\|=\sqrt{\|x-z_{max}\|^2 + y^2}} \frac{r}{\|r\|} \cdot (\nabla f)((z_{max},0)+r)\, d\sigma(r) \, dz\\
& =\uint_{z_{min}}^{z_{max}} \frac{\pa}{\pa y} g\(z,\sqrt{\|x-z\|^2 + y^2}\) \,dz&\\
& \quad\quad +\uint_{-\infty}^{z_{min}} \frac{y}{\(\sqrt{\|x-z\|^2 + y^2}\)^{n+1}} \frac{\(\sqrt{\|x-z_{min}\|^2 + y^2}\)^{n+1}}{y}\\
& \quad\quad \times \frac{\pa}{\pa y} g\(z_{min},\sqrt{\|x-z_{min}\|^2 + y^2}\) \,dz\\
& \quad\quad +\uint_{z_{max}}^\infty \frac{y}{\(\sqrt{\|x-z\|^2 + y^2}\)^{n+1}} \frac{\(\sqrt{\|x-z_{max}\|^2 + y^2}\)^{n+1}}{y}\\
& \quad\quad \times \frac{\pa}{\pa y} g\(z_{max},\sqrt{\|x-z_{max}\|^2 + y^2}\) \,dz
\end{align*}
\begin{align*}
& =\uint_{z_{min}}^{z_{max}} \frac{\pa}{\pa y} g\(z,\sqrt{\|x-z\|^2 + y^2}\) \,dz&\\
& \quad\quad +\left(\frac{\pi}{2}-\arctan\left(\frac{x-z_{min}}{y}\right)\right) \frac{\(\sqrt{\|x-z_{min}\|^2 + y^2}\)^{n+1}}{y}\\
& \quad\quad \times \frac{\pa}{\pa y} g\(z_{min},\sqrt{\|x-z_{min}\|^2 + y^2}\)\\
& \quad\quad +\left(\frac{\pi}{2}+\arctan\left(\frac{x-z_{max}}{y}\right)\right) \frac{\(\sqrt{\|x-z_{max}\|^2 + y^2}\)^{n+1}}{y}\\
& \quad\quad \times \frac{\pa}{\pa y} g\(z_{max},\sqrt{\|x-z_{max}\|^2 + y^2}\)
\end{align*}

\subsection{Comparison}

Now the images obtained via the approximate reconstruction formula derived above are compared to the results of a simple reconstruction formula that sets unmeasured data to $0$. This comparison will highlight the advantages of this new approach.

\dkbild{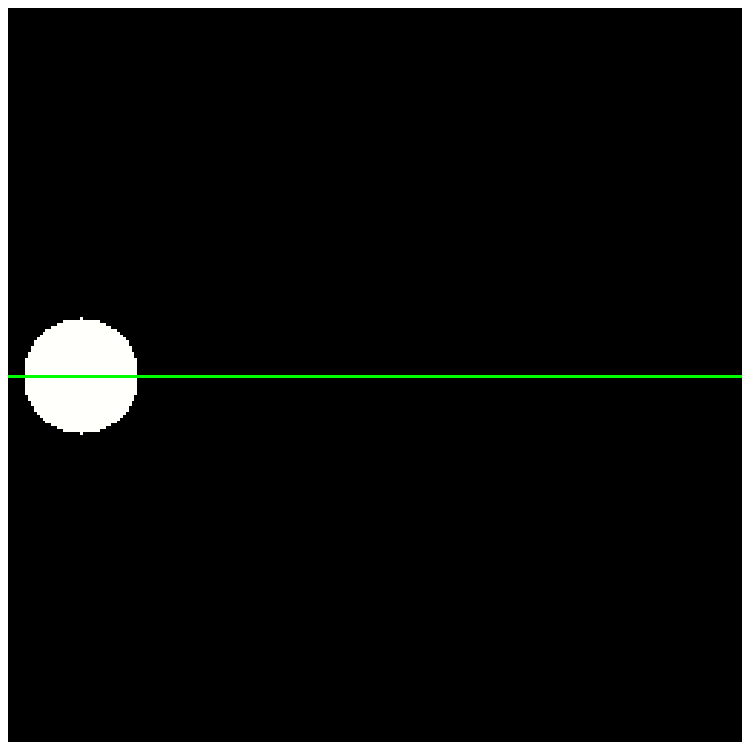}{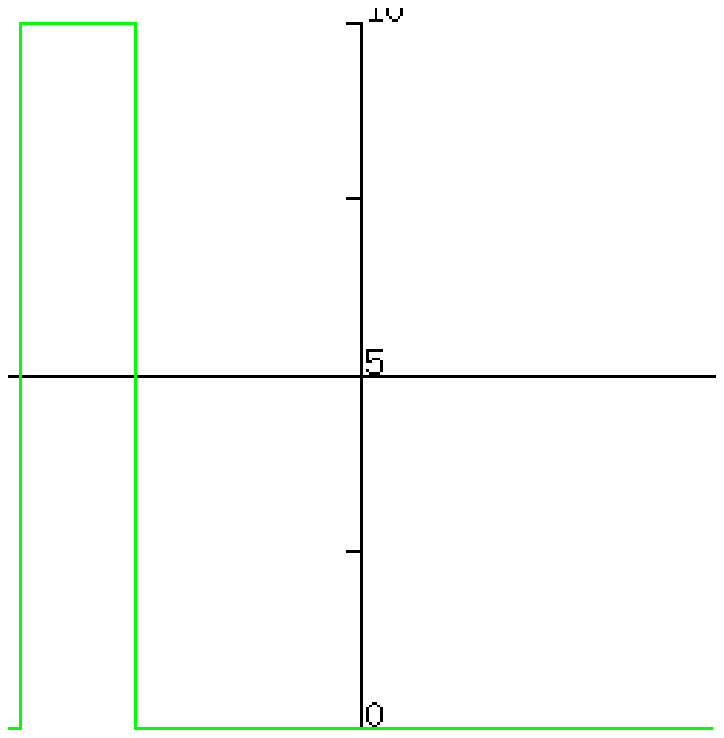}{Phantom\label{phantom}}

The underlying phantom for the following comparison is very simple. It is a circle with a radius of $20$ as depicted on the left in figure \ref{phantom}. Its center is with $(0,25)$ close to the flighttrack that runs along the left edge $(y=0)$. The circle has an amplitude of $10$, the remaining part has a reflectivity of $0$, as can be seen in the cross section on the right that is indicated in the left image by the horizontal line.
The figures shown in the following are all composed in this same way. The simulated measurements are noise free. Several reconstruction pairs are shown. The first uses data that stretches exactly as far as the reconstructed area, i. e. $0 \leq r < 256$ and $-128 \leq x < 128$. The following pairs are reconstructed each with a larger amount of data in both directions with respect to the preceding one. The reconstructions are designated accordingly. The first reconstruction in each pair is obtained by continuing the data with $0$ in the region where data is missing. The second is reconstructed using the approximation derived above.

Figures \ref{alt1} and \ref{neu1} show reconstructions that use an amount of data that is exactly as large as the data contained in the images, i. e. $0 \leq r < 256$ and $-128 \leq x < 128$. Data that is unavailable is set to $0$ for reconstruction purposes in figure \ref{alt1}. This causes two broad circular artifacts - one curved upward and one downward - that are clearly visible. The formation of the artifacts can be understood as follows. The information in the data that is actually measured causes the algorithm to reconstruct the large positive circular phantom. The only possibility however to conform to the $0$ of the continuation of the data is to form these circular negative artifacts. This is the only way to achieve that the circular integral which runs through the large positive phantom becomes $0$. Therefore the radius of these artifacts matches exactly the smallest radius that is missing in the data. A side effect is that both algorithms overshoot the amplitude of the circle. It is notable that in figure \ref{alt1} the error is large close to the object and decreases as the artifact closes to the edges of the image. Another problem is the noticeable gradient in the object's amplitude.\\
\dkdbildp{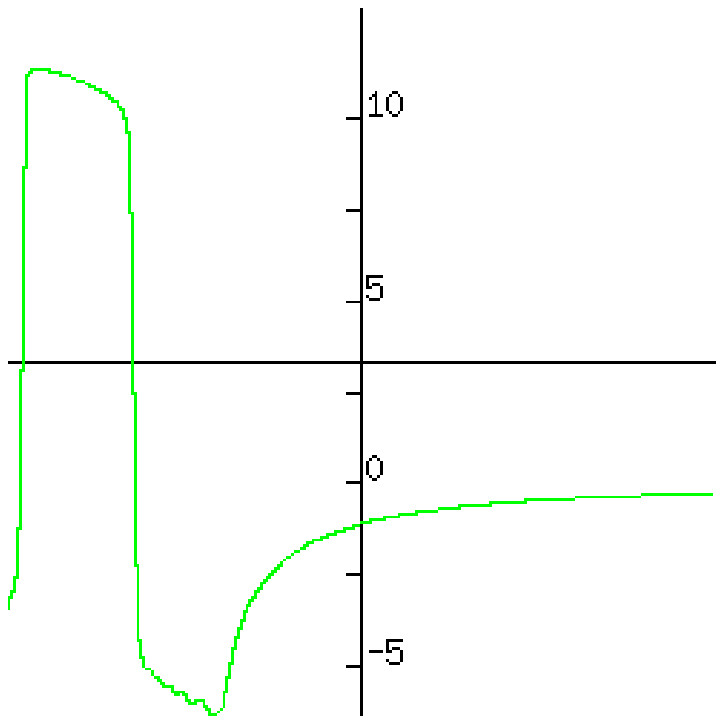}{Simple reconstruction; single data\label{alt1}}{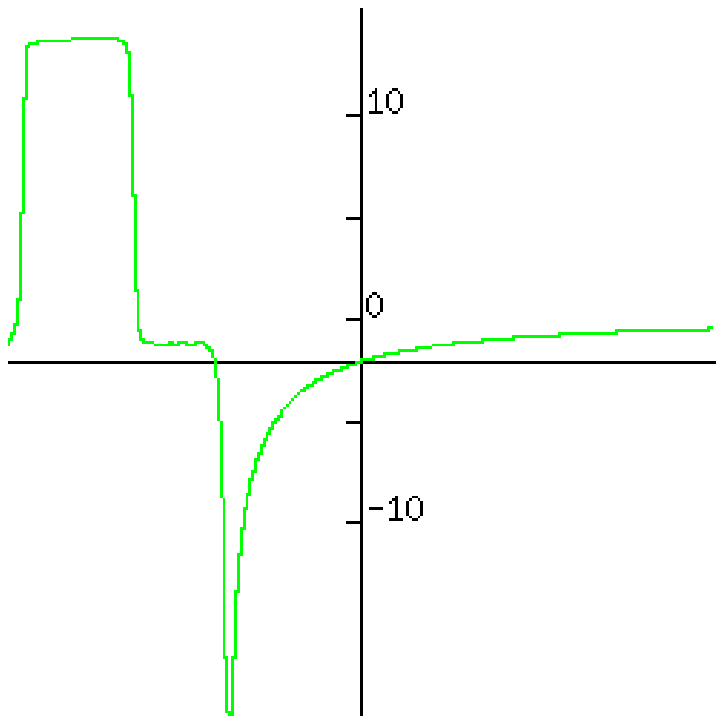}{Reconstruction using approximate continuation; single data\label{neu1}}
The missing data was dealt with as delineated in subsection \ref{approx} to obtain figure \ref{neu1}. Again this causes artifacts that are similar to the artifacts in figure \ref{alt1}. The reconstruction is very good close to the object, as can be seen in the cross section before the dip, but the artifacts get worse as they approach the edges of the image. This comes from the fact that for large values of $|x|$ the approximation gets less reliable. However, the gradient in the object's amplitude in figure \ref{alt1} is not reflected in figure \ref{neu1} where the amplitude is constant, as it should be.

\dkdbild{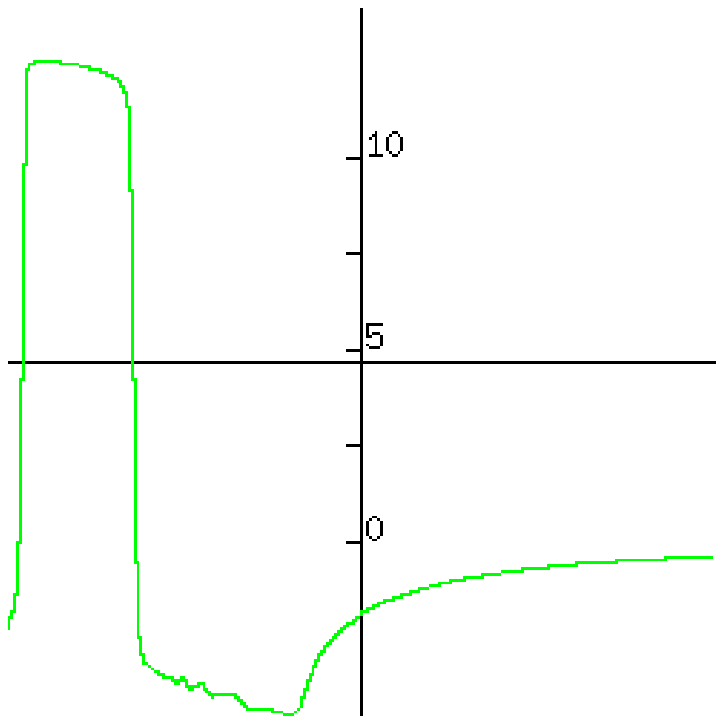}{Simple reconstruction; quadruple data\label{alt2}}{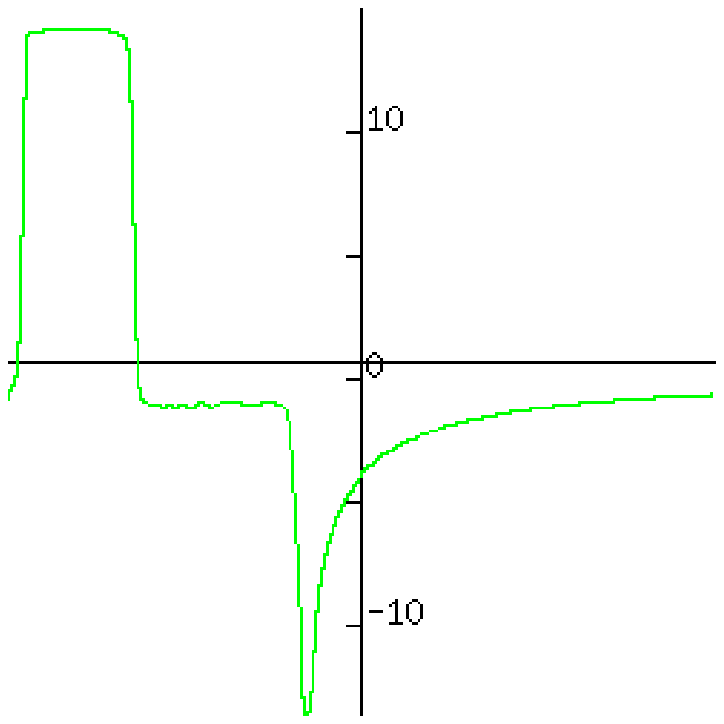}{Reconstruction using approximate continuation; quadruple data\label{neu2}}



In figures \ref{alt2} and \ref{neu2} double the amount of data in each direction is used relative to the preceding two images, i. e. $0 \leq r < 512$ and $-256 \leq x < 256$. Again the first figure shows the reconstruction using a continuation by $0$ for missing data, whereas the second figure is computed with the help of the approximation.
The artifacts seen in figures \ref{alt2} and \ref{neu2} are similar to the artifacts seen in the previous image pair, albeit less severe. The gradient of the circle in figure \ref{alt2} is less steep than in figure \ref{alt1} and the dip below $0$ is shallower but broader. A similar effect is noticeable in comparison of figures \ref{neu1} and \ref{neu2}. The area close to the object with its very good accuracy between the object and the dip, where the reflectivity stays constant around $0$, is larger and the dip is more shallow.



\dkdbildt{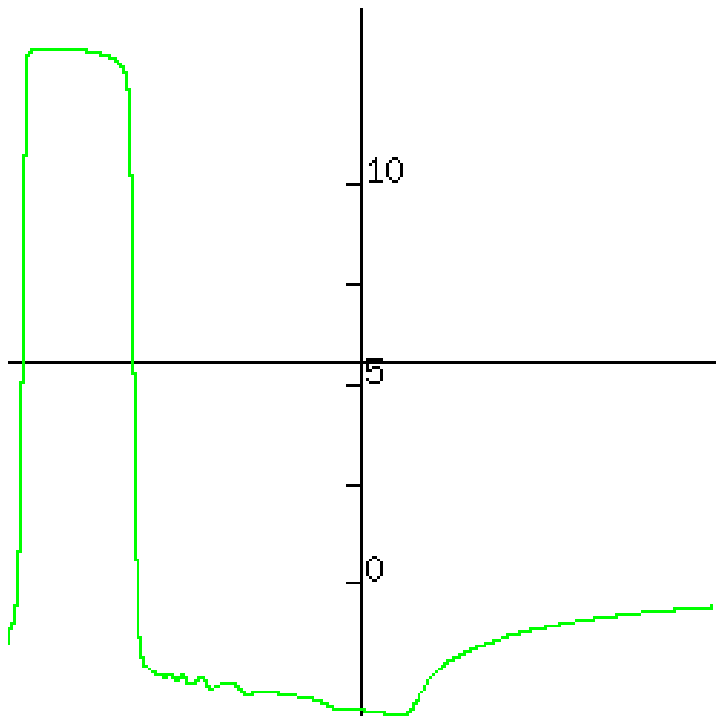}{Simple reconstruction; $4^2$ fold data\label{alt4}}{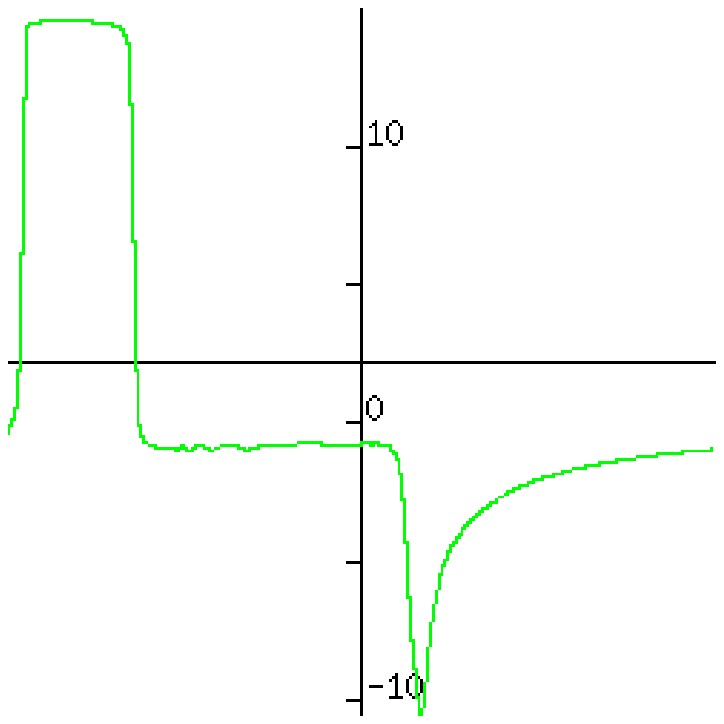}{Reconstruction using approximate continuation; $4^2$ fold data\label{neu4}}

Figures \ref{alt4} and \ref{neu4} are computed with quadruple the amount of data relative to the former pair, i. e. $0 \leq r < 1024$ and $-512 \leq x < 512$.
The trend of the preceding comparison continues in these figures. The gradient in the circle's amplitude is almost negligible in figure \ref{alt4} and the negative region behind the object is again shallower and broader. Also in figure \ref{neu4} the dip is not as deep as in figure \ref{neu2} and the area behind the circle that has an amplitude close to $0$ is broader.





\dkdbild{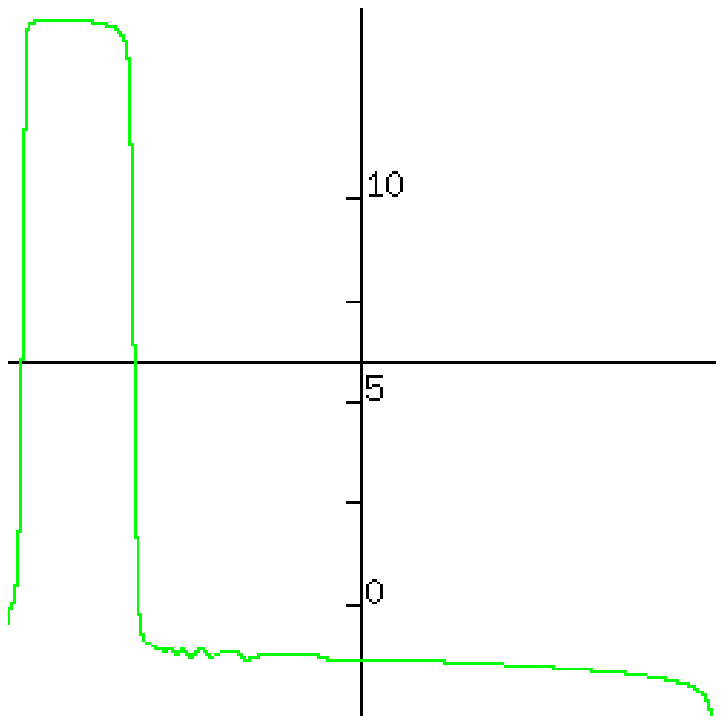}{Simple reconstruction; $16^2$ fold data\label{alt16}}{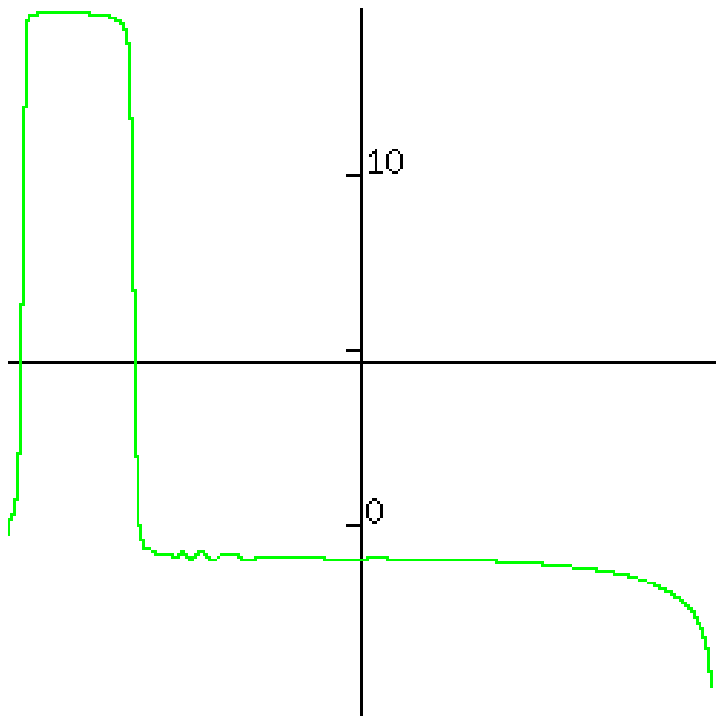}{Reconstruction using approximate continuation; $16^2$ fold data\label{neu16}}


The data used to reconstruct figures \ref{alt16} and \ref{neu16} encompasses 16 fold the length in each direction with respect to the reconstructed image, i. e. $0 \leq r < 4096$ and $-2048 \leq x < 2048$.
The gradient in the amplitude of the object in figure \ref{alt16} is nonexistent. An interesting difference between figure \ref{alt16} and the preceding images of that kind can be seen in the negative region behind the circle. The area has broadened, but at the right edge of the image there is a dip in contrast to the preceding constructions of this kind. This reminds of the reconstructions of the other type. Figure \ref{neu16} shows the same pattern as the previous reconstructions using the approximation formula. It is worth mentioning that the area where the amplitude is close to $0$ almost stretches over the whole image.


\dkdbild{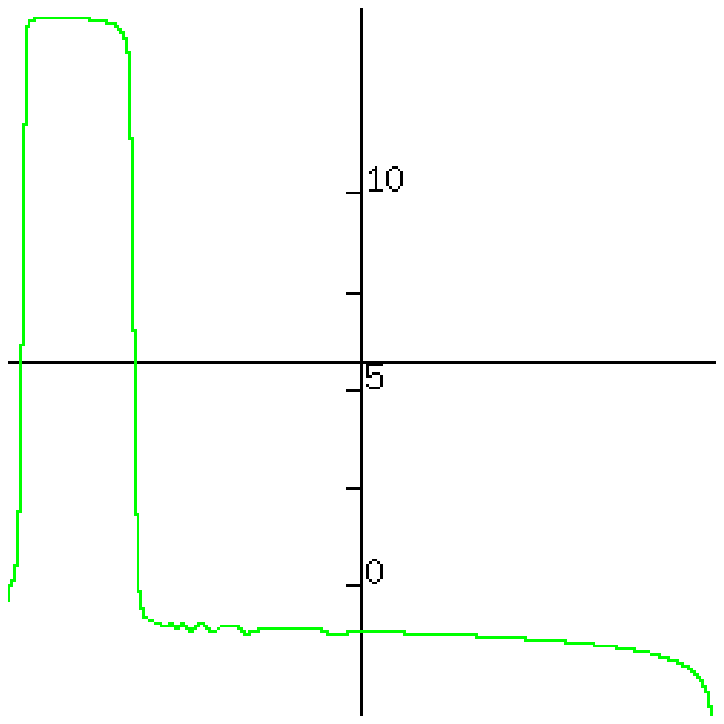}{Simple reconstruction; $32^2$ fold data\label{alt32}}{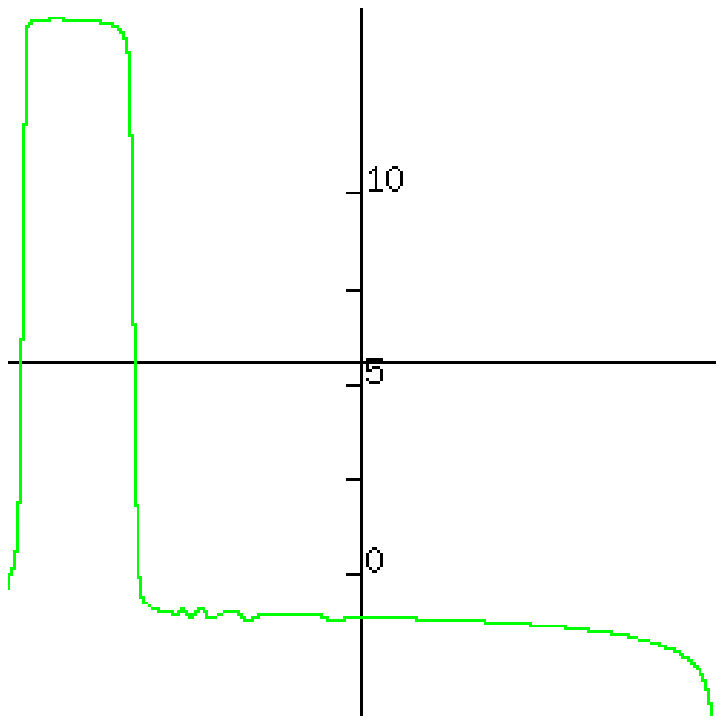}{Reconstruction using approximate continuation; $32^2$ fold data\label{neu32}}


The reconstructions in figures \ref{alt32} and \ref{neu32} use double the amount of data in each direction relative to figures \ref{alt16} and \ref{neu16}, i. e. $0 \leq r < 8192$ and $-4096 \leq x < 4096$. For this large amount of data both reconstructions look very much alike. This is due to the fact that the region where data is missing is far away and therefore the data from that region only has a very small influence on the reconstruction.\\


The comparison of the preceding images suggests that using the previously derived approximative continuation of the data is superior to a continuation by $0$ since the reconstruction quality close to the flight track is much better. As is well known from experience \cite{Klein}, the reconstruction quality away from the flight track is plagued by artifacts and therefore unreliable  anyway. Therefore the method of approximate continuation seems to be preferable. In addition the adverse effects away from the flight track could probably be alleviated by using an exponential decay as a kind of mollifier in the approximation algorithm.

\chapter{Ghosts due to limited data}

\label{Fehler}

Depending on the reconstruction algorithm there are various artifacts that appear in images obtained by the inversion of the spherical Radon transform, as seen in chapter \ref{Diplomkorrektur}. Since theorem \ref{inversionformula} states that $R$ is one-to-one, there should not be any artifacts in these reconstructions. However theorem \ref{inversionformula} requires the data to be known for all $x$ and all $r>0$. As this is not the case, uniqueness is lost. Therefore it seems obvious that one of the main origins of these artifacts is the impossibility of measuring infinitely far. This phenomenon is further examined in this chapter. It will be shown that there exist functions with support inside the measured region that do not have any effect on the measured data. Functions like these are commonly called ghosts. A complete descriptions of these ghosts will be given and some examples will be shown.\\

\section{The problem of limited data}

The reconstruction formulas in theorem \ref{inversionformula}, corollary \ref{reformulation}, and theorem \ref{modinvform} are exact within the respective assumptions. Therefore an exact reconstruction without artifacts should be possible. But for each point which is to be reconstructed all the reconstruction formulas require data from arbitrary long distances with arbitrary large radii. Of course this is impossible in reality. Even in computer simulations this can not be achieved.\\
Problems like this are known from various inverse problems. There are two kinds of difficulties associated with the implementation of the inversion formula for the spherical Radon transform. First, the data is only gathered in a discretized manner and can only be handled in a discretized manner. And secondly, the data can only be supported on a compact interval in opposition to the analytical model. For both problems there is an analogy to computerized tomography. It is known for computerized tomography that the artifacts arising from the discretization have high frequencies and do not impair reconstruction quality if properly handled \cite{Louis1}, \cite{Louis2}. Therefore this problem is not considered here, as it could probably be amended in a similar way and just as in the case of computerized tomography would probably not be as severe as the problem of limited data. The problems that are discussed in this chapter emerge solely from the fact that data is only collected over compact intervals in plane positions and circle radii. Similar problems in computerized tomography are the limited angle problem that is known to produce severe artifacts \cite{Louis3} and the exterior problem, respectively. There are differences however, as the exterior problem still preserves the uniqueness in computerized tomography, although it causes some instability \cite{Quinto}. It will be shown in the following that both limitations cause artifacts in the reconstruction of SAR-data.\\
Only the two dimensional case is considered here, but the results should be easily transferable to higher dimensions.

For a decent analysis of the effects of limited data it is crucial to find a set of orthogonal functions that represents the whole data space outside the measurable region because such a set allows to restrict the analysis of the effects to a well known set of functions. The problem is to find such a set that can also be analytically inverted by an inversion formula for the spherical Radon transform. To this end the following steps are necessary.

First, a set of orthogonal functions that are supported outside of the measurable region is constructed and some important properties are listed. It is shown that for an arbitrary, hypothetical measurement that extends over the whole half plane the information gathered by the projections of the data onto the set of orthogonal functions is sufficient to regain the data outside of some compact set that can realistically be measured.  Then the orthogonal functions, which span the whole space of data functions that could be measured outside this compact set, are inverted using the inversion formula in theorem \ref{inversionformula}. Finally some numerical examples are shown.

These examples show that the problem of limited data should not be neglected. Since the artifacts that can arise are severe and depend on the chosen reconstruction formula, as seen in the preceding chapter, it should be tried to find means to minimize these artifacts.


\section{Orthogonal functions and their transforms}

In the following, the maximal radius up to which the data is measured will be denoted with $R$ and the start- and endpoints of the flight track will be denoted with $-L$ and $L$ respectively. Therefore the data function $g(x,r)$ is only known for $0 \leq r \leq R$ and $-L \leq x \leq L$. Now some functions and their transforms are compiled. This information is later used to construct a set of orthogonal functions supported in $\R \times [0,\infty) \setminus [-L,L] \times [0,R]$.

\begin{defi}

$\bn$ denotes the Bessel function of order zero, which is defined as
$$\bn(z):=\sum_{k=0}^\infty (-1)^k \frac{z^{2k}}{2^{2k}k!\Gamma(k+1)}$$
for $|\mbox{arg } z|<\pi$.

\end{defi}

\begin{defi}

For $f: [0,\infty) \rightarrow \R$ the Hankel transform of $f$ is defined as
$$\uint_0^\infty f(r) r \bn(r\rho) \, dr \mbox{.}$$

\end{defi}

\begin{lem}

\label{costrafo}

Let $b, b' > 0$. Then
$$\frac{2}{\pi} \uint_0^\infty \cos(br) \cos(b'r) \, dr=\delta(b-b') \mbox{.}$$

\end{lem}

\begin{proof}

This is easily verified using the properties of the Fourier transform \cite{Sneddon}.

\end{proof}

\begin{theo}

\label{Hankeltransform}

If $f: [0,\infty) \rightarrow \R$ and $\sqrt{r} f(r)$ is piecewise continuous and absolutely integrable, then
$$\bar{f}(\rho)=\uint_0^\infty f(r) r \bn(r\rho) \, dr$$
exists and
$$\uint_0^\infty \bar{f}(\rho) \rho \bn(\rho r) \, d\rho=\lim_{h \rightarrow 0} \frac{1}{2}(f(r+h)+f(r-h))\mbox{.}$$

\end{theo}

\begin{proof}

\cite[5-3, Theorem 1]{Sneddon}

\end{proof}

\begin{lem}

\label{cosinusunendlich}

Let $a>0$ and
$$f(x)=\left\{ \begin{array}{lr}
\bn\(a\sqrt{x^2-L^2}\,\) &
\mbox{ for } |x|>L\\
& \\
0 & \mbox{ for } |x|<L\mbox{.}\\
\end{array} \right.$$
Then
$$\uint_0^\infty \cos(x\xi) f(x) \,dx=\left\{ \begin{array}{lr}
\frac{\exp^{-L\sqrt{a^2-\xi^2}}}{\sqrt{a^2-\xi^2}} &
\mbox{ for }0<|\xi|<a\\
& \\
-\frac{\sin\(L\sqrt{\xi^2-a^2}\)}{\sqrt{\xi^2-a^2}} & \mbox{ for }|\xi|>a\mbox{.}\\
\end{array} \right.$$

\end{lem}

\begin{proof}

\cite[I, \S 17, p. 78]{Oberhettinger1}

\end{proof}

\begin{lem}

\label{hankelunendlich}

Let
$$f(r)=\left\{ \begin{array}{lr}
\frac{\cos\(b\sqrt{r^2-R^2}\)}{\sqrt{r^2-R^2}} &
\mbox{ for } r>R\\
& \\
0 & \mbox{ for } r<R\mbox{.}\\
\end{array} \right.$$
Then
$$\uint_0^\infty r \bn(r\rho) f(r) \, dr=\left\{ \begin{array}{lr}
\frac{\cos\(R\sqrt{\rho^2-b^2}\,\)}{\sqrt{\rho^2-b^2}} &
\mbox{ for }\rho>b\\
& \\
0 & \mbox{ for }\rho<b\mbox{.}\\
\end{array} \right.$$

\end{lem}

\begin{proof}

\cite[I, 1.2, 2.56, p. 12]{Oberhettinger2}

\end{proof}

\begin{lem}

\label{hankelendlich}

Let
$$f(r)=\left\{ \begin{array}{lr}
\frac{\cos\(b\sqrt{R^2-r^2}\)}{\sqrt{R^2-r^2}} &
\mbox{ for } r<R\\
& \\
0 & \mbox{ for } r>R\mbox{.}\\
\end{array} \right.$$
Then
$$\uint_0^\infty r \bn(r\rho) f(r) \, dr=\frac{\sin\(R\sqrt{\rho^2+b^2}\)}{\sqrt{\rho^2+b^2}} \mbox{.}$$

\end{lem}

\begin{proof}

\cite[I, 1.2, 2.55, p. 12]{Oberhettinger2}

\end{proof}

\begin{lem}

\label{hankeldelta}

Let $r, r_0 >0$. Then the Hankel transform of $r_0 \bn(r_o \cdot)$ is
$$\uint_0^\infty \rho \bn(r\rho) r_0 \bn(r_0\rho)\,d\rho=\delta(r-r_0)\mbox{.}$$
The $\cdot$ is a placeholder for the variable used in the transform.
\end{lem}

\begin{proof}

This is easily verified using theorem \ref{Hankeltransform}.

\end{proof}

Now two sets of distributions will be given, $\{o_{range}^{a,b}|a \in \R, b \geq 0\}$ and \mbox{$\{o_{even}^{a,l},o_{odd}^{a,l}|a \geq 0, l \in \N_0\}$} and it is shown that they are orthogonal. The first set consists of distributions that are supported only in $r>R$. These distributions thus represent the missing information due to the limitation that the reflected waves can only be received up to some distance $R$. The second set consists of functions that are supported only in $r<R$ and $|x|>L$. These functions therefore represent the missing information due to the fact that the plane travels only a limited distance from $-L$ to $L$. Note however that in this set only the information deficit is contained that is in addition to the first case. Therefore the constraint $r<R$ is added.

\begin{defi}

Let $R, L >0$.

\begin{enumerate}

\item

Let $a \in \R$, $b \geq 0$. Define

$$o_{range}^{a,b}(x,r):=\left\{ \begin{array}{lr}
\delta(x-a)\frac{\cos\(b\sqrt{r^2-R^2}\,\)}{\sqrt{r^2-R^2}} &
\mbox{ for } r>R\\
& \\
0 & \mbox{ otherwise.}\\
\end{array} \right.$$

\item

Let $a \geq 0$, $l \in \N_0$. Define

\begin{align*}
o_{even}^{a,l}&(x,r)\\
&:=\left\{ \begin{array}{lr}
\bn\(a\sqrt{x^2-L^2}\) \frac{\cos\(l\frac{\pi}{R}\sqrt{R^2-r^2}\,\)}{\sqrt{R^2-r^2}} &
\mbox{ for } 0 \leq r < R \mbox{ and } |x|>L\\
& \\
0 & \mbox{ otherwise.}\\
\end{array} \right.
\end{align*}

\item

Let $a \geq 0$, $l \in \N_0$. Define

\begin{align*}
o_{odd}^{a,l}&(x,r)\\
&:=\left\{ \begin{array}{lr}
x \bn\(a\sqrt{x^2-L^2}\) \frac{\cos\(l\frac{\pi}{R}\sqrt{R^2-r^2}\,\)}{\sqrt{R^2-r^2}} &
\mbox{ for } 0 \leq r < R \mbox{ and } |x|>L\\
& \\
0 & \mbox{otherwise.}\\
\end{array} \right.
\end{align*}

\end{enumerate}

\end{defi}

\begin{defi}

Define Neumann's number $\epsilon_n$ \cite{Oberhettinger3}:\\
$\epsilon_0:=1$, $\epsilon_n:=2$, $n=1, 2, 3, ...$.

\end{defi}

\begin{prop}

Let $R, L >0$.

\begin{enumerate}

\item

Let $a, a' \in \R$, $b, b' \geq 0$, $a \neq a'$, $b \neq b'$. Then the distributions $o_{range}^{a,b}(x,r)$ and $o_{range}^{a',b'}(x,r)$ are orthogonal to each other with respect to the scalar product
\mbox{$<f,g>_{o_{range}}=\uint_R^\infty \uint_\R f(x,r) g(x,r) r \sqrt{r^2-R^2} \, dx \, dr$.}

\item

Let $a, a' > 0$, $l, l' \in \N_0$, $a \neq a'$, $l \neq l'$. Then the functions $o_{even}^{a,l}(x,r)$ and $o_{even}^{a',l'}(x,r)$ are orthogonal to each other with respect to the scalar product
\mbox{$<f,g>_{o_{even}}=\uint_0^R \uint_L^\infty f(x,r) g(x,r) r \sqrt{R^2-r^2} x \, dx \, dr$.}

\item

Let $a, a' > 0$, $l, l' \in \N_0$, $a \neq a'$, $l \neq l'$. Then the functions $o_{odd}^{a,l}(x,r)$ and $o_{odd}^{a',l'}(x,r)$ are orthogonal to each other with respect to the scalar product
\mbox{$<f,g>_{o_{odd}}=\uint_0^R \uint_L^\infty f(x,r) g(x,r) r \sqrt{R^2-r^2} \frac{1}{x} \, dx \, dr$.}

\end{enumerate}

\end{prop}

\begin{proof}

Let $R, L >0$.

\begin{enumerate}

\item

Let $a, a' \in \R$, $b, b' \geq 0$, $a \neq a'$, $b \neq b'$. Then
\begin{align*}
<&o_{range}^{a,b},o_{range}^{a',b'}>_{o_{range}}\\
&=\uint_R^\infty \uint_\R o_{range}^{a,b}(x,r) o_{range}^{a',b'}(x,r) r \sqrt{r^2-R^2} \, dx \, dr\\
&=\uint_R^\infty \uint_\R \delta(a-x)\frac{\cos\(b\sqrt{r^2-R^2}\,\)}{\sqrt{r^2-R^2}}\\
&\quad\times \delta(x-a')\frac{\cos\(b'\sqrt{r^2-R^2}\,\)}{\sqrt{r^2-R^2}} r \sqrt{r^2-R^2} \, dx \, dr\\
&=(\delta_0 \ast \delta_{a'})(a) \uint_R^\infty \uint_\R \cos\(b\sqrt{r^2-R^2}\) \cos\(b'\sqrt{r^2-R^2}\) \frac{r}{\sqrt{r^2-R^2}} \, dx \, dr\mbox{.}
\end{align*}
The substitution $r'=\sqrt{r^2-R^2}$ and lemma \ref{costrafo} yield
\begin{align*}
<o_{range}^{a,b},o_{range}^{a',b'}>_{o_{range}}&=(\delta_{a'})(a) \uint_0^\infty \cos(br') \cos(b'r') \, dr'\\
&=\frac{\pi}{2} \delta(a-a') \delta(b-b')\mbox{.}
\end{align*}

\item

Let $a, a' > 0$, $l, l' \in \N_0$, $a \neq a'$, $l \neq l'$. Then
\begin{align*}
<o_{even}^{a,l}&,o_{even}^{a',l'}>_{o_{even}}\\
&=\uint_0^R \uint_L^\infty o_{even}^{a,l}(x,r) o_{even}^{a',l'}(x,r) r \sqrt{R^2-r^2} \: x \, dx \, dr\\
&=\uint_0^R \uint_L^\infty \bn\(a\sqrt{x^2-L^2}\) \frac{\cos\(l\frac{\pi}{R}\sqrt{R^2-r^2}\,\)}{\sqrt{R^2-r^2}} \bn\(a'\sqrt{x^2-L^2}\)\\
&\quad\times \frac{\cos\(l'\frac{\pi}{R}\sqrt{R^2-r^2}\,\)}{\sqrt{R^2-r^2}} r \sqrt{R^2-r^2} \: x \, dx \, dr\\
&=\uint_0^R \uint_L^\infty \bn\(a\sqrt{x^2-L^2}\) \cos\(l\frac{\pi}{R}\sqrt{R^2-r^2}\) \bn\(a'\sqrt{x^2-L^2}\)\\
&\quad\times \cos\(l'\frac{\pi}{R}\sqrt{R^2-r^2}\) \frac{r}{\sqrt{R^2-r^2}} \: x \, dx \, dr\mbox{.}
\end{align*}
The substitutions $r'=\sqrt{R^2-r^2}$ and $x'=\sqrt{x^2-L^2}$ and lemma \ref{hankeldelta} yield
\begin{align*}
<o_{even}^{a,l},o_{even}^{a',l'}>_{o_{even}}&=\uint_0^R \uint_0^\infty \bn(ax') \cos(l\frac{\pi}{R}r') \bn(a'x') \cos(l'\frac{\pi}{R}r') x' \, dx' \, dr'\\
&=\frac{1}{a}\delta(a-a') \frac{R}{\epsilon_l}\delta_{ll'}
\end{align*}
with the Neumann's number $\epsilon_l$.

\item

Let $a, a' > 0$, $l, l' \in \N_0$, $a \neq a'$, $l \neq l'$. Then
\begin{align*}
<o_{odd}^{a,l}&,o_{odd}^{a',l'}>_{o_{odd}}\\
&=\uint_0^R \uint_L^\infty o_{odd}^{a,l}(x,r) o_{odd}^{a',l'}(x,r) r \sqrt{R^2-r^2} \: \frac{1}{x} \, dx \, dr\\
&=\uint_0^R \uint_L^\infty x \bn\(a\sqrt{x^2-L^2}\) \frac{\cos\(l\frac{\pi}{R}\sqrt{R^2-r^2}\,\)}{\sqrt{R^2-r^2}} x \bn\(a'\sqrt{x^2-L^2}\)\\
&\quad\times \frac{\cos\(l'\frac{\pi}{R}\sqrt{R^2-r^2}\,\)}{\sqrt{R^2-r^2}} r \sqrt{R^2-r^2} \: \frac{1}{x} \, dx \, dr\\
&=\uint_0^R \uint_L^\infty \bn\(a\sqrt{x^2-L^2}\) \cos\(l\frac{\pi}{R}\sqrt{R^2-r^2}\) \bn\(a'\sqrt{x^2-L^2}\)\\
&\quad\times \cos\(l'\frac{\pi}{R}\sqrt{R^2-r^2}\) \frac{r}{\sqrt{R^2-r^2}} \: x \, dx \, dr\mbox{.}
\end{align*}
The substitutions $r'=\sqrt{R^2-r^2}$ and $x'=\sqrt{x^2-L^2}$ and lemma \ref{hankeldelta} yield
\begin{align*}
<o_{odd}^{a,l},o_{odd}^{a',l'}>_{o_{odd}}&=\uint_0^R \uint_0^\infty \bn(ax') \cos(l\frac{\pi}{R}r') \bn(a'x') \cos(l'\frac{\pi}{R}r') x' \, dx' \, dr'\\
&=\frac{1}{a}\delta(a-a') \frac{R}{\epsilon_l}\delta_{ll'}
\end{align*}
with the Neumann's number $\epsilon_l$.

\end{enumerate}

\end{proof}

\section{Projection of unmeasurable data onto orthogonal functions}

In the following the projections of the data onto the orthogonal function sets are defined, and it is shown that the data for $r>R$ and $x>L$ can be recovered from these projections.

\begin{lem}

Let $f \in \gS$ and $g=Rf$. Then $g(x,r)=g(x,-r)$.

\end{lem}

\begin{proof}

This is easily verified using the definition of $Rf(x,r)$.

\end{proof}

%
%
%
%
%
%

\begin{defi}

Let $f \in \gS$ and $g=Rf$. Define

\begin{enumerate}

\item

$$G_{range}(a,b):=\frac{2}{\pi} <o_{range}^{a,b}(x,r),g(x,r)>_{o_{range}}$$

\item

$$G_{even}^l(a):=\frac{\epsilon_l}{2R} a <o_{even}^{a,l}(x,r),g(x,r)+g(-x,r)>_{o_{even}}$$

\item

$$G_{odd}^l(a):=\frac{\epsilon_l}{2R} a <o_{odd}^{a,l}(x,r),g(x,r)-g(-x,r)>_{o_{odd}}$$

\end{enumerate}

\end{defi}

\begin{defi}

Let $a,b \in \R \cup \{-\infty,\infty\}$, $a<b$, and $x \in \R$. Then define the characteristic function of $(a,b)$
$$\chi_{(a,b)}(x):=\left\{ \begin{array}{lr}
1 &
\mbox{ for } x \in (a,b)\\
& \\
0 & \mbox{ otherwise.}\\
\end{array} \right.$$

\end{defi}

\begin{theo}

\begin{enumerate}

\label{retrieve}

\item

$G_{range}$ is well defined if $f \in \gS$ and $g=Rf$. Moreover, for $r > R$
$$g(x,r)=\uint_0^\infty\uint_\R o_{range}^{a,b}(x,r) G_{range}(a,b) \,da\,db\mbox{.}$$

\item

$G_{even}^l(a)$ is well defined if $f \in \gS$ and $g=Rf$. Moreover, if $g$ is even in $x$, $0 < r < R$, and $x > L$, then
$$g(x,r)=\sum_{l=0}^\infty \uint_0^\infty o_{even}^{a,l}(x,r) G_{even}^l(a)\,da\mbox{.}$$

\item

$G_{odd}^l(a)$ is well defined if $f \in \gS$ and $g=Rf$. Moreover, if $g$ is odd in $x$, $0 < r < R$, and $x > L$, then
$$g(x,r)=\sum_{l=0}^\infty \uint_0^\infty o_{odd}^{a,l}(x,r) G_{odd}^l(a)\,da\mbox{.}$$

\end{enumerate}

\end{theo}

\begin{proof}

Let $f \in \gS$ and $g=Rf$.

\begin{enumerate}

\item

\begin{enumerate}

\item

First the well-definedness of $G_{range}$ will be shown.

\begin{align*}
|S^1| & \frac{\pi}{2} \left|G_{range}(a,b)\right|\\
&=|S^1| \left|<o_{range}^{a,b}(x,r),g(x,r)>_{o_{range}}\right|\\
&=|S^1|\left|\uint_R^\infty \uint_\R \delta(x-a)\frac{\cos\(b\sqrt{r^2-R^2}\,\)}{\sqrt{r^2-R^2}} g(x,r) r \sqrt{r^2-R^2} \, dx \, dr\right|\\
&=|S^1|\left|\uint_\R\uint_R^\infty \delta(a-x) r \cos\(b\sqrt{r^2-R^2}\)g(x,r)\,dr \, dx\right|\\
&=|S^1|\left|\uint_R^\infty r \cos\(b\sqrt{r^2-R^2}\)g(a,r)\,dr \, dx\right|\\
&=\left|\uint_R^\infty r \cos\(b\sqrt{r^2-R^2}\)\uint_{S^1}f(a+r\xi,r\eta) \, d S_1(\xi,\eta) \, dr\right|\mbox{.}
\end{align*}
With $\binom{r\xi}{r\eta}=y=\binom{y_1}{y_2}$
\begin{align*}
|S^1| & \frac{\pi}{2} \left|G_{range}(a,b)\right|\\
&=\left|\uint_{\|y\|>R} \cos\(b\sqrt{y^2-R^2}\)f(a+y_1,y_2) \, dy\right|\\
&\leq \uint_{\|y\|>R} |f(a+y_1,y_2)| \, dy\mbox{.}
\end{align*}
The last integral is finite, because $f \in \gS$.

\item

Now it will be shown that $g$ can be recovered for $r>R$.

\begin{align*}
\uint_0^\infty&\uint_\R o_{range}^{a,b}(x,r) G_{range}(a,b) \,da\,db\\
&=\uint_0^\infty\uint_\R G_{range}(a,b) \chi_{(R,\infty)}(r) \delta(x-a) \frac{\cos\(b\sqrt{r^2-R^2}\,\)}{\sqrt{r^2-R^2}}\,da\,db\\
&=\uint_0^\infty\uint_\R \chi_{(R,\infty)}(r) \left(\frac{2}{\pi}\uint_\R\uint_R^\infty \delta(a-x') r' \cos\(b\sqrt{r'^2-R^2}\)\right.\\
&\quad\times g(x',r')\,dr' \, dx'\Bigg) \delta(x-a) \frac{\cos\(b\sqrt{r^2-R^2}\,\)}{\sqrt{r^2-R^2}}\,da \, db\\
&=\chi_{(R,\infty)}(r) \uint_0^\infty \frac{2}{\pi} \left( \uint_R^\infty r' \cos\(b\sqrt{r'^2-R^2}\)g(x,r')\,dr' \right)\\
&\quad\times \frac{\cos\(b\sqrt{r^2-R^2}\,\)}{\sqrt{r^2-R^2}}\,db
\end{align*}
With the substitution $r''=\sqrt{r'^2-R^2}$ and $\CT\left[ \cdot g\left(x,\sqrt{\cdot^2+R^2}\right)\right]$ indicating the cosine transform of $r'' g\left(x,\sqrt{r''^2+R^2}\right)$ with respect to $r''$ it follows that
\begin{align*}
\uint_0^\infty&\uint_\R o_{range}^{a,b}(x,r) G_{range}(a,b) \,da\,db\\
&=\chi_{(R,\infty)}(r) \uint_0^\infty \frac{2}{\pi} \left( \uint_0^\infty \cos(br'') r'' g\left(x,\sqrt{r''^2+R^2}\right)\,dr'' \right)\\
&\quad\times \frac{\cos\(b\sqrt{r^2-R^2}\,\)}{\sqrt{r^2-R^2}}\,db\\
&=\frac{\chi_{(R,\infty)}(r)}{\sqrt{r^2-R^2}} \uint_0^\infty \sqrt{\frac{2}{\pi}} \CT\left[ \cdot g\left(x,\sqrt{\cdot^2+R^2}\right)\right](b)\\
&\quad\times \cos\(b\sqrt{r^2-R^2}\,\) \,db\\
&=\frac{\chi_{(R,\infty)}(r)}{\sqrt{r^2-R^2}} \sqrt{r^2-R^2} g\left(x,\sqrt{\sqrt{r^2-R^2}^2+R^2}\right) \\
&=\chi_{(R,\infty)}(r) g\left(x,\sqrt{r^2-R^2+R^2}\right)=\chi_{(R,\infty)}(r) g(x,r)\mbox{.}
\end{align*}

\end{enumerate}

\item

Assume without loss of generality that $g$ is even in $x$.

\begin{enumerate}

\item

First the well-definedness of $G_{even}$ will be shown.
\begin{align*}
|S^1|& \frac{R}{\epsilon_l a} \left|G^l_{even}(a)\right|\\
&=|S^1| \frac{1}{2} \left|<o_{even}^{a,l}(x,r),g(x,r)+g(-x,r)>_{o_{even}}\right| \\
&=|S^1| \frac{1}{2} \left|\uint_0^R \uint_L^\infty \bn\(a\sqrt{x^2-L^2}\) \frac{\cos\(l\frac{\pi}{R}\sqrt{R^2-r^2}\,\)}{\sqrt{R^2-r^2}}\right.\\
&\quad\left. \times (g(x,r)+g(-x,r)) r \sqrt{R^2-r^2} x \, dx \, dr \right|\\
&=|S^1|\left|\uint_L^\infty x \bn\(a\sqrt{x^2-L^2}\)\right.\\
&\quad\left. \times \uint_0^R \cos\(l\frac{\pi}{R}\sqrt{R^2-r^2}\) r g(x,r)\,dr \, dx\right|
\end{align*}
\begin{align*}
&=\left|\uint_L^\infty x \bn\(a\sqrt{x^2-L^2}\)\uint_0^R \cos\(l\frac{\pi}{R}\sqrt{R^2-r^2}\) r \right.\\
&\quad\left.\times \left(\uint_{S^1}f(x+r\xi,r\eta) \, d S_1(\xi,\eta)\right)\,dr \, dx\right|\\
&\leq C \uint_0^R \uint_{S^1} \uint_L^\infty |x f(x+r\xi,r\eta)| \, dx\,d S_1(\xi,\eta)\,dr
\end{align*}
with an appropriate constant $C>0$. The change in the order of integration is valid and the last integral is finite because $f \in \gS$.

\item

Now it will be shown that the even part of $g$ can be recovered for $0 \leq r < R$ and $|x|>L$.

\begin{align*}
&\sum_{l=0}^\infty \uint_0^\infty o_{even}^{a,l}(x,r) G_{even}^l(a)\,da\\
&=\left(\chi_{(-\infty,-L)}(x)+\chi_{(L,\infty)}(x)\right)\chi_{(0,R)}(r)\sum_{l=0}^\infty \uint_0^\infty \bn\(a\sqrt{x^2-L^2}\)\\
&\quad \times \frac{\cos\(l\frac{\pi}{R}\sqrt{R^2-r^2}\,\)}{\sqrt{R^2-r^2}} \left(\frac{\epsilon_l}{2R}\uint_L^\infty a x'  \bn\(a\sqrt{x'^2-L^2}\) \right.\\
&\quad\left. \times\uint_0^R \cos\(l\frac{\pi}{R}\sqrt{R^2-r'^2}\) r' (g(x',r')+g(-x',r')) \,dr'\,dx'\right)da\\
&=\left(\chi_{(-\infty,-L)}(x)+\chi_{(L,\infty)}(x)\right)\chi_{(0,R)}(r)\sum_{l=0}^\infty \uint_0^\infty \bn\(a\sqrt{x^2-L^2}\) \\
&\quad \times \frac{\cos\(l\frac{\pi}{R}\sqrt{R^2-r^2}\,\)}{\sqrt{R^2-r^2}} \left(\frac{\epsilon_l}{R}\uint_L^\infty a x' \bn\(a\sqrt{x'^2-L^2}\) \right.\\
&\quad\left. \times \uint_0^R \cos\(l\frac{\pi}{R}\sqrt{R^2-r'^2}\) r' g(x',r') \,dr'\,dx'\right)da\mbox{.}
\end{align*}
The substitutions $x''=\sqrt{x'^2-L^2}$ and $r''=\sqrt{R^2-r'^2}$ yield with $g_l\(\sqrt{x''^2+L^2}\)$ representing the Fourier coefficients corresponding to $r'' g\(\sqrt{x''^2+L^2},\sqrt{R^2-r''^2}\)$ and $\overline{g_l\(\sqrt{\cdot^2+L^2}\)}$ indicating the Hankel transform of $g_l\(\sqrt{x''^2+L^2}\)$ with respect to $x''$
\begin{align*}
&\sum_{l=0}^\infty \uint_0^\infty o_{even}^{a,l}(x,r) G_{even}^l(a)\,da\\
&=\left(\chi_{(-\infty,-L)}(x)+\chi_{(L,\infty)}(x)\right)\chi_{(0,R)}(r)\sum_{l=0}^\infty \uint_0^\infty \bn\(a\sqrt{x^2-L^2}\) \\
&\quad \times \frac{\cos\(l\frac{\pi}{R}\sqrt{R^2-r^2}\,\)}{\sqrt{R^2-r^2}} \frac{\epsilon_l}{R}\uint_0^\infty a x'' \bn(ax'') \\
&\quad \times \uint_0^R \cos\(l\frac{\pi}{R}r''\) r'' g\(\sqrt{x''^2+L^2},\sqrt{R^2-r''^2}\)\,dr''\,dx''\,da
\end{align*}
\begin{align*}
&=\left(\chi_{(-\infty,-L)}(x)+\chi_{(L,\infty)}(x)\right)\chi_{(0,R)}(r)\sum_{l=0}^\infty \frac{\cos\(l\frac{\pi}{R}\sqrt{R^2-r^2}\,\)}{\sqrt{R^2-r^2}} \\
&\quad \times \uint_0^\infty a \bn\(a\sqrt{x^2-L^2}\) \uint_0^\infty x'' \bn(ax'') g_l\(\sqrt{x''^2+L^2}\)\,dx''\,da
\end{align*}
\begin{align*}
&=\frac{\left(\chi_{(-\infty,-L)}(x)+\chi_{(L,\infty)}(x)\right)\chi_{(0,R)}(r)}{\sqrt{R^2-r^2}} \sum_{l=0}^\infty \cos\(l\frac{\pi}{R}\sqrt{R^2-r^2}\,\) \\
&\quad \times \uint_0^\infty a \bn\(a\sqrt{x^2-L^2}\) \overline{g_l\(\sqrt{\cdot^2+L^2}\)}(a)\,da\\
&=\frac{\left(\chi_{(-\infty,-L)}(x)+\chi_{(L,\infty)}(x)\right)\chi_{(0,R)}(r)}{\sqrt{R^2-r^2}}\\
&\quad \times \sum_{l=0}^\infty \cos\(l\frac{\pi}{R}\sqrt{R^2-r^2}\,\) g_l\(\sqrt{\sqrt{x^2-L^2}^2+L^2}\)\\
&=\frac{\left(\chi_{(-\infty,-L)}(x)+\chi_{(L,\infty)}(x)\right)\chi_{(0,R)}(r)}{\sqrt{R^2-r^2}} \\
&\quad \times \sum_{l=0}^\infty \cos\(l\frac{\pi}{R}\sqrt{R^2-r^2}\,\) g_l(x)
\end{align*}
\begin{align*}
&=\frac{\left(\chi_{(-\infty,-L)}(x)+\chi_{(L,\infty)}(x)\right)\chi_{(0,R)}(r)}{\sqrt{R^2-r^2}} \\
&\quad \times \sqrt{R^2-r^2} \: g\left(x,\sqrt{R^2-\sqrt{R^2-r^2}^2}\right)\\
&=\left(\chi_{(-\infty,-L)}(x)+\chi_{(L,\infty)}(x)\right)\chi_{(0,R)}(r)g(x,r)\mbox{.}
\end{align*}

\end{enumerate}

\item

Assume without loss of generality that $g$ is odd in $x$.

\begin{enumerate}

\item

First the well-definedness of $G_{odd}$ will be shown.
\begin{align*}
|S^1|& \frac{R}{\epsilon_l a} \left|G^l_{odd}(a)\right|\\
&=|S^1| \frac{1}{2} \left|<o_{odd}^{a,l}(x,r),g(x,r)-g(-x,r)>_{o_{odd}}\right|\\
&=|S^1| \frac{1}{2} \left|\uint_0^R \uint_L^\infty x \bn\(a\sqrt{x^2-L^2}\) \frac{\cos\(l\frac{\pi}{R}\sqrt{R^2-r^2}\,\)}{\sqrt{R^2-r^2}}\right.\\
&\quad\left. \times (g(x,r)-g(-x,r)) r \sqrt{R^2-r^2} \frac{1}{x} \, dx \, dr \right| \\
&=|S^1|\left|\uint_L^\infty x \bn\(a\sqrt{x^2-L^2}\)\right.\\
&\quad\left. \times \uint_0^R \cos\(l\frac{\pi}{R}\sqrt{R^2-r^2}\) \frac{1}{x} r g(x,r)\,dr \, dx\right|\\
&=\left|\uint_L^\infty \bn\(a\sqrt{x^2-L^2}\)\uint_0^R \cos\(l\frac{\pi}{R}\sqrt{R^2-r^2}\) r \right.\\
&\quad\left.\times \left(\uint_{S^1}f(x+r\xi,r\eta) \, d S_1(\xi,\eta)\right)\,dr \, dx\right|\\
&\leq C \uint_0^R \uint_{S^1} \uint_L^\infty |f(x+r\xi,r\eta)| \, dx\,d S_1(\xi,\eta)\,dr
\end{align*}
with an appropriate constant $C>0$. The change in the order of integration is valid and the last integral is finite because $f \in \gS$.

\item

Now it will be shown that the odd part of $g$ can be recovered for $0 \leq r < R$ and $|x|>L$.

\begin{align*}
&\sum_{l=0}^\infty \uint_0^\infty o_{odd}^{a,l}(x,r) G_{odd}^l(a)\,da\\
&=\left(\chi_{(-\infty,-L)}(x)+\chi_{(L,\infty)}(x)\right)\chi_{(0,R)}(r)\sum_{l=0}^\infty \uint_0^\infty x \bn\(a\sqrt{x^2-L^2}\) \\
&\quad \times \frac{\cos\(l\frac{\pi}{R}\sqrt{R^2-r^2}\,\)}{\sqrt{R^2-r^2}} \left(\frac{\epsilon_l}{2R}\uint_L^\infty a x' \bn\(a\sqrt{x'^2-L^2}\) \right.\\
&\quad\left. \times \uint_0^R \cos\(l\frac{\pi}{R}\sqrt{R^2-r'^2}\) r' \frac{1}{x'} (g(x',r')-g(-x',r')) \,dr'\,dx'\right)da\\
&=\left(\chi_{(-\infty,-L)}(x)+\chi_{(L,\infty)}(x)\right)\chi_{(0,R)}(r)\sum_{l=0}^\infty \uint_0^\infty x \bn\(a\sqrt{x^2-L^2}\) \\
&\quad \times \frac{\cos\(l\frac{\pi}{R}\sqrt{R^2-r^2}\,\)}{\sqrt{R^2-r^2}} \left(\frac{\epsilon_l}{R}\uint_L^\infty a \bn\(a\sqrt{x'^2-L^2}\) \right.\\
&\quad\left. \times \uint_0^R \cos\(l\frac{\pi}{R}\sqrt{R^2-r'^2}\) r' g(x',r') \,dr'\,dx'\right)da\mbox{.}
\end{align*}
With the substitutions $x''=\sqrt{x'^2-L^2}$, $r''=\sqrt{R^2-r'^2}$, $g_l\(\sqrt{x''^2+L^2}\)$ indicating the Fourier coefficients corresponding to $r'' g\(\sqrt{x''^2+L^2},\sqrt{R^2-r''^2}\)$ and $\overline{\frac{1}{\sqrt{\cdot^2+L^2}} g_l\(\sqrt{\cdot^2+L^2}\)}$ representing the Hankel transform of $\frac{1}{\sqrt{x''+L^2}} g_l\(\sqrt{x''^2+L^2}\)$ it follows that
\begin{align*}
&\sum_{l=0}^\infty \uint_0^\infty o_{odd}^{a,l}(x,r) G_{odd}^l(a)\,da\\
&=\left(\chi_{(-\infty,-L)}(x)+\chi_{(L,\infty)}(x)\right)\chi_{(0,R)}(r)\sum_{l=0}^\infty \uint_0^\infty x \bn\(a\sqrt{x^2-L^2}\) \\
&\quad \times \frac{\cos\(l\frac{\pi}{R}\sqrt{R^2-r^2}\,\)}{\sqrt{R^2-r^2}} \frac{\epsilon_l}{R}\uint_0^\infty a \frac{x''}{\sqrt{x''+L^2}} \bn(ax'')
\end{align*}
\begin{align*}
&\quad \times \uint_0^R \cos\(l\frac{\pi}{R}r''\) r'' g\(\sqrt{x''^2+L^2},\sqrt{R^2-r''^2}\)\,dr''\,dx''\,da\\
&=\left(\chi_{(-\infty,-L)}(x)+\chi_{(L,\infty)}(x)\right)\chi_{(0,R)}(r)\sum_{l=0}^\infty \frac{\cos\(l\frac{\pi}{R}\sqrt{R^2-r^2}\,\)}{\sqrt{R^2-r^2}} \\
&\quad \times x \uint_0^\infty a \bn\(a\sqrt{x^2-L^2}\) \uint_0^\infty x'' \bn(ax'') \\
&\quad \times \frac{1}{\sqrt{x''+L^2}} g_l\(\sqrt{x''^2+L^2}\)\,dx''\,da\\
&=\left(\chi_{(-\infty,-L)}(x)+\chi_{(L,\infty)}(x)\right)\chi_{(0,R)}(r)\sum_{l=0}^\infty \frac{\cos\(l\frac{\pi}{R}\sqrt{R^2-r^2}\,\)}{\sqrt{R^2-r^2}} \\
&\quad \times x \uint_0^\infty a \bn\(a\sqrt{x^2-L^2}\) \overline{\frac{1}{\sqrt{\cdot^2+L^2}} g_l\(\sqrt{\cdot^2+L^2}\)}(a)\,da\\
&=\frac{\left(\chi_{(-\infty,-L)}(x)+\chi_{(L,\infty)}(x)\right)\chi_{(0,R)}(r)}{\sqrt{R^2-r^2}} x \sum_{l=0}^\infty \cos\(l\frac{\pi}{R}\sqrt{R^2-r^2}\,\) \\
&\quad \times \frac{1}{\sqrt{\sqrt{x^2-L^2}^2+L^2}} g_l\(\sqrt{\sqrt{x^2-L^2}^2+L^2}\)\\
&=\frac{\left(\chi_{(-\infty,-L)}(x)+\chi_{(L,\infty)}(x)\right)\chi_{(0,R)}(r)}{\sqrt{R^2-r^2}} \\
&\quad \times \sum_{l=0}^\infty \cos\(l\frac{\pi}{R}\sqrt{R^2-r^2}\,\) g_l(x)\\
&=\frac{\left(\chi_{(-\infty,-L)}(x)+\chi_{(L,\infty)}(x)\right)\chi_{(0,R)}(r)}{\sqrt{R^2-r^2}} \\
&\quad \times \sqrt{R^2-r^2} \: g\(x,\sqrt{R^2-\sqrt{R^2-r^2}^2}\)\\
&=\left(\chi_{(-\infty,-L)}(x)+\chi_{(L,\infty)}(x)\right)\chi_{(0,R)}(r)g(x,r)\mbox{.}
\end{align*}

\end{enumerate}

\end{enumerate}

\end{proof}

Note that this theorem uses essentially the cosine and the Hankel transform and it looks like an overly complicated formulation. However it is essential to choose the functions like this to compute the inversions analytically.

\section{Reconstruction of the orthogonal functions}

The previous chapter showed that the missing data from the unmeasurable region can be projected onto the set of orthogonal functions $\{o_{range}^{a,b}: a \in \R, b \geq 0\} \cup \{o_{even}^{a,l}, o_{odd}^{a,l}: a > 0, l \in \N_0\}$ and can be retrieved again. In the following it is therefore sufficient to consider only the orthogonal functions to examine the error that arises from limited data. Now the reconstructions of these functions are performed.

\begin{defi}

Let $f \in \gS(\R^n \times \R)$, $a \in \R^n$, and $b \in \R$. Then
$(\tau_{(a,\:b)} f)(x,y)=f(x+a,y+b)$.

\end{defi}

\begin{theo}

\label{ghosts}

\begin{enumerate}

\item

Let $a \in \R, b\geq 0$, and
$$g(x,r)=o_{range}^{a,b}(x,r)\mbox{.}$$
Then $Rf=g$ with
$$\hat{f}(\xi,\eta)=\frac{1}{\sqrt{8\pi}}|\eta|\exp^{-ia\xi}\left\{ \begin{array}{lr}
\frac{\cos\(R\sqrt{\xi^2+\eta^2-b^2}\,\)}{\sqrt{\xi^2+\eta^2-b^2}} &
\mbox{ for } \sqrt{\xi^2+\eta^2}>b\\
& \\
0 & \mbox{ for } \sqrt{\xi^2+\eta^2}<b\\
\end{array} \right.$$
and
\begin{align*}
f(x,y)=&\frac{1}{\sqrt{8\pi}} \mbox{H}_y \frac{\pa}{\pa y}\\
&\quad\times \left\{ \begin{array}{lr}
\frac{\cos\(b\sqrt{(x-a)^2+y^2-R^2}\,\)}{\sqrt{(x-a)^2+y^2-R^2}} &
\mbox{ for } \sqrt{(x-a)^2+y^2}>R\\
& \\
0 & \mbox{ for } \sqrt{(x-a)^2+y^2}<R\mbox{,}\\
\end{array} \right.
\end{align*}
where $\mbox{H}_y$ refers to the Hilbert transform in y.

\item

Let $a \geq 0$, $l \in \N_0$, and
$$g(x,r)=o_{even}^{a,l}(x,r)\mbox{.}$$
Then $Rf=g$ with
\begin{align*}
\hat{f}(\xi,\eta)=&\frac{1}{\sqrt{2\pi}}|\eta|\frac{\sin\(R\sqrt{\xi^2+\eta^2+(l\frac{\pi}{R})^2}\,\)}{\sqrt{\xi^2+\eta^2+(l\frac{\pi}{R})^2}}\\
&\quad\times \left\{ \begin{array}{lr}
\frac{\exp^{-L\sqrt{a^2-\xi^2}}}{\sqrt{a^2-\xi^2}} &
\mbox{ for }0<|\xi|<a\\
& \\
-\frac{\sin\(L\sqrt{\xi^2-a^2}\,\)}{\sqrt{\xi^2-a^2}} & \mbox{ for }|\xi|>a\\
\end{array} \right.
\end{align*}
and
$$f(x,y)=\sqrt{\frac{2}{\pi}} \mbox{H}_y \frac{\pa}{\pa y} \uint_\R \left( \left\{ \begin{array}{lr}
\bn\(a\sqrt{t^2-L^2}\,\) &
\mbox{ for } |t|>L\\
& \\
0 & \mbox{ for } |t|<L\\
\end{array} \right\}\right.$$
$$\times \left. \left\{ \begin{array}{lr}
\frac{\cos\(l\frac{\pi}{R}\sqrt{R^2-(x-t)^2-y^2}\,\)}{\sqrt{R^2-(x-t)^2-y^2}} &
\mbox{ for } (x-t)^2+y^2<R^2\\
& \\
0 & \mbox{ for } (x-t)^2+y^2>R^2\\
\end{array} \right\} \right) dt\mbox{.}$$

\item

Let $a \geq 0$, $l \in \N_0$, and
$$g(x,r)=o_{odd}^{a,l}(x,r)\mbox{.}$$
Then $Rf=g$ with
\begin{align*}
\hat{f}(\xi,\eta)=&-\frac{1}{\sqrt{2\pi}}|\eta|\frac{\pa}{\pa \xi}\frac{\sin\(R\sqrt{\xi^2+\eta^2+(l\frac{\pi}{R})^2}\,\)}{\sqrt{\xi^2+\eta^2+(l\frac{\pi}{R})^2}}\\
&\quad\times \left\{ \begin{array}{lr}
\frac{\exp^{-L\sqrt{a^2-\xi^2}}}{\sqrt{a^2-\xi^2}} &
\mbox{ for }0<|\xi|<a\\
& \\
-\frac{\sin\(L\sqrt{\xi^2-a^2}\,\)}{\sqrt{\xi^2-a^2}} & \mbox{ for }|\xi|>a\\
\end{array} \right.
\end{align*}
and
\begin{align*}
f(x,y)=&\sqrt{\frac{2}{\pi}} x \mbox{H}_y \frac{\pa}{\pa y} \uint_\R \left( \left\{ \begin{array}{lr}
\bn\(a\sqrt{t^2-L^2}\,\) &
\mbox{ for } |t|>L\\
& \\
0 & \mbox{ for } |t|<L\\
\end{array} \right\}\right.\\
&\quad\times \left. \left\{ \begin{array}{lr}
\frac{\cos\(l\frac{\pi}{R}\sqrt{R^2-(x-t)^2-y^2}\,\)}{\sqrt{R^2-(x-t)^2-y^2}} &
\mbox{ for } (x-t)^2+y^2<R^2\\
& \\
0 & \mbox{ for } (x-t)^2+y^2>R^2\\
\end{array} \right\} \right) dt\mbox{.}
\end{align*}

\end{enumerate}

\end{theo}

\begin{proof}

\begin{enumerate}

\item

Let $a \in \R, b\geq 0$, and
$$g(x,r)=o_{range}^{a,b}(x,r)\mbox{.}$$
Then the three dimensional Fourier transform of $g$ is given by
$$\hat{g}(\xi,\rho)=\frac{1}{\sqrt{2\pi}}\uint_\R \uint_0^\infty \exp^{-ix\xi}r\bn(r\rho)g(x,r)\,dr \, dx$$
and with lemma \ref{hankelunendlich}
$$\hat{g}(\xi,\rho)=\frac{1}{\sqrt{2\pi}}\exp^{-ia\xi}\left\{ \begin{array}{lr}
\frac{\cos\(R\sqrt{\rho^2-b^2}\,\)}{\sqrt{\rho^2-b^2}} &
\mbox{ for } \rho>b\\
& \\
0 & \mbox{ for } \rho<b \mbox{.}\\
\end{array} \right.$$
Therefore, according to theorem \ref{inversionformula}
\begin{align*}
\hat{f}(\xi,\eta)&=\frac{1}{2}|\eta|\hat{g}\(\xi,\sqrt{\xi^2+\eta^2}\)\\
&=\frac{1}{\sqrt{8\pi}}|\eta|\exp^{-ia\xi}\left\{ \begin{array}{lr}
\frac{\cos\(R\sqrt{\xi^2+\eta^2-b^2}\,\)}{\sqrt{\xi^2+\eta^2-b^2}} &
\mbox{ for } \sqrt{\xi^2+\eta^2}>b\\
& \\
0 & \mbox{ for } \sqrt{\xi^2+\eta^2}<b \mbox{.}\\
\end{array} \right.
\end{align*}
Hence
\begin{align*}
f(x&,y)=\frac{1}{2\pi} \uint_\R \uint_\R \exp^{ix\xi} \exp^{iy\eta} \hat{f}(\xi,\eta)\,d\eta \, d\xi\\
&=\frac{1}{2\pi} \uint_\R \uint_\R \exp^{ix\xi} \exp^{iy\eta}\\
&\quad\times \left(\frac{1}{\sqrt{8\pi}}|\eta|\exp^{-ia\xi}\frac{\cos\(R\sqrt{\xi^2+\eta^2-b^2}\,\)}{\sqrt{\xi^2+\eta^2-b^2}}\chi_{(b^2,\infty)}(\xi^2+\eta^2)\right)\,d\eta \, d\xi\\
&=\frac{1}{2\pi} \frac{1}{\sqrt{8\pi}} \tau_{(-a,0)} \mbox{H}_y \frac{\pa}{\pa y} \uint_\R \uint_\R \exp^{ix\xi} \exp^{iy\eta}\\
&\quad\times \frac{\cos\(R\sqrt{\xi^2+\eta^2-b^2}\,\)}{\sqrt{\xi^2+\eta^2-b^2}}\chi_{(b^2,\infty)}(\xi^2+\eta^2)\,d\eta \, d\xi\\
&=\frac{1}{\sqrt{8\pi}} \tau_{(-a,0)} \mbox{H}_y \frac{\pa}{\pa y} \uint_b^\infty \rho \bn(\rho\sqrt{x^2+y^2}) \frac{\cos\(R\sqrt{\rho^2-b^2}\,\)}{\sqrt{\rho^2-b^2}}\,d\rho\mbox{.}
\end{align*}
With lemma \ref{hankelunendlich} and because the Hankel transform is its own inverse, this leads to
\begin{align*}
f(x,y)&=\frac{1}{\sqrt{8\pi}} \tau_{(-a,0)}\mbox{H}_y \frac{\pa}{\pa y}  \left\{ \begin{array}{lr}
\frac{\cos\(b\sqrt{x^2+y^2-R^2}\,\)}{\sqrt{x^2+y^2-R^2}} &
\mbox{ for } \sqrt{x^2+y^2}>R\\
& \\
0 & \mbox{ for } \sqrt{x^2+y^2}<R\\
\end{array} \right.\\
&=\frac{1}{\sqrt{8\pi}} \mbox{H}_y \frac{\pa}{\pa y} \left\{ \begin{array}{lr}
\frac{\cos\(b\sqrt{(x-a)^2+y^2-R^2}\,\)}{\sqrt{(x-a)^2+y^2-R^2}} &
\mbox{ for } \sqrt{(x-a)^2+y^2}>R\\
& \\
0 & \mbox{ for } \sqrt{(x-a)^2+y^2}<R\mbox{.}\\
\end{array} \right.
\end{align*}

\item

Let $a \geq 0$, $l \in \N_0$, and
$$g(x,r)=o_{even}^{a,l}(x,r)\mbox{.}$$
Then the three dimensional Fourier transform of $g$ is given by
$$\hat{g}(\xi,\rho)=\frac{1}{\sqrt{2\pi}}\uint_\R \uint_0^\infty \exp^{-ix\xi}r\bn(r\rho)g(x,r)\,dr \, dx$$
and with lemmas \ref{cosinusunendlich} and \ref{hankelendlich}
$$\hat{g}(\xi,\rho)=\sqrt{\frac{2}{\pi}}\frac{\sin\(R\sqrt{\rho^2+(l\frac{\pi}{R})^2}\,\)}{\sqrt{\rho^2+(l\frac{\pi}{R})^2}}\left\{ \begin{array}{lr}
\frac{\exp^{-L\sqrt{a^2-\xi^2}}}{\sqrt{a^2-\xi^2}} &
\mbox{ for }0<|\xi|<a\\
& \\
-\frac{\sin\(L\sqrt{\xi^2-a^2}\,\)}{\sqrt{\xi^2-a^2}} & \mbox{ for }|\xi|>a \mbox{.}\\
\end{array} \right.$$
Thus, with theorem \ref{inversionformula} it follows that
\begin{align*}
\hat{f}(\xi,\eta)&=\frac{1}{2}|\eta|\hat{g}\(\xi,\sqrt{\xi^2+\eta^2}\)\\
&=\frac{1}{\sqrt{2\pi}}|\eta|\frac{\sin\(R\sqrt{\xi^2+\eta^2+(l\frac{\pi}{R})^2}\,\)}{\sqrt{\xi^2+\eta^2+(l\frac{\pi}{R})^2}}\\
&\quad\times \left\{ \begin{array}{lr}
\frac{\exp^{-L\sqrt{a^2-\xi^2}}}{\sqrt{a^2-\xi^2}} &
\mbox{ for }0<|\xi|<a\\
& \\
-\frac{\sin\(L\sqrt{\xi^2-a^2}\,\)}{\sqrt{\xi^2-a^2}} & \mbox{ for }|\xi|>a \mbox{.}\\
\end{array} \right.
\end{align*}
Therefore
$$f(x,y)=\frac{1}{2\pi} \uint_\R \uint_\R \exp^{ix\xi} \exp^{iy\eta} \hat{f}(\xi,\eta)\,d\eta \, d\xi\mbox{.}$$
Applying the Fourier convolution theorem leads to
\begin{align*}
f&(x,y)=\frac{1}{(2\pi)^{3/2}} \mbox{H}_y \frac{\pa}{\pa y} \uint_\R\\
&\quad\times \left( \uint_\R \exp^{it\xi}\left\{ \begin{array}{lr}
\frac{\exp^{-L\sqrt{a^2-\xi^2}}}{\sqrt{a^2-\xi^2}} &
\mbox{ for }0<|\xi|<a\\
& \\
-\frac{\sin\(L\sqrt{\xi^2-a^2}\,\)}{\sqrt{\xi^2-a^2}} & \mbox{ for }|\xi|>a\\
\end{array} \right\}\right.\,d\xi\\
&\quad\times \left. \uint_\R \uint_\R \exp^{i(x-t)\xi} \exp^{iy\eta} \frac{\sin\(R\sqrt{\xi^2+\eta^2+(l\frac{\pi}{R})^2}\,\)}{\sqrt{\xi^2+\eta^2+(l\frac{\pi}{R})^2}}\,d\xi \, d\eta \right) dt\\
&=\sqrt{\frac{2}{\pi}} \mbox{H}_y \frac{\pa}{\pa y} \uint_\R \left( \uint_0^\infty \cos(t\xi)\left\{ \begin{array}{lr}
\frac{\exp^{-L\sqrt{a^2-\xi^2}}}{\sqrt{a^2-\xi^2}} &
\mbox{ for }0<|\xi|<a\\
& \\
-\frac{\sin\(L\sqrt{\xi^2-a^2}\,\)}{\sqrt{\xi^2-a^2}} & \mbox{ for }|\xi|>a\\
\end{array} \right\}\right.\,d\xi\\
&\quad\times \left. \uint_0^\infty \bn\(\rho\sqrt{(x-t)^2+y^2}\) \frac{\sin\(R\sqrt{\rho^2+(l\frac{\pi}{R})^2}\,\)}{\sqrt{\rho^2+(l\frac{\pi}{R})^2}}\,d\rho \right) dt\mbox{.}
\end{align*}
With lemmas \ref{cosinusunendlich} and \ref{hankelendlich} and bearing in mind that the Hankel transform and the cosine transform are their respective inverses this leads to
\begin{align*}
f(x,y)=&\sqrt{\frac{2}{\pi}} \mbox{H}_y \frac{\pa}{\pa y} \uint_\R \left( \left\{ \begin{array}{lr}
\bn\(a\sqrt{t^2-L^2}\,\) &
\mbox{ for } |t|>L\\
& \\
0 & \mbox{ for } |t|<L\\
\end{array} \right\}\right.\\
&\quad\times \left. \left\{ \begin{array}{lr}
\frac{\cos\(l\frac{\pi}{R}\sqrt{R^2-(x-t)^2-y^2}\,\)}{\sqrt{R^2-(x-t)^2-y^2}} &
\mbox{ for } (x-t)^2+y^2<R^2\\
& \\
0 & \mbox{ for } (x-t)^2+y^2>R^2\\
\end{array} \right\} \right) dt\mbox{.}
\end{align*}

\item

Let $a \geq 0$, $l \in \N_0$, and
$$g(x,r)=o_{odd}^{a,l}(x,r)\mbox{.}$$
Then the proof is analogous to 2. using $x\sin(x\xi)=-\frac{\pa}{\pa \xi}\cos(x\xi)$.

\end{enumerate}

\end{proof}

\section{Numerical simulations}

To imbue the unwieldy formulas from theorem \ref{ghosts} with life, in the following some of the ghosts - functions that are in the null space of the measurement operator due to limited data - derived in this theorem are shown. All of the following examples can therefore be added to the reflectivity function $f$ without changing the measured data. In all cases it is assumed that $L=R=1$, and the region $(x,y) \in [0,1] \times [0,1]$ is shown.
At first, examples for $g \in \{o_{range}^{a,b} : a \in \R, b \geq 0\}$ are shown, then examples for $g \in \{o_{even}^{a,l} : a > 0, l \in \N_0\}$. The examples for $g \in \{o_{odd}^{a,l} : a > 0, l \in \N_0\}$ are omitted because they differ from the examples of $g \in \{o_{even}^{a,l} : a > 0, l \in \N_0\}$ only by an additional factor $x$.

\subsection{Reconstructions for $g \in \{o_{range}^{a,b} : a \in \R, b \geq 0\}$}

Now two series of images are shown to demonstrate the effects of the parameters $a$ and $b$ for ghosts derived from functions  $g \in \{o_{range}^{a,b} : a \in \R, b \geq 0\}$. These are the artifacts that can arise due to the limitation that the echoes from the emitted waves can only be received up to the limited distance $R$. In the first series $b$ is constant at $0.25$, and $a$ is varied from $-0.1$ to $1.2$. That range is sufficient because reconstructions for different values of $a$ can be obtained via translations in $x$, as can be seen from the formula in theorem \ref{ghosts},1. In the second series $a$ is constant at $0.6$, and $b$ takes the values of $0.25$ and $1$. Since the functions $g \in \{o_{range}^{a,b} : a \in \R, b \geq 0\}$ all have a singularity for $r=R$ that would usually not be reflected in real data and that would only obscure the interesting features, only the differences $o_{range}^{a,b}-o_{range}^{a,0}$ are reconstructed. Thereby the singularity is removed and that should lead to ghosts similar to those that can be expected in reconstructions from real data. The maximum displayed in the following images is capped at a reasonable level to enhance the visibility of the features. Otherwise a couple of large peaks would distort the image. These peaks are located at the points where the course of the depicted circles is close to a grid point. It should be noted that the following reconstructions are only approximative because the Hilbert transform from theorem \ref{ghosts} does only exist in the distributional sense, even using the differences mentioned above. This should not be a problem however since for real data the integration with the kernel $\bar{g}_r$ as in theorem \ref{retrieve} should lead to a regular function.


\fbildr{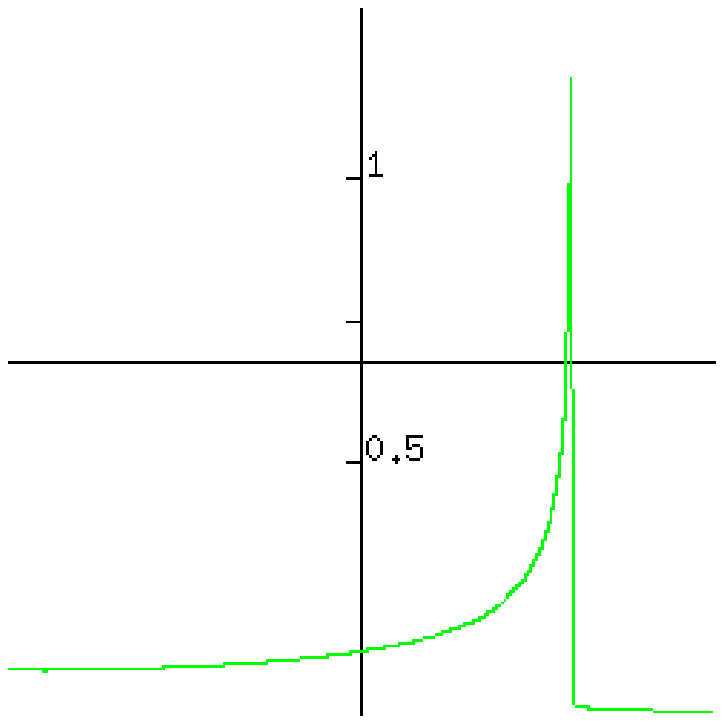}{Reconstructed $f$ for $g=Rf=o_{range}^{-0.1,0.25}-o_{range}^{-0.1,0}$\label{fr--0.1-0.25}}
\fbildr{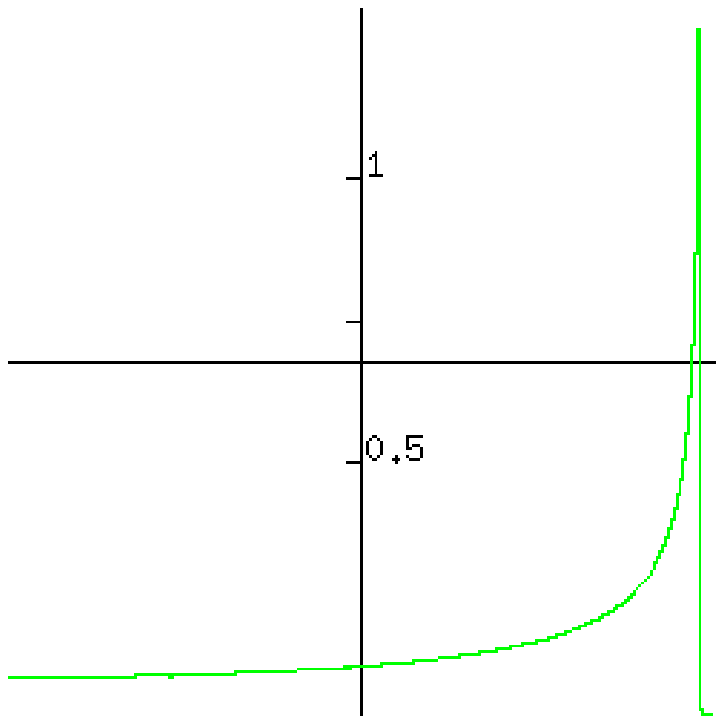}{Reconstructed $f$ for $g=Rf=o_{range}^{0.3,0.25}-o_{range}^{0.3,0}$\label{fr-0.3-0.25}}
\fbildr{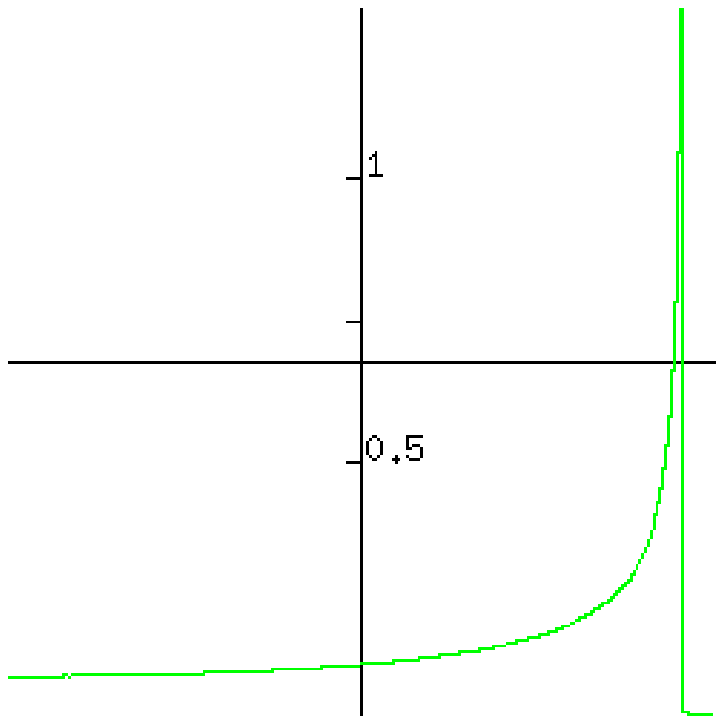}{Reconstructed $f$ for $g=Rf=o_{range}^{0.8,0.25}-o_{range}^{0.8,0}$\label{fr-0.8-0.25}}
\fbildr{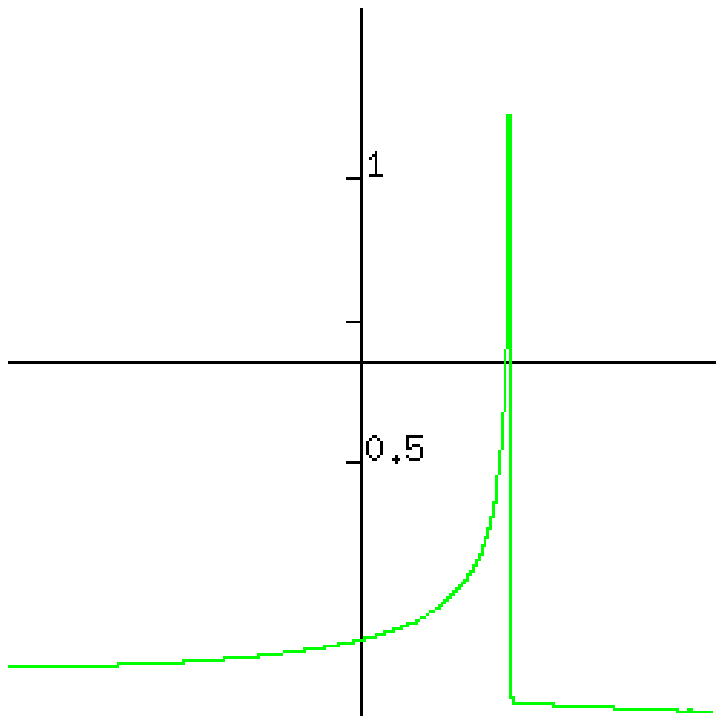}{Reconstructed $f$ for $g=Rf=o_{range}^{1.2,0.25}-o_{range}^{1.2,0}$\label{fr-1.2-0.25}}

Figures \ref{fr--0.1-0.25} to \ref{fr-1.2-0.25} each show a reconstruction of $o_{range}^{a,0.25}-o_{range}^{a,0}$. Figure \ref{fr--0.1-0.25} displays the result for $a=-0.1$, figure \ref{fr-0.3-0.25} for $a=0.3$, figure \ref{fr-0.8-0.25} for $a=0.8$, and figure \ref{fr-1.2-0.25} for $a=1.2$. It can be clearly seen in figures \ref{fr--0.1-0.25} to \ref{fr-1.2-0.25} that the missing data in this parameter range causes a circular artifact that is centered at $(a,0)$. These artifacts remind of the artifacts seen in chapter \ref{Diplomkorrektur}, so this is a possible explanation. It should be noted that the reconstructed functions are not $0$ for $(x-a)^2 + y^2 > R^2$ although some reconstructions could lead to this assumption, but this is only an effect of the overshadowing singularity. In the cross section of figure \ref{fr-1.2-0.25} it is clearly visible that the reconstruction is greater than $0$ for $(x-a)^2 + y^2 > R^2$.


\fbildr{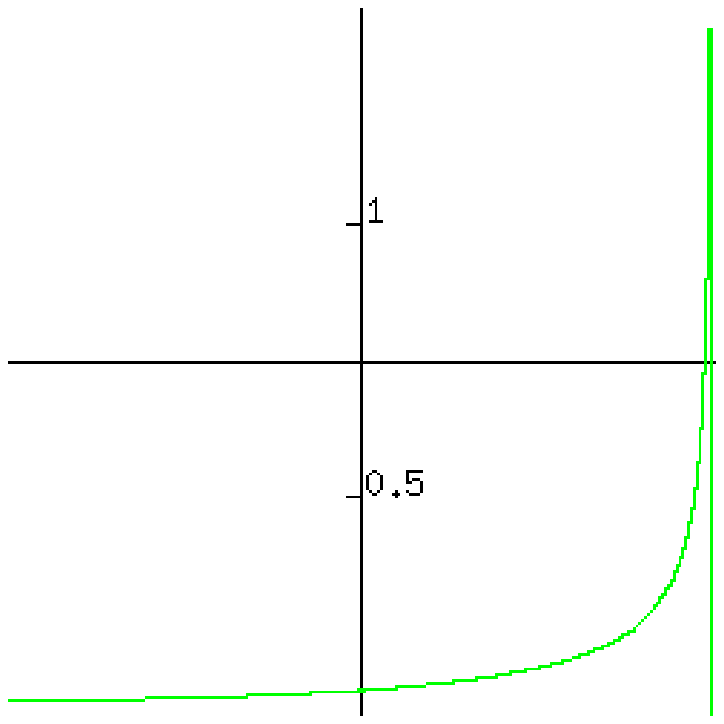}{Reconstructed $f$ for $g=Rf=o_{range}^{0.6,0.25}-o_{range}^{0.6,0}$\label{fr-0.6-0.25}}
\fbildr{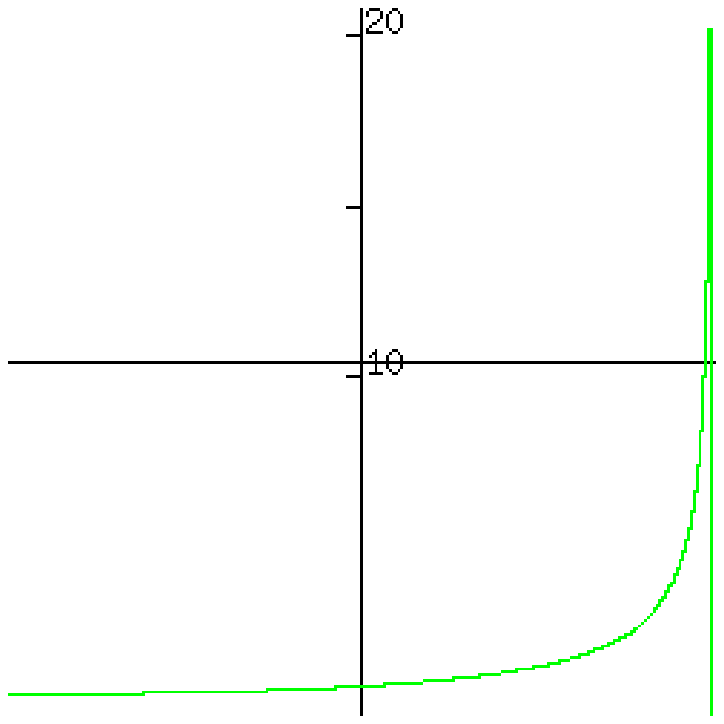}{Reconstructed $f$ for $g=Rf=o_{range}^{0.6,1}-o_{range}^{0.6,0}$\label{fr-0.6-1}}

Figures \ref{fr-0.6-0.25} and \ref{fr-0.6-1} each show a reconstruction of $o_{range}^{0.6,b}-o_{range}^{0.6,0}$. Figure \ref{fr-0.6-0.25} delineates the result for $b=0.25$ and figure \ref{fr-0.6-1} for $b=1$. As can be seen, apart from the scaling figures \ref{fr-0.6-0.25} and \ref{fr-0.6-1} are very similar. They depict a circle with center $(0.6,0)$ and radius $r$. It can be seen in the cross section that the amplitude changes due to the variation in the second parameter.



\subsection{Reconstructions for $g \in \{o_{even}^{a,l} : a > 0, l \in \N_0\}$}


In the following two groups of images, functions $g \in \{o_{even}^{a,l} : a > 0, l \in \N_0\}$ are shown. The first group, for $l=1$, comprises examples for $a \in \{1,4,16\}$ as does the second group, however for $l=16$. As in the preceding subsection, the functions $g \in \{o_{even}^{a,l} : a > 0, l \in \N_0\}$ all have a singularity for $r=R$ that would usually not be reflected in real data and that would only obscure the interesting features. This problem is solved in the same way so that the shown examples should lead to ghosts that can be expected in real data. As in the previous subsection it should be noted that only an approximation is computed. But again this should not be a problem as the same circumstances hold.

\fbildt{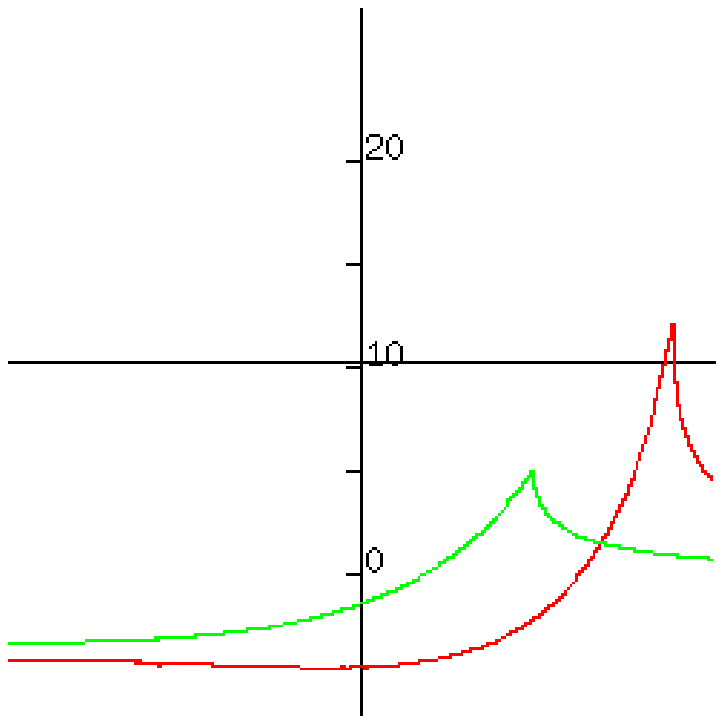}{Reconstructed $f$ for $g=Rf=o_{even}^{1,1}-o_{even}^{1,0}$\label{ft-1-1}}
\fbildt{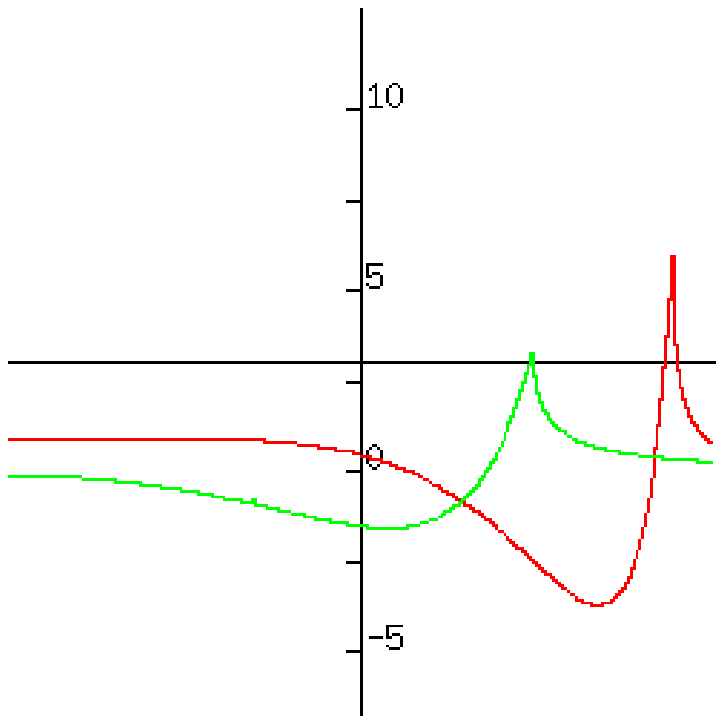}{Reconstructed $f$ for $g=Rf=o_{even}^{4,1}-o_{even}^{4,0}$\label{ft-4-1}}
\fbildt{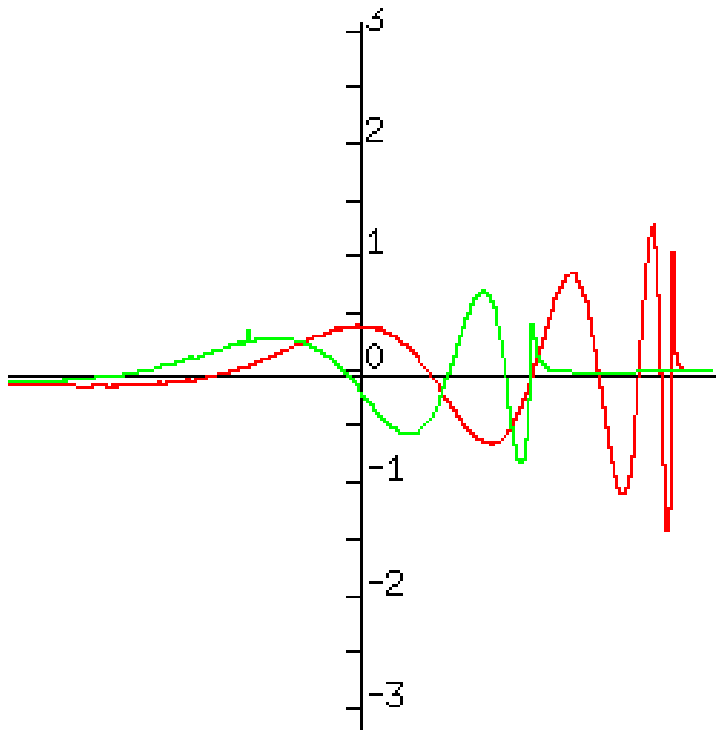}{Reconstructed $f$ for $g=Rf=o_{even}^{16,1}-o_{even}^{16,0}$\label{ft-16-1}}

Figures \ref{ft-1-1} to \ref{ft-16-1} each show a reconstruction of $o_{even}^{a,1}-o_{even}^{a,0}$. Figure \ref{ft-1-1} depicts the result for $a=1$, figure \ref{ft-4-1} for $a=4$, and figure \ref{ft-16-1} for $a=16$. Figures \ref{ft-1-1} to \ref{ft-16-1} again show circular artifacts, but of a different kind. These artifacts show the biggest variation in amplitude around $x^2+y^2=L^2$. However the artifacts clearly extend into the region $x^2+y^2<L^2$. Figures \ref{ft-1-1} and \ref{ft-4-1} remind of artifacts that appear in reconstructions of objects close to the edges of the flight track whereas figure \ref{ft-16-1} exhibits a higher frequency phenomenon that was not encountered so far. In comparing figures \ref{ft-1-1} to \ref{ft-16-1} one can conclude that with larger $a$ the amplitude of the artifact gets smaller and the frequency of the artifacts becomes larger. This is understandable as the parameter $a$ affects the frequency in the Bessel function.

\fbildt{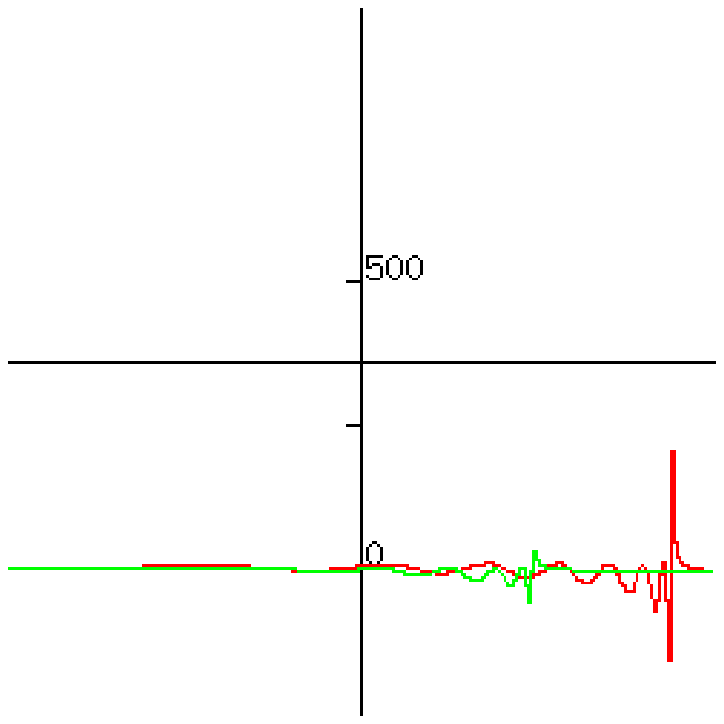}{Reconstructed $f$ for $g=Rf=o_{even}^{1,16}-o_{even}^{1,0}$\label{ft-1-16}}
\fbildt{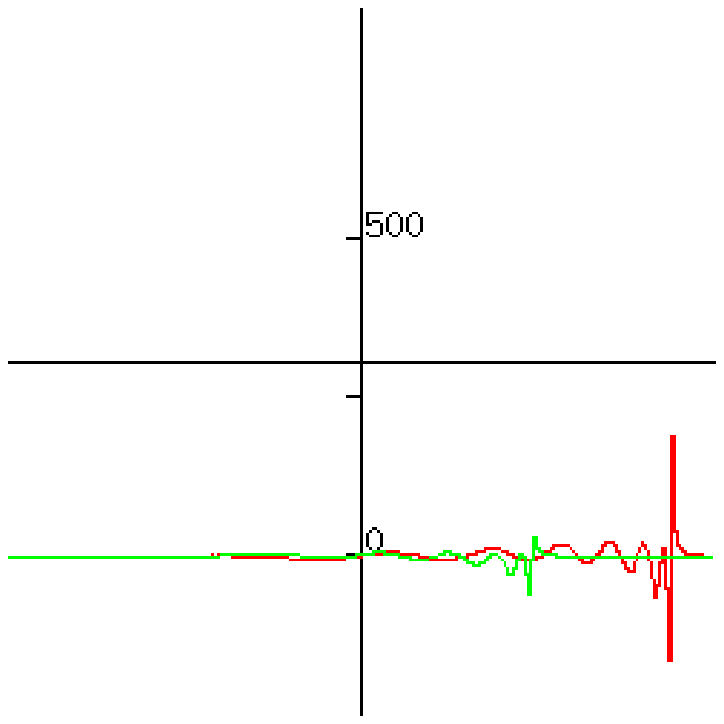}{Reconstructed $f$ for $g=Rf=o_{even}^{4,16}-o_{even}^{4,0}$\label{ft-4-16}}
\fbildt{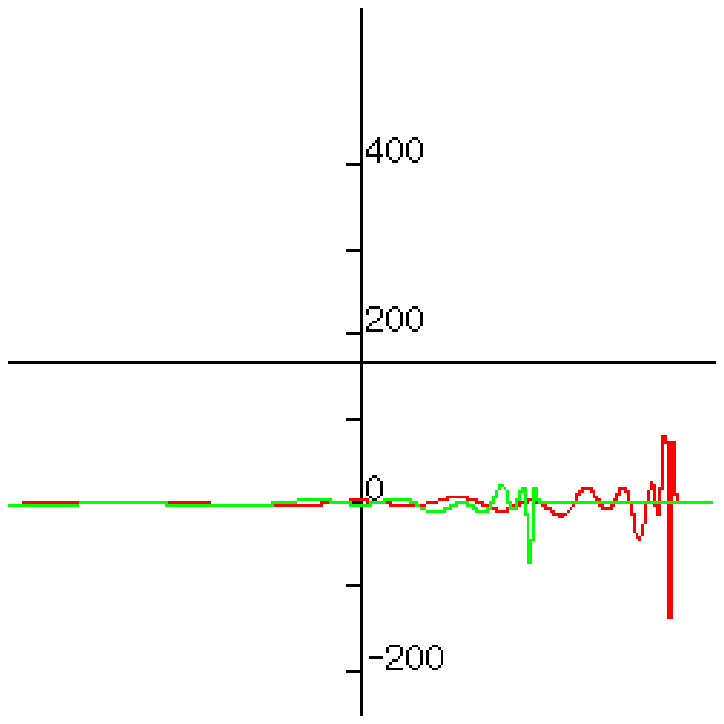}{Reconstructed $f$ for $g=Rf=o_{even}^{16,16}-o_{even}^{16,0}$\label{ft-16-16}}

Figures \ref{ft-1-16} to \ref{ft-16-16} each show a reconstruction of $o_{even}^{a,16}-o_{even}^{a,0}$. Figure \ref{ft-1-16} displays the result for $a=1$, figure \ref{ft-4-16} for $a=4$, and figure \ref{ft-16-16}  for $a=16$. For these figures the same holds as for figures \ref{ft-1-1} to \ref{ft-16-1}. Again, similar to the preceding subsection, an increase of the second parameter leads to a higher amplitude as demonstrated by figures \ref{ft-1-16} to \ref{ft-16-16}. Additionally, the artifacts seem to get sharper.\\
No reconstructions for the functions $g \in \{o_{odd}^{a,l} : a > 0, l \in \N_0\}$ are shown since they would differ from the examples in this subsection only by an additional gradient in $x$ direction as can be clearly seen in theorem \ref{ghosts}.\\

The computed ghosts partially show a close resemblance to the artifacts known from reconstructions. Unfortunately the analysis of this chapter shows that the common and unavoidable problem of limited data is not restricted to high frequency artifacts. This urges a thorough examination of how these artifacts can be avoided, e. g. by a regularization using information from the physical background. The knowledge about the nature of the ghosts that was derived in this chapter should hopefully simplify this task.

\chapter{A new approach to invert the spherical Radon transform}

\label{Fast}

With the ideas from the last chapter, it is possible to reconstruct the recoverable part of $f$ using orthogonal functions. In the following, a set of two dimensional, orthogonal functions with a compact support will be introduced. Then it will be shown that for functions that fulfill the properties of a measurement the information gathered by the projections of the data onto the orthogonal functions is sufficient to regain the data. Finally, the inversions of the orthogonal functions will be calculated. With these results an alternative way to reconstruct the images is obtained: First, the data is projected onto the orthogonal functions, and then the inversions of the orthogonal functions are summed using the coefficients obtained from the projections.

\section{Orthogonal functions and their transforms}

Now a set of functions is defined, and it is shown that they are orthogonal.

\begin{lem}

\label{cosinusendlich}

Let $L > 0, 0 < x < L$, and $a \geq 0$. Then with the Bessel function of order zero $\bn$
$$\uint_0^L \frac{\cos\(a\sqrt{L^2-x^2}\,\)}{\sqrt{L^2-x^2}}\cos(x\xi)\,dx=\frac{\pi}{2}\bn\(L\sqrt{a^2+\xi^2}\)\mbox{.}$$

\end{lem}

\begin{proof}

\cite[I, \S 5, p. 30]{Oberhettinger1}

\end{proof}

\begin{defi}

Let $R, L >0$ and $k, l \in \N_0$. Define

\begin{enumerate}

\item

$$i_{even}^{k,l}(x,r):=\chi_{(0,L)}(x) \chi_{(0,R)}(r) \frac{\cos\(k\frac{\pi}{L}\sqrt{L^2-x^2}\,\)}{\sqrt{L^2-x^2}} \frac{\cos\(l\frac{\pi}{R}\sqrt{R^2-r^2}\,\)}{\sqrt{R^2-r^2}} \mbox{.}$$

\item

$$i_{odd}^{k,l}(x,r):=\chi_{(0,L)}(x) \chi_{(0,R)}(r) x \frac{\cos\(k\frac{\pi}{L}\sqrt{L^2-x^2}\,\)}{\sqrt{L^2-x^2}} \frac{\cos\(l\frac{\pi}{R}\sqrt{R^2-r^2}\,\)}{\sqrt{R^2-r^2}} \mbox{.}$$

\end{enumerate}

\end{defi}

\begin{prop}

Let $R, L >0$, $k, k', l, l' \in \N_0$, $k \neq k'$, and $l \neq l'$. Then

\begin{enumerate}

\item

the functions $i_{even}^{k,l}(x,r)$ and $i_{even}^{k',l'}(x,r)$ are orthogonal to each other with respect to the scalar product
$<f,g>_{i_{even}}=\uint_0^R \uint_0^L f(x,r) g(x,r) r \sqrt{R^2-r^2} \: x \: \sqrt{L^2-x^2} \, dx \, dr$.

\item

the functions $i_{odd}^{k,l}(x,r)$ and $i_{odd}^{k',l'}(x,r)$ are orthogonal to each other with respect to the scalar product
$<f,g>_{i_{odd}}=\uint_0^R \uint_0^L f(x,r) g(x,r) r \sqrt{R^2-r^2} \frac{\sqrt{L^2-x^2}}{x} \, dx \, dr$.

\end{enumerate}

\end{prop}

\begin{proof}

Let $R, L >0$ and $k, k', l, l' \in \N_0$. Then

\begin{enumerate}

\item

\begin{align*}
<i_{even}^{k,l},&i_{even}^{k',l'}>_{i_{even}}=\uint_0^R \uint_0^L i_{even}^{k,l}(x,r) i_{even}^{k',l'}(x,r) r \sqrt{R^2-r^2} \: x \: \sqrt{L^2-x^2} \, dx \, dr\\
&=\uint_0^R \uint_0^L \frac{\cos\(k\frac{\pi}{L}\sqrt{L^2-x^2}\,\)}{\sqrt{L^2-x^2}} \frac{\cos\(l\frac{\pi}{R}\sqrt{R^2-r^2}\,\)}{\sqrt{R^2-r^2}}  \frac{\cos\(k'\frac{\pi}{L}\sqrt{L^2-x^2}\,\)}{\sqrt{L^2-x^2}}\\
&\quad\times \frac{\cos\(l'\frac{\pi}{R}\sqrt{R^2-r^2}\,\)}{\sqrt{R^2-r^2}} \: r \: \sqrt{R^2-r^2} \: x \: \sqrt{L^2-x^2} \, dx \, dr\\
&=\uint_0^R \uint_0^L \cos\(k\frac{\pi}{L}\sqrt{L^2-x^2}\) \cos\(l\frac{\pi}{R}\sqrt{R^2-r^2}\) \cos\(k'\frac{\pi}{L}\sqrt{L^2-x^2}\)\\
&\quad\times \cos\(l'\frac{\pi}{R}\sqrt{R^2-r^2}\) \frac{r}{\sqrt{R^2-r^2}} \frac{x}{\sqrt{L^2-x^2}} \, dx \, dr\mbox{.}
\end{align*}
With the substitutions $r'=\sqrt{R^2-r^2}$ and $x'=\sqrt{L^2-x^2}$
\begin{align*}
<i_{even}^{k,l},i_{even}^{k',l'}>_{i_{even}}&=\uint_0^R \uint_0^L  \cos(k\frac{\pi}{L}x') \cos(l\frac{\pi}{R}r') \cos(k'\frac{\pi}{L}x') \cos(l'\frac{\pi}{R}r') \, dx' \, dr'\\
&=\frac{L}{\epsilon_k}\delta_{kk'} \frac{R}{\epsilon_l}\delta_{ll'}
\end{align*}
with the Neumann's numbers $\epsilon_k, \epsilon_l$.

\item

\begin{align*}
<&i_{odd}^{k,l},i_{odd}^{k',l'}>_{i_{odd}}=\uint_0^R \uint_0^L i_{odd}^{k,l}(x,r) i_{odd}^{k',l'}(x,r) r \sqrt{R^2-r^2} \frac{\sqrt{L^2-x^2}}{x} \, dx \, dr\\
&=\uint_0^R \uint_0^L x \frac{\cos\(k\frac{\pi}{L}\sqrt{L^2-x^2}\,\)}{\sqrt{L^2-x^2}} \: \frac{\cos\(l\frac{\pi}{R}\sqrt{R^2-r^2}\,\)}{\sqrt{R^2-r^2}} \: x \: \frac{\cos\(k'\frac{\pi}{L}\sqrt{L^2-x^2}\,\)}{\sqrt{L^2-x^2}}\\
&\quad\times \frac{\cos\(l'\frac{\pi}{R}\sqrt{R^2-r^2}\,\)}{\sqrt{R^2-r^2}} r \sqrt{R^2-r^2} \frac{\sqrt{L^2-x^2}}{x} \, dx \, dr\\
&=\uint_0^R \uint_0^L \cos\(k\frac{\pi}{L}\sqrt{L^2-x^2}\) \cos\(l\frac{\pi}{R}\sqrt{R^2-r^2}\) \cos\(k'\frac{\pi}{L}\sqrt{L^2-x^2}\)\\
&\quad\times \cos\(l'\frac{\pi}{R}\sqrt{R^2-r^2}\) \frac{r}{\sqrt{R^2-r^2}} \frac{x}{\sqrt{L^2-x^2}} \, dx \, dr \mbox{.}
\end{align*}
With the substitutions $r'=\sqrt{R^2-r^2}$ and $x'=\sqrt{L^2-x^2}$
\begin{align*}
<i_{odd}^{k,l},i_{odd}^{k',l'}>_{i_{odd}}&=\uint_0^R \uint_0^L  \cos(k\frac{\pi}{L}x') \cos(l\frac{\pi}{R}r') \cos(k'\frac{\pi}{L}x') \cos(l'\frac{\pi}{R}r') \, dx' \, dr'\\
&=\frac{L}{\epsilon_k}\delta_{kk'} \frac{R}{\epsilon_l}\delta_{ll'}
\end{align*}
with the Neumann's numbers $\epsilon_k, \epsilon_l$.

\end{enumerate}

\end{proof}

\section{Projection of the data onto the orthogonal functions}

It will be shown that the measured data can be projected onto the set of orthogonal functions and that it is completely recoverable from these projections.

\begin{defi}

Let $f \in \gS$, $g=Rf$, and let $\epsilon_k, \epsilon_l$ denote the Neumann's numbers.

\begin{enumerate}

\item

$$G_{even}^{k,l}:=\frac{\epsilon_k\epsilon_l}{2LR} <i_{even}^{k,l}(x,r),g(x,r)+g(-x,r)>_{i_{even}}$$

\item

$$G_{odd}^{k,l}:=\frac{\epsilon_k\epsilon_l}{2LR} <i_{odd}^{k,l}(x,r),g(x,r)-g(-x,r)>_{i_{odd}}$$

\end{enumerate}

\end{defi}

\begin{theo}

\label{grecoverable}

\begin{enumerate}

\item

$G_{even}^{k,l}$ is well defined, and if $f \in \gS$, $g=Rf$, $0 < r < R$, $0 < x < L$, and $g$ is even in $x$, then
$$g(x,r)=\sum_{k,l=0}^\infty i_{even}^{k,l}(x,r) G_{even}^{k,l}\mbox{.}$$

\item

$G_{odd}^{k,l}$ is well defined, and if $f \in \gS$, $g=Rf$, $0 < r < R$, $0 < x < L$, and $g$ is odd in $x$, then
$$g(x,r)=\sum_{k,l=0}^\infty i_{odd}^{k,l}(x,r) G_{odd}^{k,l}\mbox{.}$$

\end{enumerate}

\end{theo}

\begin{proof}

Let $f \in \gS$, $g=Rf$, $0 < r < R$, and $0 < x < L$.

\begin{enumerate}

\item

$G_{even}^{k,l}$ is obviously well defined since the integrals are finite. Assume without loss of generality that $g$ is even in $x$:
\begin{align*}
\sum_{k,l=0}^\infty& i_{even}^{k,l}(x,r) G_{even}^{k,l}\\
&=\chi_{(0,R)}(r) \chi_{(0,L)}(x) \sum_{k,l=0}^\infty \frac{\cos\(k\frac{\pi}{L}\sqrt{L^2-x^2}\,\)}{\sqrt{L^2-x^2}} \frac{\cos\(l\frac{\pi}{R}\sqrt{R^2-r^2}\,\)}{\sqrt{R^2-r^2}} \\
&\quad \times \frac{\epsilon_k\epsilon_l}{2LR} <i_{even}^{k,l}(x,r),g(x,r)+g(-x,r)>_{i_{even}}\\
&=\chi_{(0,R)}(r) \chi_{(0,L)}(x) \sum_{k,l=0}^\infty \frac{\cos\(k\frac{\pi}{L}\sqrt{L^2-x^2}\,\)}{\sqrt{L^2-x^2}} \frac{\cos\(l\frac{\pi}{R}\sqrt{R^2-r^2}\,\)}{\sqrt{R^2-r^2}}\\
&\quad\times \Biggl(\frac{\epsilon_k\epsilon_l}{2LR} \uint_0^L \uint_0^R \cos\(k\frac{\pi}{L}\sqrt{L^2-x'^2}\) \cos\(l\frac{\pi}{R}\sqrt{R^2-r'^2}\)\\
&\quad\times x' r' (g(x',r')+g(-x',r'))\, dx' \, dr'\Biggr)\\
&=\chi_{(0,R)}(r) \chi_{(0,L)}(x) \sum_{k,l=0}^\infty \frac{\cos\(k\frac{\pi}{L}\sqrt{L^2-x^2}\,\)}{\sqrt{L^2-x^2}} \frac{\cos\(l\frac{\pi}{R}\sqrt{R^2-r^2}\,\)}{\sqrt{R^2-r^2}}\\
&\quad\times \Biggl(\frac{\epsilon_k\epsilon_l}{LR} \uint_0^L \uint_0^R \cos\(k\frac{\pi}{L}\sqrt{L^2-x'^2}\) \cos\(l\frac{\pi}{R}\sqrt{R^2-r'^2}\)\\
&\quad\times x' r' g(x',r')\, dx' \, dr'\Biggr)
\end{align*}
With the substitutions $x''=\sqrt{L^2-x'^2}$ and $r''=\sqrt{R^2-r'^2}$ and $g_{k,l}$ representing the Fourier coefficients corresponding to $\frac{\epsilon_k\epsilon_l}{LR} x'' r'' g(\sqrt{L^2-x''^2},\sqrt{R^2-r''^2})$ it follows that
\begin{align*}
\sum_{k,l=0}^\infty& i_{even}^{k,l}(x,r) G_{even}^{k,l}\\
&=\chi_{(0,R)}(r) \chi_{(0,L)}(x) \sum_{k,l=0}^\infty \frac{\cos\(k\frac{\pi}{L}\sqrt{L^2-x^2}\,\)}{\sqrt{L^2-x^2}} \frac{\cos\(l\frac{\pi}{R}\sqrt{R^2-r^2}\,\)}{\sqrt{R^2-r^2}}\\
&\quad\times \Biggl(\frac{\epsilon_k\epsilon_l}{LR} \uint_0^L \uint_0^R \cos\(k\frac{\pi}{L}x''\) \cos\(l\frac{\pi}{R}r''\)
\end{align*}
\begin{align*}
&\quad\times x'' r'' g\left(\sqrt{L^2-x''^2},\sqrt{R^2-r''^2}\right)\, dx'' \, dr''\Biggr)\\
&=\chi_{(0,R)}(r) \chi_{(0,L)}(x) \frac{1}{\sqrt{L^2-x^2}} \frac{1}{\sqrt{R^2-r^2}} \\
&\quad \times \sum_{k,l=0}^\infty \cos\(k\frac{\pi}{L}\sqrt{L^2-x^2}\,\) \cos\(l\frac{\pi}{R}\sqrt{R^2-r^2}\,\) g_{k,l}\\
&=\chi_{(0,R)}(r) \chi_{(0,L)}(x) \frac{1}{\sqrt{L^2-x^2}} \frac{1}{\sqrt{R^2-r^2}} \\
&\quad \times \sqrt{L^2-x^2} \sqrt{R^2-r^2} g\left(\sqrt{L^2-\sqrt{L^2-x^2}^2},\sqrt{R^2-\sqrt{R^2-r^2}2}\right)\\
&=\chi_{(0,R)}(r) \chi_{(0,L)}(x) g(x,r)\mbox{.}
\end{align*}

\item

$G_{odd}^{k,l}$ is obviously well defined since the integrals are finite. Assume without loss of generality that $g$ is odd in $x$:
\begin{align*}
\sum_{k,l=0}^\infty& i_{odd}^{k,l}(x,r) G_{odd}^{k,l}\\
&=\chi_{(0,R)}(r) \chi_{(0,L)}(x) \sum_{k,l=0}^\infty x \frac{\cos\(k\frac{\pi}{L}\sqrt{L^2-x^2}\,\)}{\sqrt{L^2-x^2}} \frac{\cos\(l\frac{\pi}{R}\sqrt{R^2-r^2}\,\)}{\sqrt{R^2-r^2}} \\
&\quad \times \frac{\epsilon_k\epsilon_l}{2LR} <i_{odd}^{k,l}(x,r),g(x,r)-g(-x,r)>_{i_{odd}}\\
&=\chi_{(0,R)}(r) \chi_{(0,L)}(x) \sum_{k,l=0}^\infty x \frac{\cos\(k\frac{\pi}{L}\sqrt{L^2-x^2}\,\)}{\sqrt{L^2-x^2}} \frac{\cos\(l\frac{\pi}{R}\sqrt{R^2-r^2}\,\)}{\sqrt{R^2-r^2}}\\
&\quad\times \Biggl(\frac{\epsilon_k\epsilon_l}{2LR} \uint_0^L \uint_0^R \cos\(k\frac{\pi}{L}\sqrt{L^2-x'^2}\) \cos\(l\frac{\pi}{R}\sqrt{R^2-r'^2}\)\\
&\quad\times r' (g(x',r')-g(-x',r'))\, dx' \, dr'\Biggr)
\end{align*}
\begin{align*}
&=\chi_{(0,R)}(r) \chi_{(0,L)}(x) \sum_{k,l=0}^\infty x \frac{\cos\(k\frac{\pi}{L}\sqrt{L^2-x^2}\,\)}{\sqrt{L^2-x^2}} \frac{\cos\(l\frac{\pi}{R}\sqrt{R^2-r^2}\,\)}{\sqrt{R^2-r^2}}\\
&\quad\times \Biggl(\frac{\epsilon_k\epsilon_l}{LR} \uint_0^L \uint_0^R \cos\(k\frac{\pi}{L}\sqrt{L^2-x'^2}\) \cos\(l\frac{\pi}{R}\sqrt{R^2-r'^2}\)\\
&\quad\times r' g(x',r')\, dx' \, dr'\Biggr)
\end{align*}
With the substitutions $x''=\sqrt{L^2-x'^2}$ and $r''=\sqrt{R^2-r'^2}$ and $g_{k,l}$ representing the Fourier coefficients corresponding to $\frac{\epsilon_k\epsilon_l}{LR} \frac{x''}{\sqrt{L^2-x''^2}} r'' g(\sqrt{L^2-x''^2},\sqrt{R^2-r''^2})$ it follows that
\begin{align*}
\sum_{k,l=0}^\infty& i_{odd}^{k,l}(x,r) G_{odd}^{k,l}\\
&=\chi_{(0,R)}(r) \chi_{(0,L)}(x) \sum_{k,l=0}^\infty x \frac{\cos\(k\frac{\pi}{L}\sqrt{L^2-x^2}\,\)}{\sqrt{L^2-x^2}} \frac{\cos\(l\frac{\pi}{R}\sqrt{R^2-r^2}\,\)}{\sqrt{R^2-r^2}}\\
&\quad\times \Biggl(\frac{\epsilon_k\epsilon_l}{LR} \uint_0^L \uint_0^R \cos\(k\frac{\pi}{L}x''\) \cos\(l\frac{\pi}{R}r''\)\\
&\quad\times \frac{x''}{\sqrt{L^2-x''^2}} r'' g(\sqrt{L^2-x''^2},\sqrt{R^2-r''^2})\, dx'' \, dr''\Biggr)\\
&=\chi_{(0,R)}(r) \chi_{(0,L)}(x) x \frac{1}{\sqrt{L^2-x^2}} \frac{1}{\sqrt{R^2-r^2}} \\
&\quad \times \sum_{k,l=0}^\infty \cos\(k\frac{\pi}{L}\sqrt{L^2-x^2}\,\) \cos\(l\frac{\pi}{R}\sqrt{R^2-r^2}\,\) g_{k,l}\\
&=\chi_{(0,R)}(r) \chi_{(0,L)}(x) x \frac{1}{\sqrt{L^2-x^2}} \frac{1}{\sqrt{R^2-r^2}} \\
&\quad \times \frac{\sqrt{L^2-x^2}}{\sqrt{L^2-\sqrt{L^2-x^2}^2}} \sqrt{R^2-r^2} \\
&\quad \times g(\sqrt{L^2-\sqrt{L^2-x^2}^2},\sqrt{R^2-\sqrt{R^2-r^2}2})\\
&=\chi_{(0,R)}(r) \chi_{(0,L)}(x) g(x,r)\mbox{.}
\end{align*}

\end{enumerate}

\end{proof}

\section{Reconstruction of the orthogonal functions}

As seen above, the measurable data can be projected onto the set of orthogonal functions $\{i_{even}^{k,l}, i_{odd}^{k,l}: k, l \in \N_0\}$ and can be recovered again. In the following it is therefore sufficient to perform the reconstruction only for the orthogonal functions.

\begin{theo}

\label{reconstructions}

\begin{enumerate}

\item

Let $k, l \in \N_0$ and
$$g(x,r)=i_{even}^{k,l}(x,r)\mbox{.}$$
Then
$$\hat{f}(\xi,\eta)=\sqrt{\frac{\pi}{8}} \: |\eta|\bn\left(L\sqrt{(k\frac{\pi}{L})^2+\xi^2}\,\right)\frac{\sin\(R\sqrt{(l\frac{\pi}{R})^2+\xi^2+\eta^2}\,\)}{\sqrt{(l\frac{\pi}{R})^2+\xi^2+\eta^2}}$$
and
\begin{align*}
f(x,y)&=\sqrt{\frac{\pi}{2}} \: \mbox{H}_y \frac{\pa}{\pa y} \uint_\R \left( \left\{ \begin{array}{lr}
\frac{\cos\(k\frac{\pi}{L}\sqrt{L^2-t^2}\,\)}{\sqrt{L^2-t^2}} &
\mbox{ for } |t|<L\\
& \\
0 & \mbox{ otherwise}\\
\end{array} \right\}\right.\\
&\quad\times \left. \left\{ \begin{array}{lr}
\frac{\cos\(l\frac{\pi}{R}\sqrt{R^2-(x-t)^2-y^2}\,\)}{\sqrt{R^2-(x-t)^2-y^2}} &
\mbox{ for } (x-t)^2+y^2<R^2\\
& \\
0 & \mbox{ otherwise}\\
\end{array} \right\} \right) dt
\end{align*}
and $\mbox{H}_y$ refers again to the Hilbert transform in $y$.

\item

Let $k, l \in \N_0$ and
$$g(x,r)=i_{odd}^{k,l}(x,r)\mbox{.}$$
Then
$$\hat{f}(\xi,\eta)=-\sqrt{\frac{\pi}{8}} \: |\eta|\frac{\pa}{\pa \xi}\bn\left(L\sqrt{(k\frac{\pi}{L})^2+\xi^2}\,\right)\frac{\sin\(R\sqrt{(l\frac{\pi}{R})^2+\xi^2+\eta^2}\,\)}{\sqrt{(l\frac{\pi}{R})^2+\xi^2+\eta^2}}$$
and
\begin{align*}
f(x,y)&=\sqrt{\frac{\pi}{2}} \: x \mbox{H}_y \frac{\pa}{\pa y} \uint_\R \left( \left\{ \begin{array}{lr}
\frac{\cos\(k\frac{\pi}{L}\sqrt{L^2-t^2}\,\)}{\sqrt{L^2-t^2}} &
\mbox{ for } |t|<L\\
& \\
0 & \mbox{ otherwise}\\
\end{array} \right\}\right.\\
&\quad\times \left. \left\{ \begin{array}{lr}
\frac{\cos\(l\frac{\pi}{R}\sqrt{R^2-(x-t)^2-y^2}\,\)}{\sqrt{R^2-(x-t)^2-y^2}} &
\mbox{ for } (x-t)^2+y^2<R^2\\
& \\
0 & \mbox{ otherwise}\\
\end{array} \right\} \right) dt\mbox{.}
\end{align*}

\end{enumerate}

\end{theo}

\begin{proof}

\begin{enumerate}

\item

Let $k, l \in \N_0$ and
$$g(x,r)=i_{even}^{k,l}(x,r)\mbox{.}$$
Then the three dimensional Fourier transform of $g$ is given by
$$\hat{g}(\xi,\rho)=\frac{1}{\sqrt{2\pi}}\uint_\R \uint_0^\infty \exp^{-ix\xi}r\bn(r\rho)g(x,r)\,dr \, dx$$
and with lemmas \ref{cosinusendlich} and \ref{hankelendlich}
$$\hat{g}(\xi,\rho)=\sqrt{\frac{\pi}{2}} \; \bn\left(L\sqrt{(k\frac{\pi}{L})^2+\xi^2}\,\right)\frac{\sin\(R\sqrt{(l\frac{\pi}{R})^2+\rho^2}\,\)}{\sqrt{(l\frac{\pi}{R})^2+\rho^2}}\mbox{.}$$
With theorem \ref{inversionformula} this results in
\begin{align*}
\hat{f}(\xi,\eta)&=\frac{1}{2}|\eta|\hat{g}\(\xi,\sqrt{\xi^2+\eta^2}\)\\
&=\sqrt{\frac{\pi}{8}} \: |\eta| \bn\left(L\sqrt{(k\frac{\pi}{L})^2+\xi^2}\,\right)\frac{\sin\(R\sqrt{(l\frac{\pi}{R})^2+\xi^2+\eta^2}\,\)}{\sqrt{(l\frac{\pi}{R})^2+\xi^2+\eta^2}}\mbox{.}
\end{align*}
Therefore
$$f(x,y)=\frac{1}{2\pi} \uint_\R \uint_\R \exp^{ix\xi} \exp^{iy\eta} \hat{f}(\xi,\eta)\,d\eta \, d\xi$$
and applying the Fourier convolution theorem leads to
\begin{align*}
f(x,y)&=\frac{1}{2\pi}\sqrt{\frac{\pi}{8}} \: \mbox{H}_y \frac{\pa}{\pa y} \uint_\R \left( \uint_\R \exp^{it\xi} \bn\left(L\sqrt{(k\frac{\pi}{L})^2+\xi^2}\,\right)\,d\xi \right.\\
&\quad\times \left. \uint_\R \uint_\R \exp^{i(x-t)\xi} \exp^{iy\eta} \frac{\sin\(R\sqrt{(l\frac{\pi}{R})^2+\xi^2+\eta^2}\,\)}{\sqrt{(l\frac{\pi}{R})^2+\xi^2+\eta^2}}\,d\xi \, d\eta \right) dt\\
&=\sqrt{\frac{\pi}{2}} \: \mbox{H}_y \frac{\pa}{\pa y} \uint_\R \left( \uint_0^\infty \cos(t\xi) \bn\left(L\sqrt{(k\frac{\pi}{L})^2+\xi^2}\,\right)\right.\,d\xi\\
&\quad\times \left. \uint_0^\infty \bn\(\rho\sqrt{(x-t)^2+y^2}\) \frac{\sin\(R\sqrt{(l\frac{\pi}{R})^2+\rho^2}\,\)}{\sqrt{(l\frac{\pi}{R})^2+\rho^2}}\,d\rho \right) dt\mbox{.}
\end{align*}
With lemmas \ref{cosinusendlich} and \ref{hankelendlich} and because the Hankel transform and the cosine transform are their respective inverses, this can be written as:
\begin{align*}
f(x,y)&=\sqrt{\frac{\pi}{2}} \: \mbox{H}_y \frac{\pa}{\pa y} \uint_\R \left( \left\{ \begin{array}{lr}
\frac{\cos\(k\frac{\pi}{L}\sqrt{L^2-t^2}\,\)}{\sqrt{L^2-t^2}} &
\mbox{ for } |t|<L\\
& \\
0 & \mbox{ otherwise}\\
\end{array} \right\}\right.\\
&\quad\times \left. \left\{ \begin{array}{lr}
\frac{\cos\(l\frac{\pi}{R}\sqrt{R^2-(x-t)^2-y^2}\,\)}{\sqrt{R^2-(x-t)^2-y^2}} &
\mbox{ for } (x-t)^2+y^2<R^2\\
& \\
0 & \mbox{ otherwise}\\
\end{array} \right\} \right) dt\mbox{.}
\end{align*}

\item

Let $k, l \in \N_0$ and
$$g(x,r)=i_{odd}^{k,l}(x,r)\mbox{.}$$
Then the proof is analogous to 1. using $x\sin(x\xi)=-\frac{\pa}{\pa \xi}\cos(x\xi)$.

\end{enumerate}

\end{proof}

\begin{rem}

$G_{even}^{k,l}$ and $G_{odd}^{k,l}$ can be computed via a fast two dimensional non-equidistant cosine transform as described in \cite{Fourmont}. As according to theorem \ref{grecoverable}
$$g(x,r)=\sum_{k,l=0}^\infty (1+x) \frac{\cos\(k\frac{\pi}{L}\sqrt{L^2-x^2}\,\)}{\sqrt{L^2-x^2}} \frac{\cos\(l\frac{\pi}{R}\sqrt{R^2-r^2}\,\)}{\sqrt{R^2-r^2}} (G_{even}^{k,l}+G_{odd}^{k,l})\mbox{,}$$
it follows with the results of theorem \ref{reconstructions} that
\begin{align*}
f(x,y)&=\sum_{k,l=0}^\infty (G_{even}^{k,l}+G_{odd}^{k,l})\\
&\quad\times \sqrt{\frac{\pi}{2}} \: (1+x) \mbox{H}_y \frac{\pa}{\pa y} \uint_\R \left( \left\{ \begin{array}{lr}
\frac{\cos\(k\frac{\pi}{L}\sqrt{L^2-t^2}\,\)}{\sqrt{L^2-t^2}} &
\mbox{ for } |t|<L\\
& \\
0 & \mbox{ otherwise}\\
\end{array} \right\}\right.\\
&\quad\times \left. \left\{ \begin{array}{lr}
\frac{\cos\(l\frac{\pi}{R}\sqrt{R^2-(x-t)^2-y^2}\,\)}{\sqrt{R^2-(x-t)^2-y^2}} &
\mbox{ for } (x-t)^2+y^2<R^2\\
& \\
0 & \mbox{ otherwise}\\
\end{array} \right\} \right) dt
\mbox{.}
\end{align*}
This summation of the precomputed solutions of the orthogonal functions with the appropriate weights $G_{even}^{k,l}$ and $G_{odd}^{k,l}$ is not faster than the direct implementation of the inversion formulas presented so far. However two improvements are possible.

\begin{enumerate}

\item

This inversion formula is adapted to the fact that the measured data is limited. Therefore it is probably easier to identify means to diminish the artifacts caused by this limitation.

\item

Depending on the application, the coefficients representing higher frequencies are probably prone to noise. Therefore it could be worthwhile to neglect them in the summation to obtain a faster algorithm since these coefficients do not contain much dependable information.

\end{enumerate}

\end{rem}

\chapter{Tackling the left-right-ambiguity}

\label{Antisym}

The so called problem of the left-right ambiguity is one that exists since the invention of SAR and continues to be an issue, albeit alleviated, till this day \cite{Ulander2}. It refers to the fact that because of geometrical reasons the reconstruction formulas in this field are not able to distinguish between points that lie symmetrically to the left and right of the flight path. At first this problem seemed inevitable \cite{Fawcett}, \cite{Andersson}, but later, approaches to solve it were proposed. The method of beamforming \cite{Cheney}, for example, tries to overcome the ambiguity by directing the energy emitted from the antenna as much toward one side of the flight track as possible. In this way the echo from the illuminated side is much stronger than the echo from the other side. Unfortunately these weak signals still lead to shadows in the reconstruction for highly reflecting objects. A different solution was proposed in \cite{CheneyNatterer}. The idea consisted of using data from two parallel flight tracks. However, singularities arising from noisy measurements present a serious problem in this approach.\\
Inspired by the idea of gathering and then combining data from two slightly different positions, two post-processing formulas are given in the following that manipulate the data measured by an airplane equipped with two or more antennas. It is shown that it is possible to recover the odd part of the reflectivity function from measurements with at least two antennas. The reconstruction of the original image, including the odd part, is difficult for only two antennas, but it will be shown that since the low frequencies are attenuated in the inversions of the spherical Radon transform, it is possible. As will be seen, for a good reconstruction quality it is nevertheless recommendable to use more than two antennas. This has the additional advantage that no regularization is needed and therefore an exact formula can be used.\\
The chapter is structured as follows. The first section will present analytic considerations, which prove that it is possible to retrieve the original, asymmetric image given at least two even images. This corresponds to a measurement with at least two antennas. In the next section the two post-processing formulas are introduced. The first requires a sophisticated regularization whereas the second needs at least three source-images as input. Afterwards the results of numerical simulations are displayed.

\section{Using two data sets}

At first some definitions are necessary. Then the important trick is shown separately, which makes the following central theorem fairly easy. At the end some conclusions derived from the theorem are mentioned.

\begin{defi}

Let $f \in \gS(\R^n \times \R)$ and $b \in \R$. Then

\begin{enumerate}

\item

$f^{F,I}$ and $f^{I,F}$ denote the Fourier transform of $f$ in the first $n$ variables and the last variable respectively. $f^{I,I}=f$ and $f^{F,F}=\hat{f}$.

\item

$f^{(C,I)}(x,y)=f(-x,y)$ and $f^{(I,C)}(x,y)=f(x,-y)$.

\item

$f^e_b(x,y)=\frac{1}{2}[f(x,y+b)+f(x,-y+b)]$

\vspace{.1cm}\hspace{1.5cm}$=\frac{1}{2}[(\tau_{(0,\:b)}f)(x,y)+(\tau_{(0,\:b)}f)^{(I,C)}(x,y)]$.

\end{enumerate}

\end{defi}

\begin{lem}

\label{trick}

Let $f \in \gS(\R^n \times \R)$ and $b \in \R$. Then $(\tau_{(0,\:b)}f)^{(I,C)}(x,y)=\tau_{(0,\:-b)}f^{(I,C)}$.

\end{lem}

\begin{proof}

\begin{align*}
(\tau_{(0,\:b)}f)^{(I,C)}(x,y)&=(\tau_{(0,\:b)}f)(x,-y)=f(x,-y+b)\\
&=f^{(I,C)}(x,y-b)=\tau_{(0,\:-b)}f^{(I,C)}(x,y)\mbox{.}
\end{align*}

\end{proof}

Inspired by the ideas in \cite{CheneyNatterer}, the following theorem is formulated.

\begin{theo}

Let $f \in \gS(\R^n \times \R)$, $b, \eta \in \R$, and $x \in \R^{n}$. Then
$$\sin(b\eta)f^{I,F}(x,\eta)=\frac{1}{i}[f^e_b-\tau_{(0,\:-b)}f^e_0]^{I,F}(x,\eta)\mbox{.}$$

\end{theo}

\begin{proof}

\begin{align*}
\sin(b\eta)f^{I,F}(x,\eta)&=\frac{1}{2i}[\exp^{ib\eta}f^{I,F}(x,\eta)-\exp^{-ib\eta}f^{I,F}(x,\eta)]\\
&=\frac{1}{2i}[\tau_{(0,\:b)}f-\tau_{(0,\:-b)}f]^{I,F}(x,\eta)
\end{align*}
$$=\frac{1}{2i}[\tau_{(0,\:b)}f-\tau_{(0,\:-b)}f+(\tau_{(0,\:b)}f)^{(I,C)}-(\tau_{(0,\:b)}f)^{(I,C)}]^{I,F}(x,\eta)\mbox{.}$$
Using lemma \ref{trick},
\begin{align*}
\sin(b\eta)&f^{I,F}(x,\eta)\\
&=\frac{1}{2i}[\tau_{(0,\:b)}f+(\tau_{(0,\:b)}f)^{(I,C)}-\tau_{(0,\:-b)}f-(\tau_{(0,\:-b)}f^{(I,C)})]^{I,F}(x,\eta)\\
&=\frac{1}{2i}[2f^e_b-\tau_{(0,\:-b)}(f+f^{(I,C)})]^{I,F}(x,\eta)=\frac{1}{i}[f^e_b-\tau_{(0,\:-b)}f^e_0]^{I,F}(x,\eta)\mbox{.}
\end{align*}

\end{proof}

This means that the data from two antennas with distance $b$ is sufficient to recover $\sin(b\eta)f^{I,F}$.

\begin{rem}

\begin{enumerate}

\item

$f^{I,F}(x,\eta)$ cannot be recovered with two antennas with distance $b$ for \mbox{$\eta \in \frac{\pi}{b} \Z$}.

\item

$f^{I,F}(x,0)=(f^e_b)^{I,F}(x,0) \, \forall b \in \R$.



\end{enumerate}

\end{rem}

\begin{cor}

Let $f \in \gS(\R^n \times \R)$. If $f^e_b$ is known for all $b \geq 0$, then $f$ is completely determined.

\end{cor}

\begin{proof}

If $f^e_0$ is given, the even part of $f$ is known.
The odd part is determined for $b \geq 0$ by
\begin{align*}
\frac{1}{2}[f(x,b)-f(x,-b)]&=\frac{1}{2\pi} \uint_\R \sin(b\eta) \uint_\R \sin(y\eta) f(x,y) \, dy \, d\eta\\
&=\frac{-i}{\sqrt{2\pi}} \uint_\R \sin(b\eta)f^{I,F}(x,\eta) \, d\eta \\
&=\frac{-1}{\sqrt{2\pi}} \uint_\R [f^e_b-\tau_{(0,\:-b)}f^e_0]^{I,F}(x,\eta) \, d\eta \mbox{.}
\end{align*}

\end{proof}

\section{Post-processing formulas}

\begin{defi}

Let $f \in \gS(\R^n \times \R)$. Define
$$h_b (\eta):= \sin(b\eta) f^{I,F}(x,\eta)\mbox{.}$$

\end{defi}

It turns out that a noise resistant algorithm to recover the odd part of the reflectivity function $f$ in addition to the even part is not so easily found due to the many zeros of the sine. Some obvious regularizations such as
$$f^{I,F}(x,\eta) \approx \left\{ \begin{array}{lr}
(f^e_{a})^{I,F}(x,0) &
\mbox{ for }\eta=0 \mbox{ with } a \in \{0,b\}\\
& \\
\frac{\sin(b\eta)h_b(\eta)}{\sin^2(b\eta) + \epsilon} & \mbox{ otherwise,}\\
\end{array} \right.$$
$$f^{I,F}(x,\eta) \approx \left\{ \begin{array}{lr}
(f^e_{a})^{I,F}(x,0) &
\mbox{ for }\eta=0 \mbox{ with } a \in \{0,b\}\\
& \\
\frac{\sin(b\eta)h_b(\eta)}{\sin^2(b\eta) + \epsilon\cos^{2k}(b\eta)} & \mbox{ otherwise}\\
\end{array} \right.$$
or
$$f^{I,F}(x,\eta) \approx \left\{ \begin{array}{lr}
(f^e_{a})^{I,F}(x,0) &
\mbox{ for }\eta=0 \mbox{ with } a \in \{0,b\}\\
& \\
\frac{\sin(b\eta)h_b(\eta)}{\sin^2(b\eta) + \epsilon\cos^{2k}(b\eta)\mbox{\footnotesize exp}(d\eta^{2l})} & \mbox{ otherwise}\\
\end{array} \right.$$
with $k, l \in \N$ and $\epsilon > 0$ yield only high-pass-filtered images. The regularization
$$f^{I,F}(x,\eta) \approx \left\{ \begin{array}{lr}
(f^e_{a})^{I,F}(x,0) &
\mbox{ for }\eta=0 \mbox{ with } a \in \{0,b\} ,\\
& \\
\frac{\sin(b\eta)h_b(\eta)}{\sin^2(b\eta) + \epsilon\cos^{2k}(b\eta)(\mbox{\footnotesize exp}(d\eta^{2l})-1)} & \mbox{ otherwise}\\
\end{array} \right.$$
did not provide satisfactory results either.
To avoid the problem of a vanishing denominator, two satisfactory solutions exist which will be given in the following theorem.

\begin{theo}

\label{antisymformel}

Let $f \in \gS(\R^n \times \R)$.

\begin{enumerate}

\item

For two data sets, i.e. only one $b \in \R$, a regularization is necessary:
$$f^{I,F}(x,\eta) \approx \left\{ \begin{array}{lr}
(f^e_{a})^{I,F}(x,0) &
\mbox{ for }\eta=0 \mbox{ with } a \in \{0,b\}\\
& \\
\frac{\sin(b\eta)h_b(\eta)}{\sin^2(b\eta) + H_\frac{\pi}{2b}(\eta)\epsilon\cos^{2k}(b\eta)} & \mbox{ otherwise}\\
\end{array} \right.$$
with $\epsilon >0$, $k \in \N$, and the Heaveside function $H=\chi_{(\frac{\pi}{2b},\infty)}$.

\item

For $m>2$ data sets, i.e. a set $b_1, ..., b_{m(m-1)/2} \in \R$, it is possible to choose the $b_k$ in a way that the denominator does not vanish, e.g. $b_1=1, b_2=\pi$. Therefore an analytically exact formula is achievable:
$$f^{I,F}(x,\eta)=\left\{ \begin{array}{lr}
(f^e_{b_k})^{I,F}(x,0) &
\mbox{ for }\eta=0 \mbox{ with}\hspace{1.8cm}\\
& k \in \{1, ..., m(m-1)/2\}\\
& \\
\frac{\usum_{k=1}^{m(m-1)/2}\sin(b_k\eta)h_{b_k}(\eta)}{\usum_{k=1}^{m(m-1)/2}\sin^2(b_k\eta)} & \mbox{ otherwise.}\\
\end{array} \right.$$

\end{enumerate}

\end{theo}

\begin{proof}


\begin{enumerate}

\item

Let $\eta \not\in \frac{\pi}{2b} \Z$, $k \in \N$, and $\epsilon >0$. Then
\begin{align*}
&\frac{\sin(b\eta)h_b(\eta)}{\sin^2(b\eta) + H_\frac{\pi}{2b}(\eta)\epsilon\cos^{2k}(b\eta)}\\
&=\frac{\sin^2(b\eta)f^{I,F}(x,\eta)}{\sin^2(b\eta) + H_\frac{\pi}{2b}(\eta)\epsilon\cos^{2k}(b\eta)} \geg{\epsilon}{0} f^{I,F}(x,\eta)\mbox{.}
\end{align*}

\item

Let $\eta \neq 0$. Then
\begin{align*}
&\frac{\usum_{k=1}^{m(m-1)/2}\sin(b_k\eta)h_{b_k}(\eta)}{\usum_{k=1}^{m(m-1)/2}\sin^2(b_k\eta)}\\
&=\frac{\usum_{k=1}^{m(m-1)/2}\sin^2(b_k\eta)f^{I,F}(x,\eta)}{\usum_{k=1}^{m(m-1)/2}\sin^2(b_k\eta)}=f^{I,F}(x,\eta)\mbox{.}
\end{align*}

\end{enumerate}

\end{proof}

\begin{rem}

\begin{enumerate}

\item

The same as in theorem \ref{antisymformel} could be shown for the other aforementioned regularizations although they do not yield satisfactory results. However this property is necessary to be a viable approach.

\item

It is essential not to regularize everywhere. Otherwise only a high-pass filtered image is obtained.



\item

Considering the sampling-theorem, the formula for two data sets is exact if supp $h_b \subset (-\frac{\pi}{b},-\frac{\pi}{b})$, i.e. if the bandwidth of $f$ in $y$-direction is smaller than $b$.

\item

An advantage of using more than two data sets is a sizable reduction of noise-effects by the weighted summation of the data.

\end{enumerate}

\end{rem}

\section{Numerical simulations}

In the following, two different kinds of numerical simulations will be shown. First, to demonstrate the ability of the proposed algorithms in solving the left-right ambiguity, a set of left-right symmetric images will be created that serve as input for the algorithms. In this way artifacts arising from the inversion of the spherical Radon transform are omitted, and it is possible to isolate the effects of the post-processing formulas. Second, to show what results can be expected in reality, measurements with multiple antenna positions will be simulated with the spherical Radon transform that serve as data for reconstructions according to theorem \ref{modinvform}. Then the acquired left-right symmetric images are fed into the algorithms to solve the left-right ambiguity.\\
In the following figures antenna distances will be given in multiples of the sampling length (pixels). The noise is applied as a pixelwise multiplication with the corresponding percentage multiplied with a completely uncorrelated (with regard to pixels as well as pictures), normally distributed deviate with zero mean and unit variance. Additionally a similarly structured absolute error is applied.\\
Figure \ref{reffunction} shows on the left hand side a simple circle with five smaller circles arranged in a cross formation in it. The flight track, i. e. the leftmost antenna,  runs in a vertical line right in the middle of the figure with the additional antennas on the right. In the following this will serve as the phantom for the attempts to reconstruct the reflectivity function $f$. The right hand side of figure \ref{reffunction} shows the corresponding cross section.
\asbilds{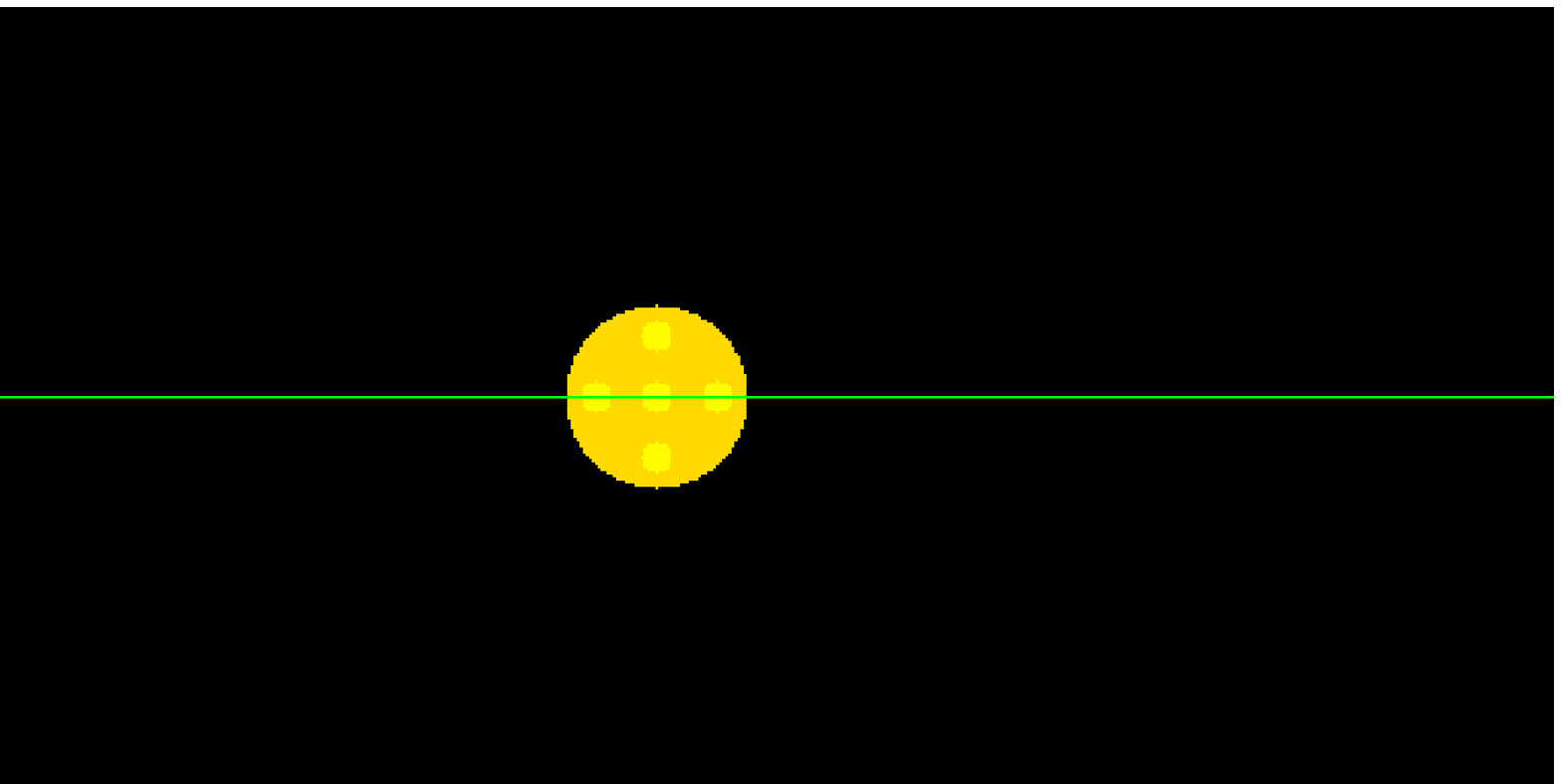}{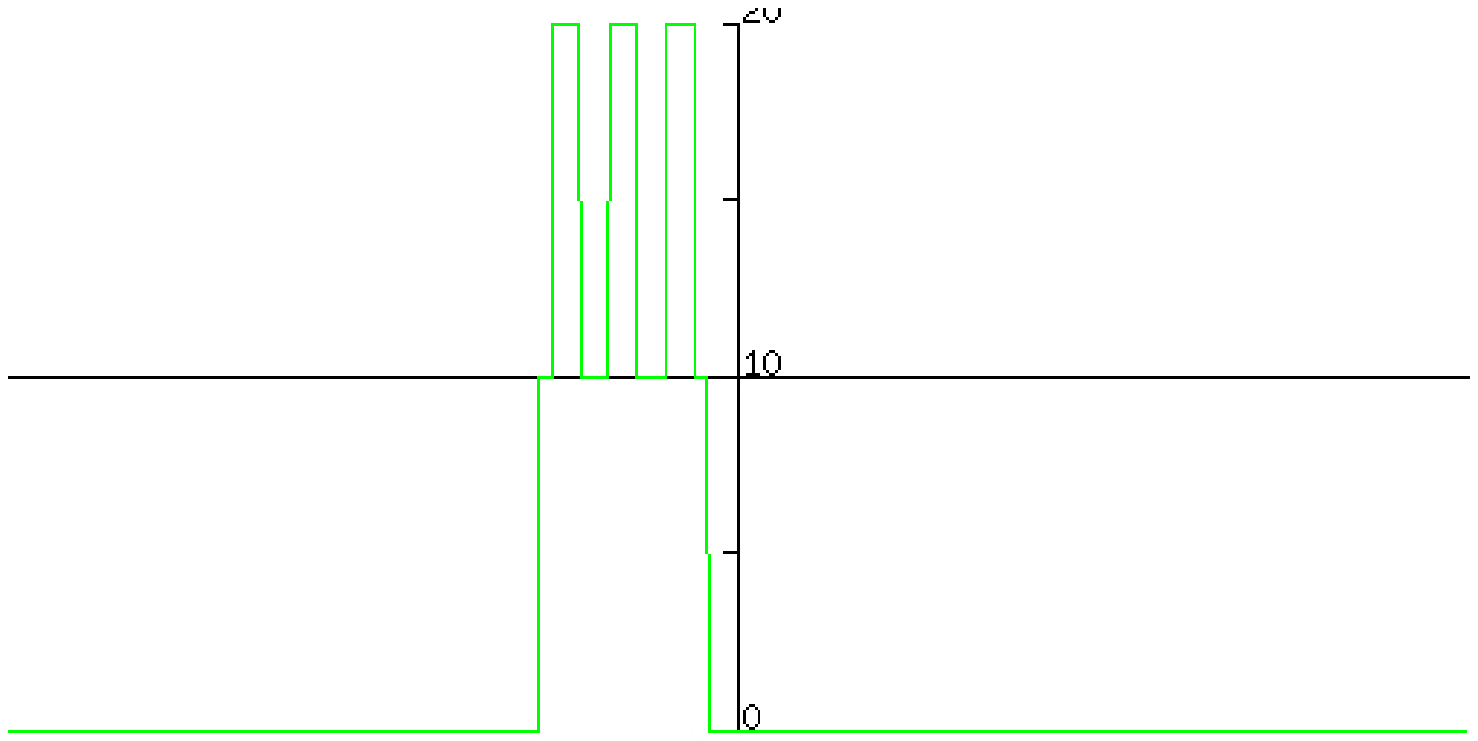}{Reflectivity function $f$\label{reffunction}}

\subsection{Solving the left-right ambiguity with computed symmetric images\label{Typ2}}

\label{symmetric}

In the following examples the even part of $f$ with respect to several flight tracks is computed directly. These even images are combined to reconstruct the ground reflectivity function $f$ according to the corresponding formulas in theorem \ref{antisymformel}. As an example, figure \ref{measuredtype2} shows the symmetric input data, where the number of circle pairs depends on the number of antennas employed and the distance between the two circles of the pair is associated with the antenna distance. The noise level is reflected in the intensity of the deviations and will be kept at a level of $10 \%$ in this subsection.

\asbildhg{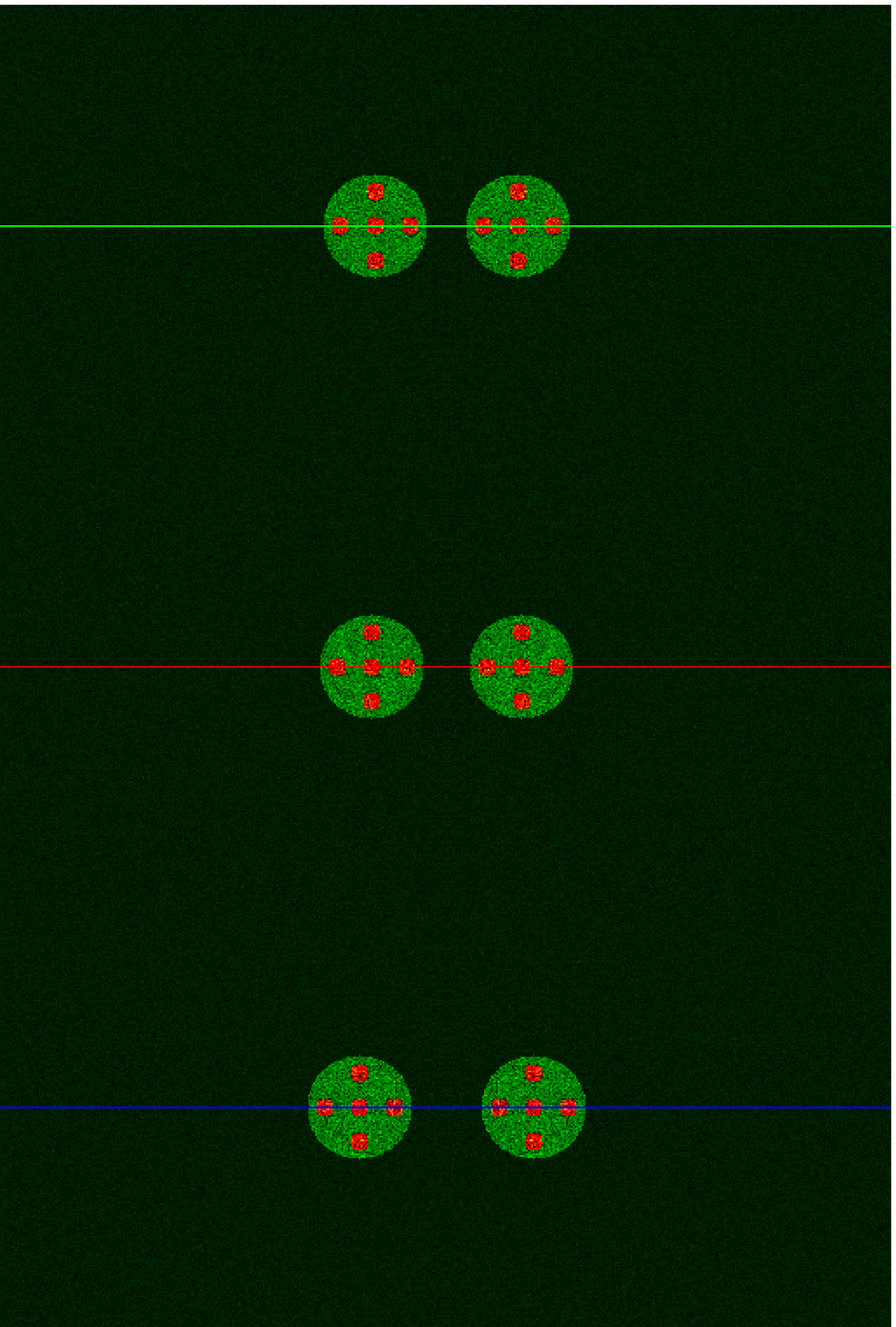}{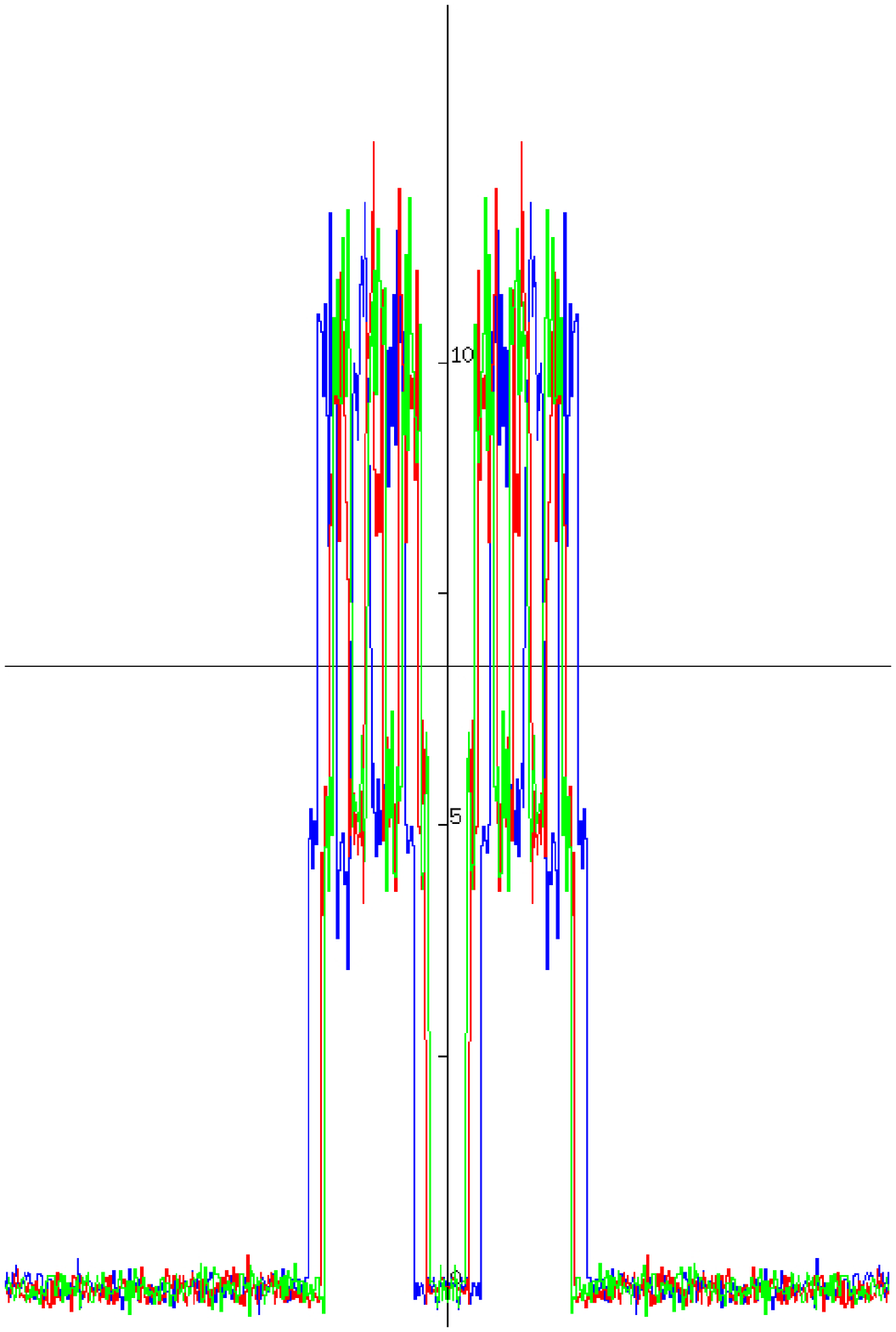}{The even part of $f$ with respect to three different flight tracks with 10\% noise\label{measuredtype2}}

First, the influence of the antenna distance is demonstrated. Figures \ref{2Antennen1Typ2} to \ref{2Antennen9Typ2FS} contrast the achievable results for reconstructions with a small antenna distance to reconstructions with a large antenna distance.

\asbildz{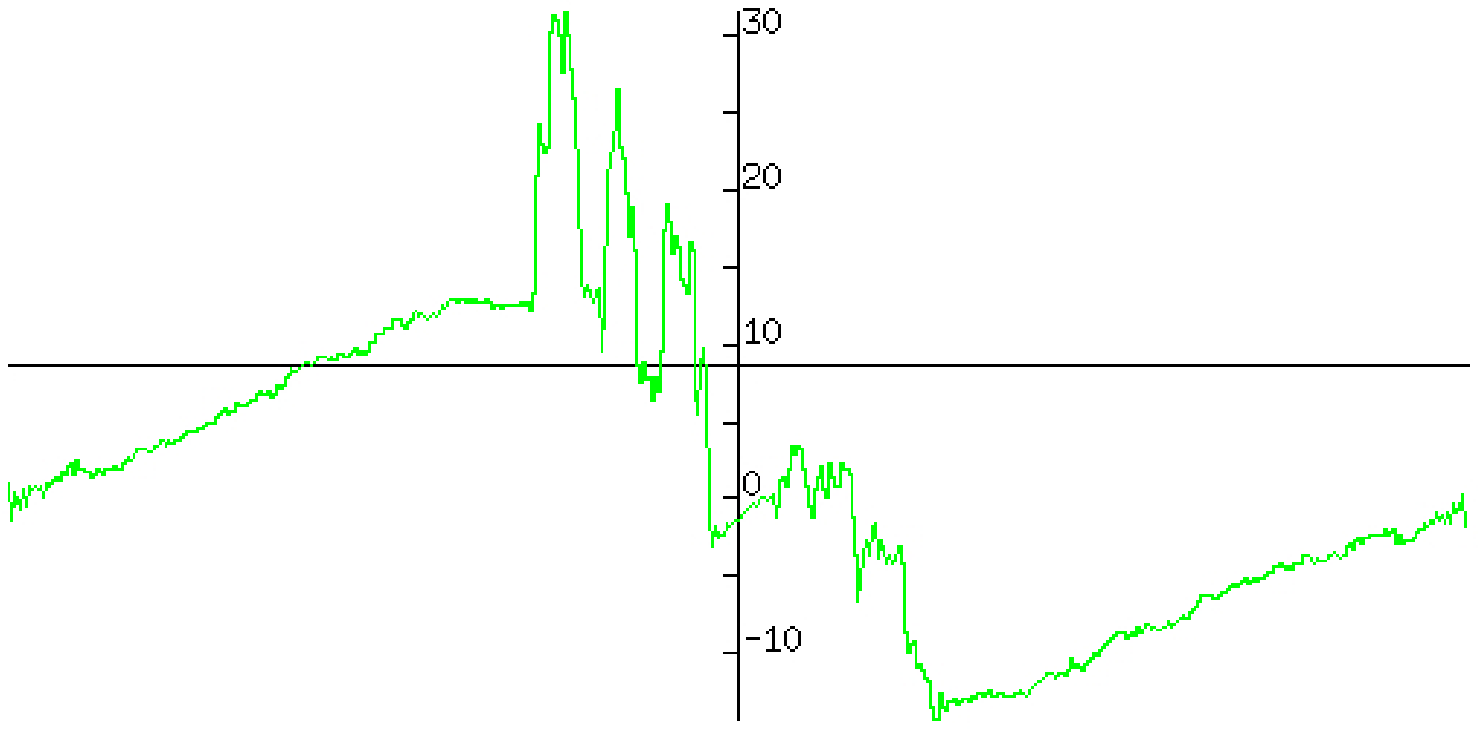}{Two antennas with distance 1, 10\% noise, spatial domain\label{2Antennen1Typ2}}{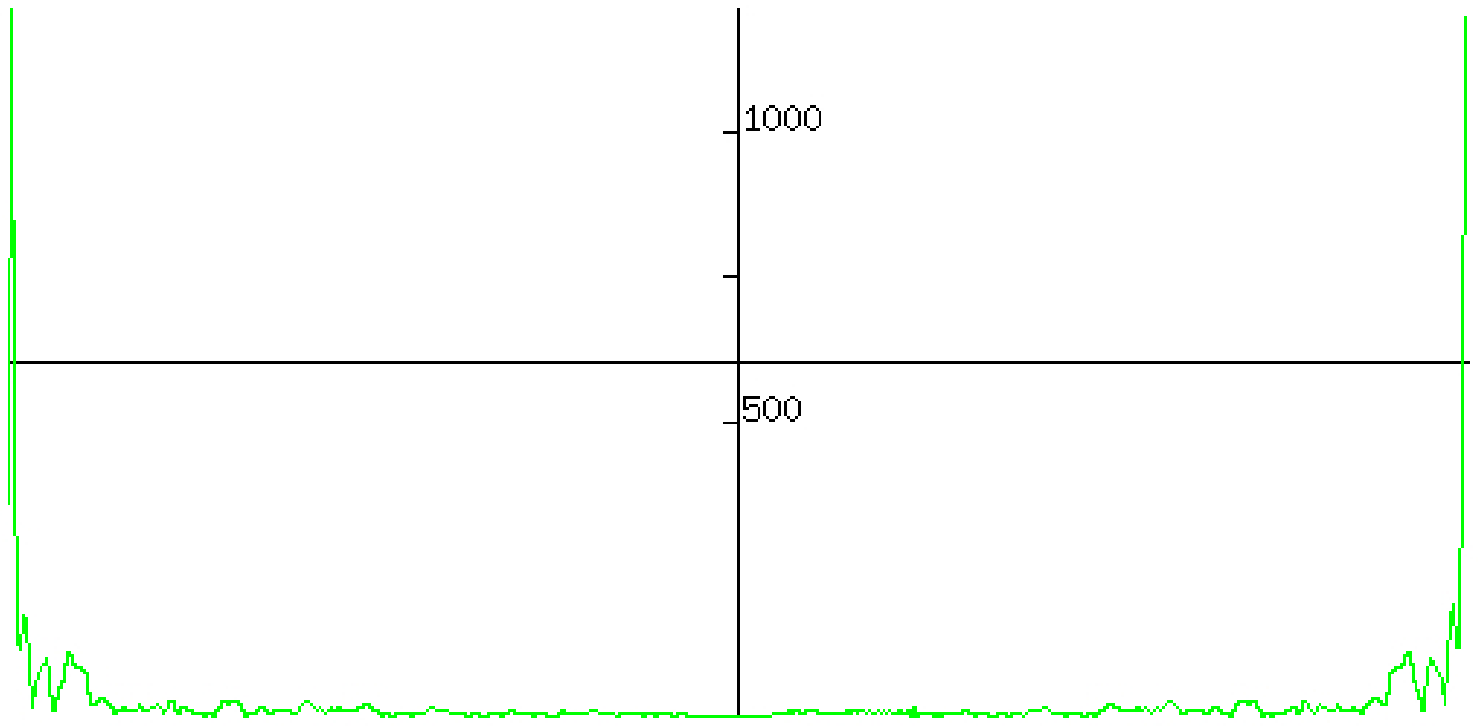}{Two antennas with distance 1, 10\% noise, frequency domain. Left: Logarithmic scale\label{2Antennen1Typ2FS}}
\asbildz{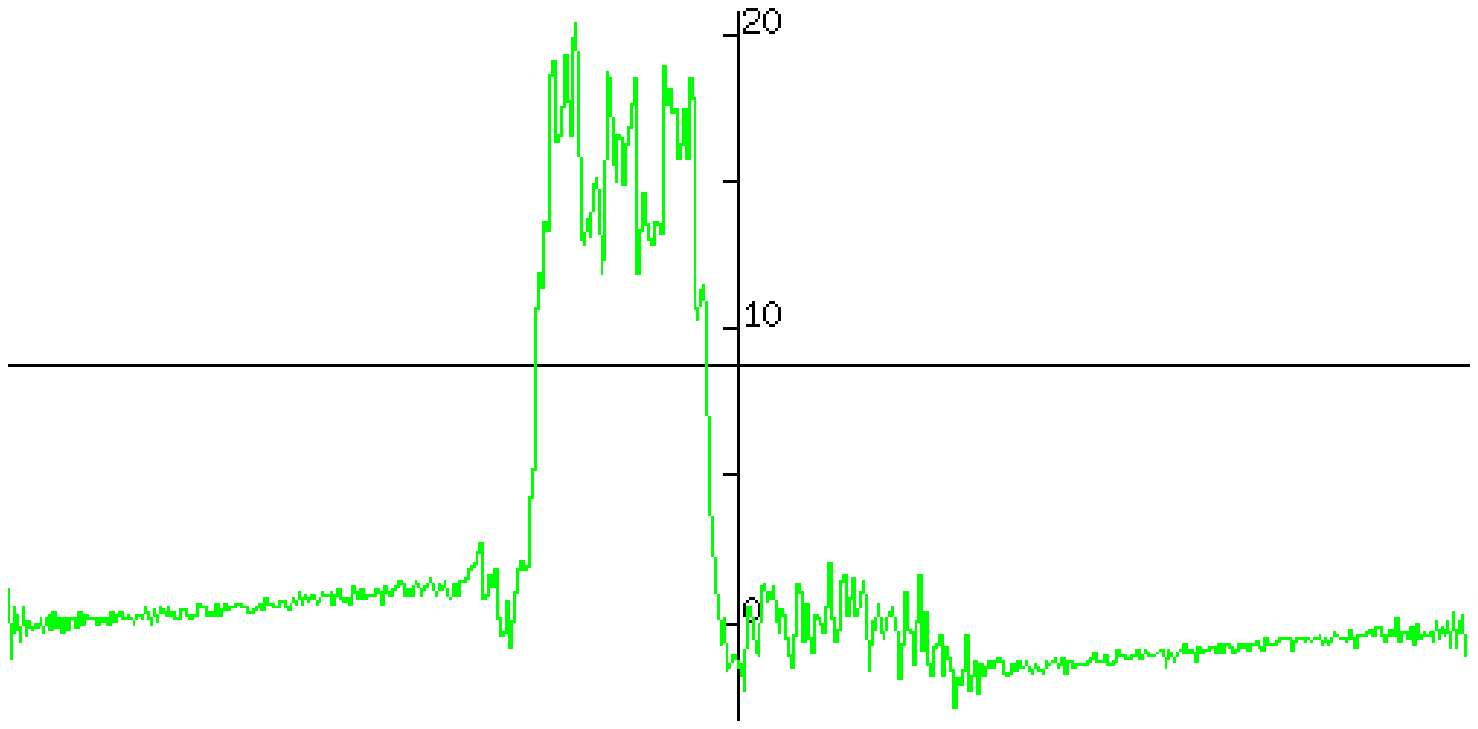}{Two antennas with distance 9, 10\% noise, spatial domain\label{2Antennen9Typ2}}{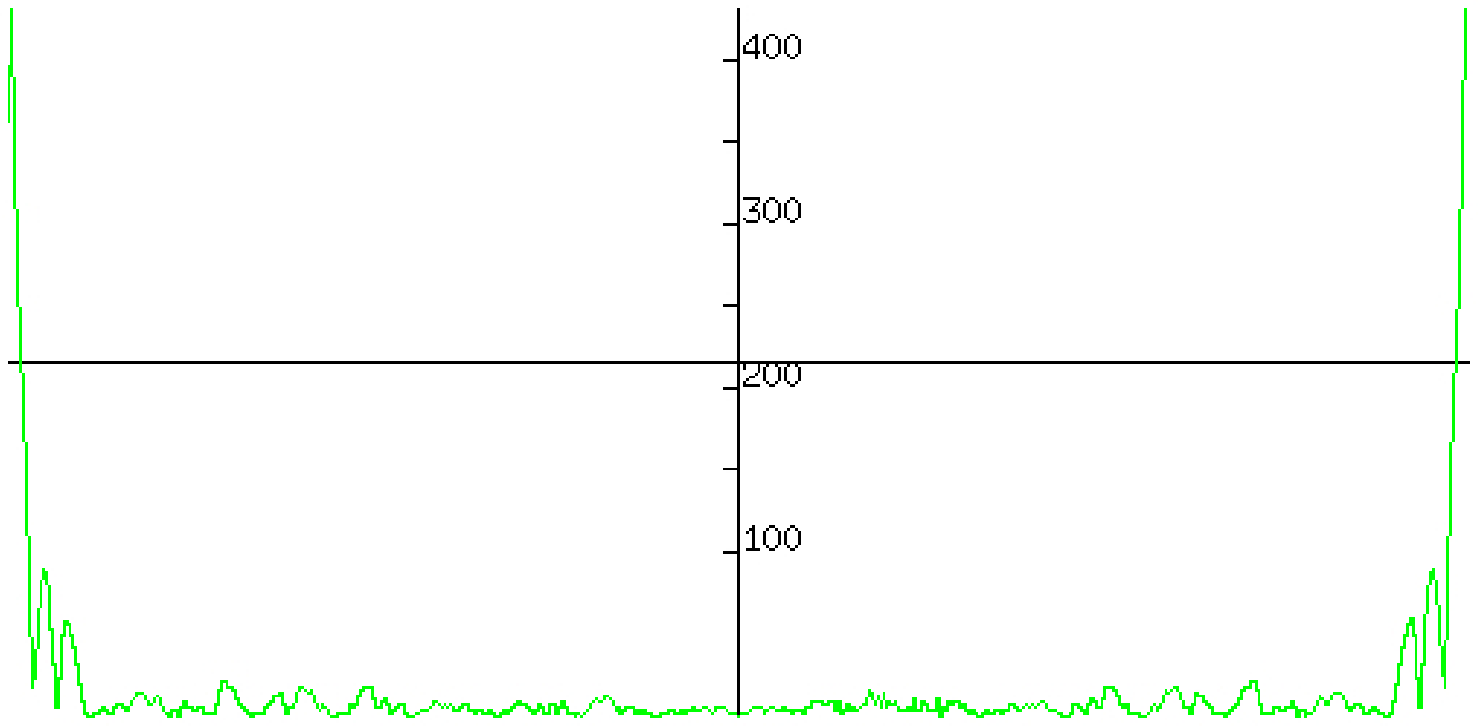}{Two antennas with distance 9, 10\% noise, frequency domain. Left: Logarithmic scale\label{2Antennen9Typ2FS}}
In figure \ref{2Antennen1Typ2} the reconstruction for two antennas with a distance of one is depicted with the corresponding image in the Fourier space given in figure \ref{2Antennen1Typ2FS}. The phantom is acceptably reconstructed and the mirror image is decently suppressed. Unfortunately an error in the low frequencies is annoying that is clearly visible from the large values of $f^{I,F}$ around $\eta=0$ in the cross section of figure \ref{2Antennen1Typ2FS} and the notable slope in the cross section of figure \ref{2Antennen1Typ2}. Figures \ref{2Antennen9Typ2} and \ref{2Antennen9Typ2FS} show the reconstruction and its Fourier transform for two antennas with a distance of nine. In figure \ref{2Antennen9Typ2} the phantom seems to be out of focus with shadows to both sides. Nevertheless the suppression of the mirror image is satisfactory. As can be seen from the steep slope in the cross section in figure \ref{2Antennen1Typ2} and the large values of $f^{I,F}$ around $\eta=0$ in the frequency domain in figure \ref{2Antennen1Typ2FS} in comparison with figure \ref{2Antennen9Typ2FS}, small antenna distances lead to problems with the low frequencies. For large antenna distances there are many frequencies missing on periodical, vertical stripes, as becomes evident from figure \ref{2Antennen9Typ2FS}. This is an effect of the regularization delineated in theorem \ref{antisymformel},1. An advantage is that there are almost no problems with low frequencies.

This points to the possible benefits of a combination of more than two antennas, as it is possible to choose one small and one large distance. The achievable improvements are demonstrated in figures \ref{3Antennen25Typ2} to \ref{4Antennen149Typ2FS} with three and four antennas.
\asbildz{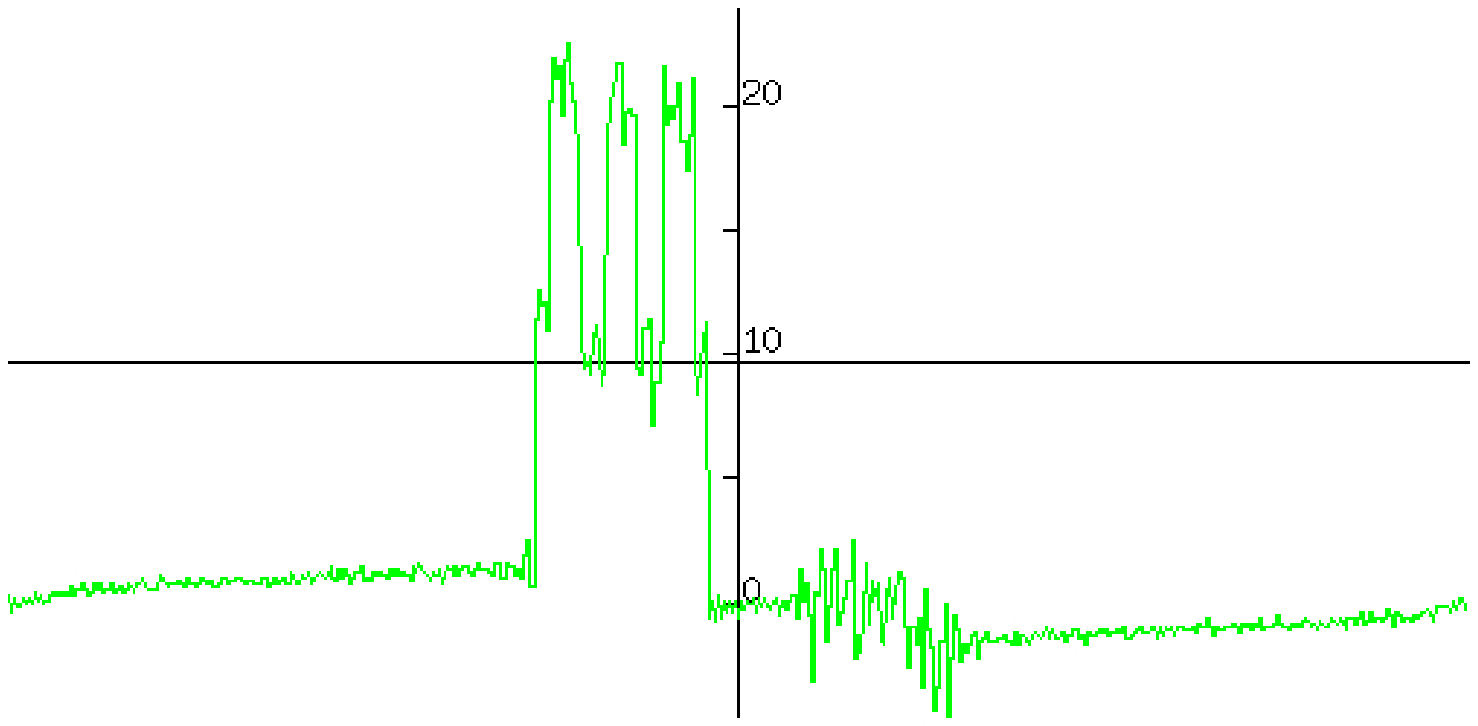}{Three antennas with positions 0, 2, and 5, 10\% noise, spatial domain\label{3Antennen25Typ2}}{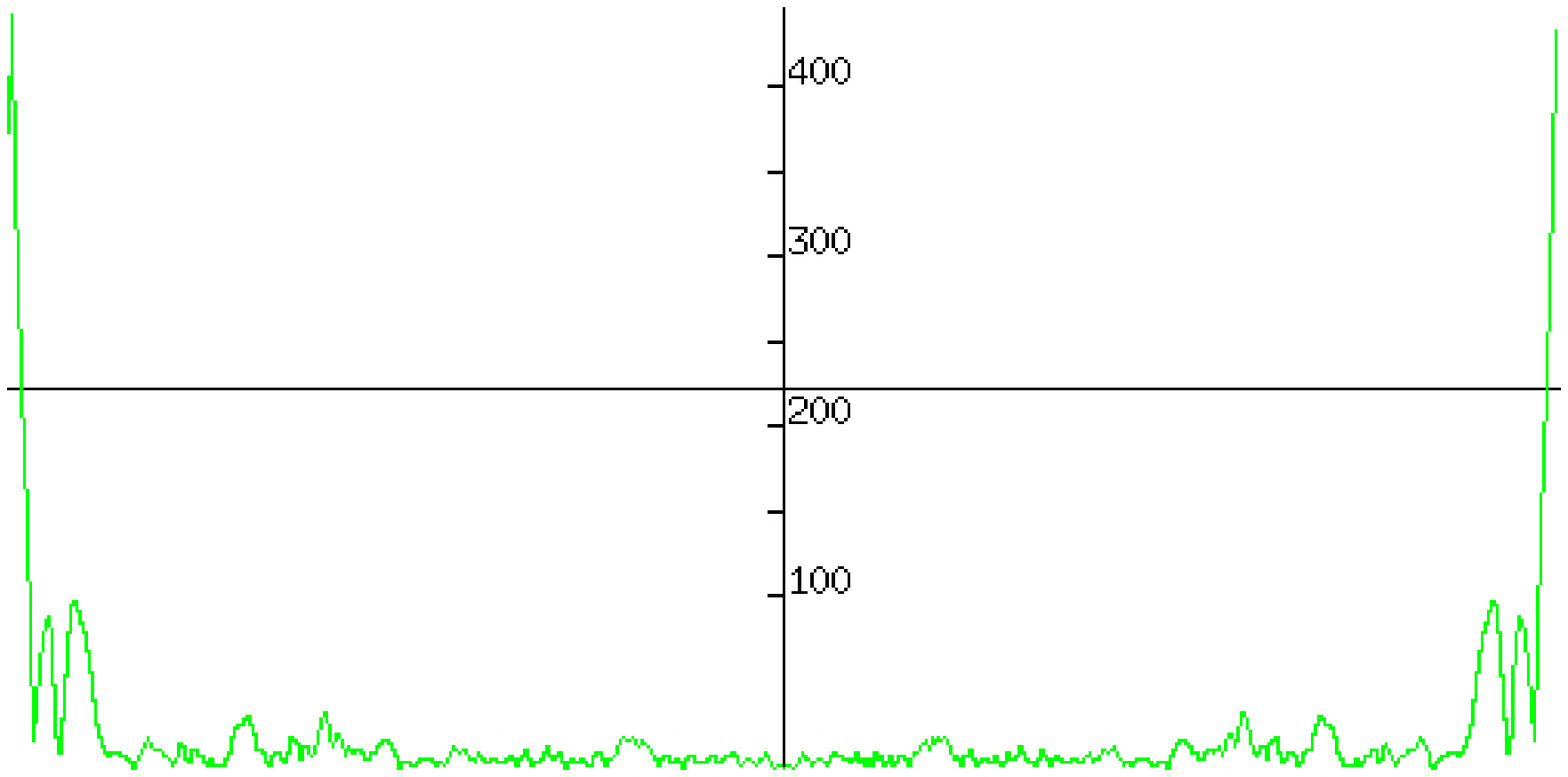}{Three antennas with positions 0, 2, and 5, 10\% noise, frequency domain. Left: Logarithmic scale\label{3Antennen25Typ2FS}}
\asbildz{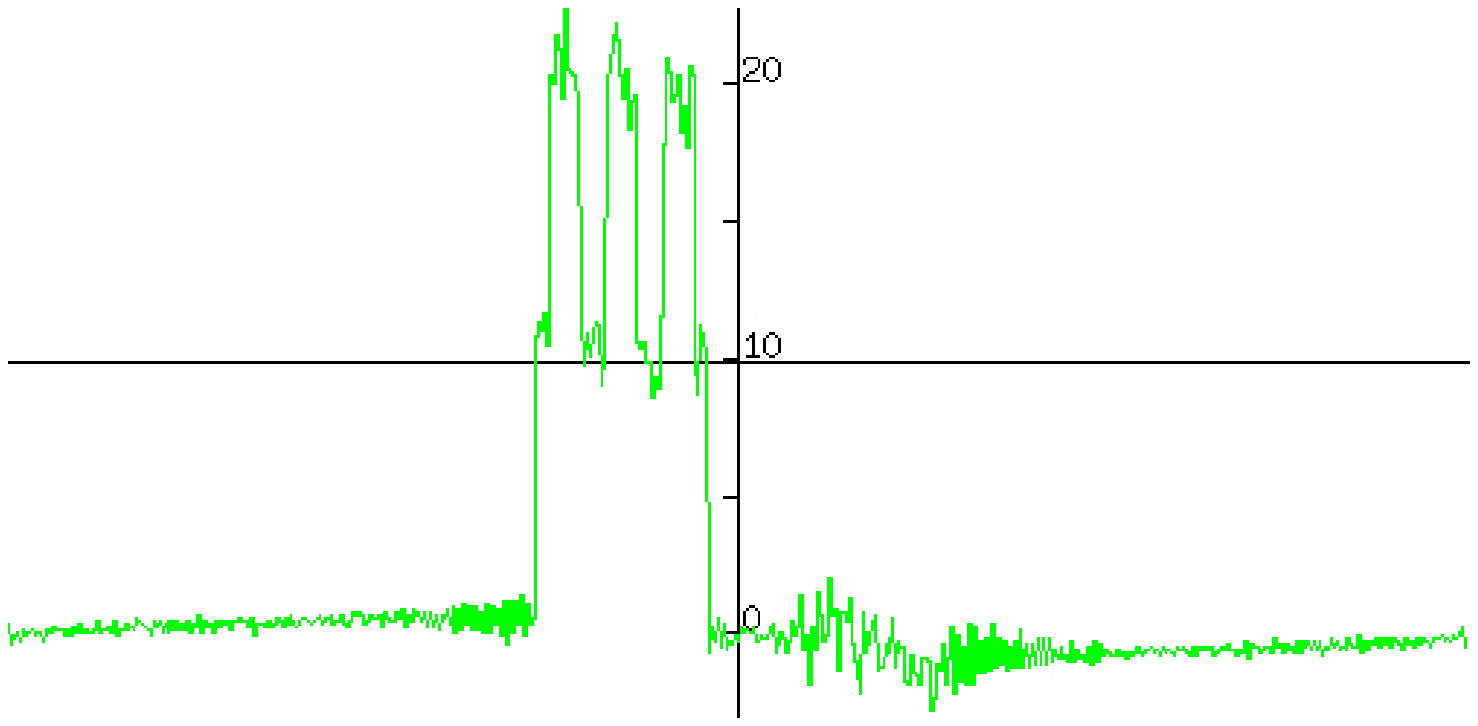}{Four antennas with positions 0, 1, 4, and 9, 10\% noise, spatial domain\label{4Antennen149Typ2}}{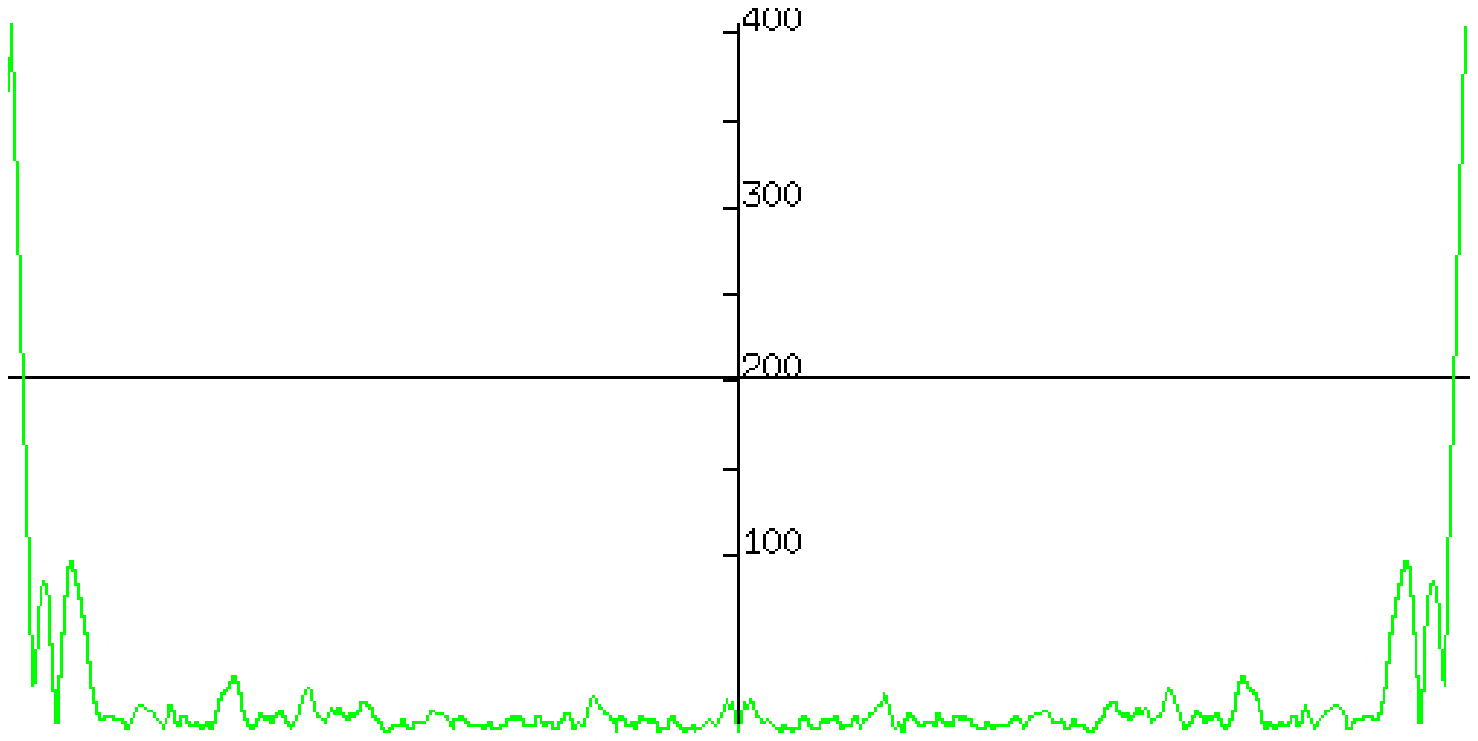}{Four antennas with positions 0, 1, 4, and 9, 10\% noise, frequency domain. Left: Logarithmic scale\label{4Antennen149Typ2FS}}
In figures \ref{3Antennen25Typ2} and \ref{3Antennen25Typ2FS} the reconstruction and its Fourier transform are depicted for three antennas with distances two and five. The reconstruction quality improved in comparison with the preceding figures, but there is still a slight slope visible in the cross section of figure \ref{3Antennen25Typ2}. However,  there are no frequencies missing since a regularization is no longer necessary due to data from more than two antennas (see theorem \ref{antisymformel}). In figures \ref{4Antennen149Typ2} and \ref{4Antennen149Typ2FS} the result from four antennas with distances one, four, and nine can be seen. The quality of these figures is very good, since the phantom is reconstructed very accurately and the mirror image is only slightly visible. As can be seen from the preceding images, the reconstruction quality improves with more antennas. Since there are problems, if only two antennas are used, regardless of the distance, it is advisable to employ at least three antennas to avoid these problems. The effects of more antennas with regards to noise will be further discussed in subsection \ref{noise}.

\subsection{Solving the left-right ambiguity with reconstructed symmetric images}

\label{spherical}

In the following examples measurements for several parallel flight tracks are simulated using the spherical Radon transform. Then the corresponding even images are formed with the formula in theorem \ref{modinvform}. Finally, the images are combined to obtain the ground reflectivity function $f$ using the formulas in theorem \ref{antisymformel}.\\
As an example, figure \ref{measuredtype1} shows the data from the simulated measurements where the number of boomerang-shapes depends on the number of antennas employed and the positions of the shapes are related to the antenna distances. The noise level will be at a level of $10 \%$ throughout this subsection.

\asbildg{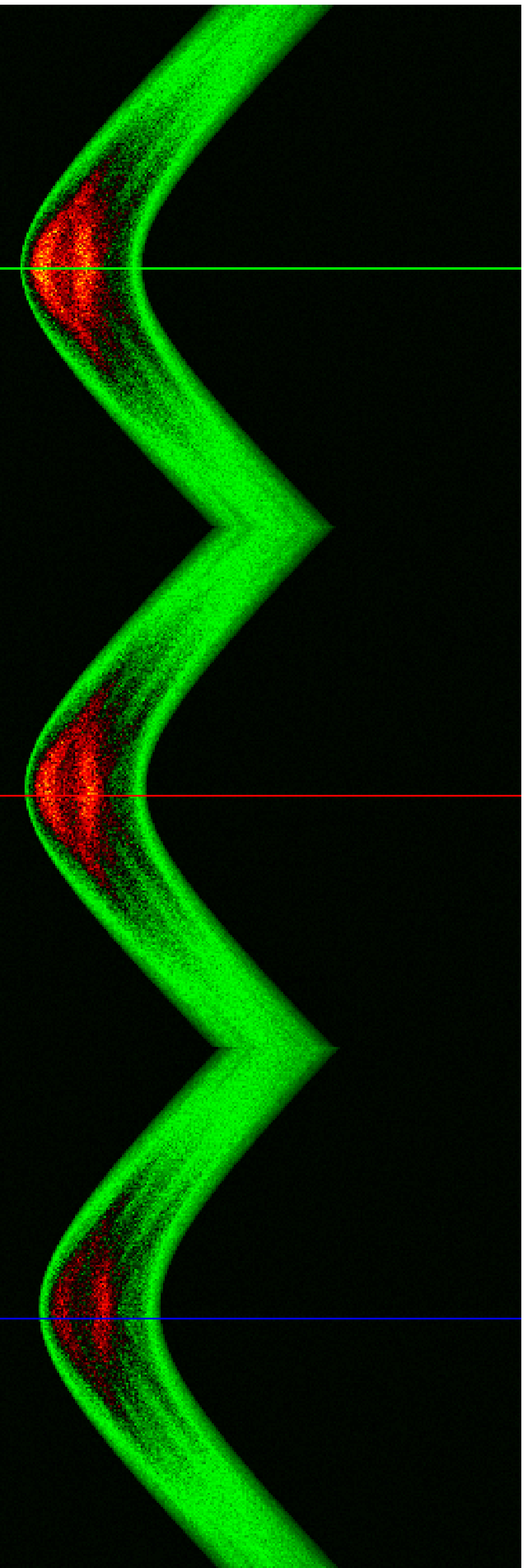}{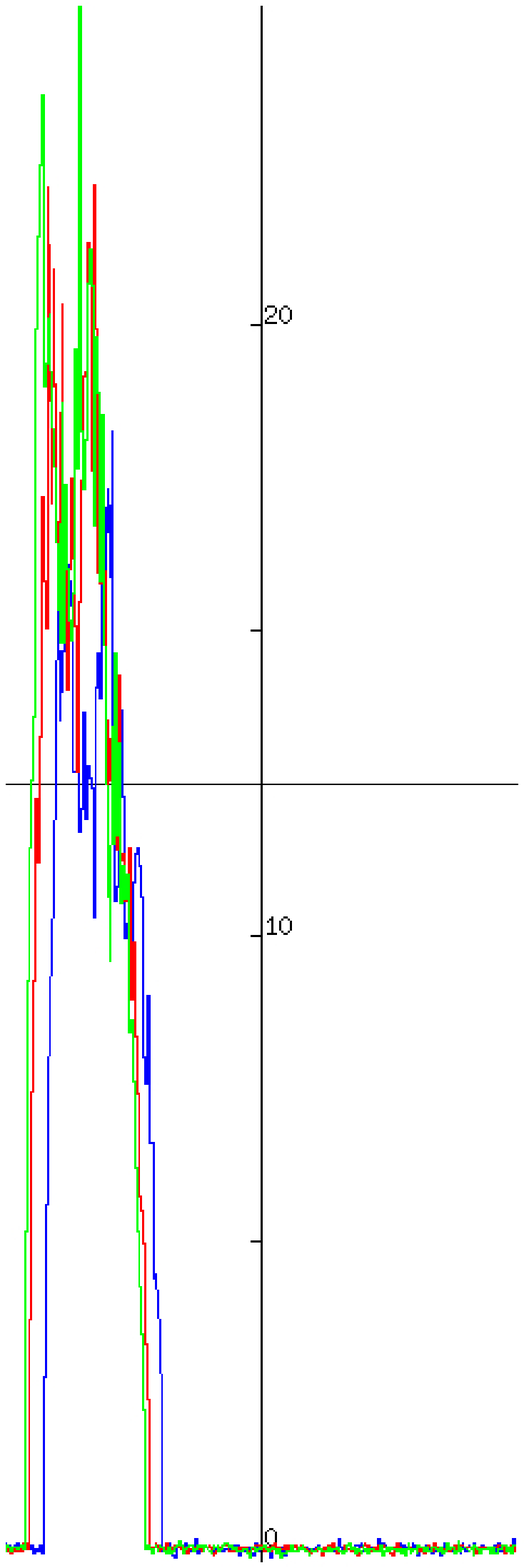}{Measured data from three antennas with 10\% noise\label{measuredtype1}}


In the following, the effect of the antenna distance on the reconstruction quality is examined. Figures \ref{2Antennen1Typ1} to \ref{2Antennen9Typ1FS} contrast the reconstruction with a small antenna distance to a reconstruction with a large antenna distance.

\asbildz{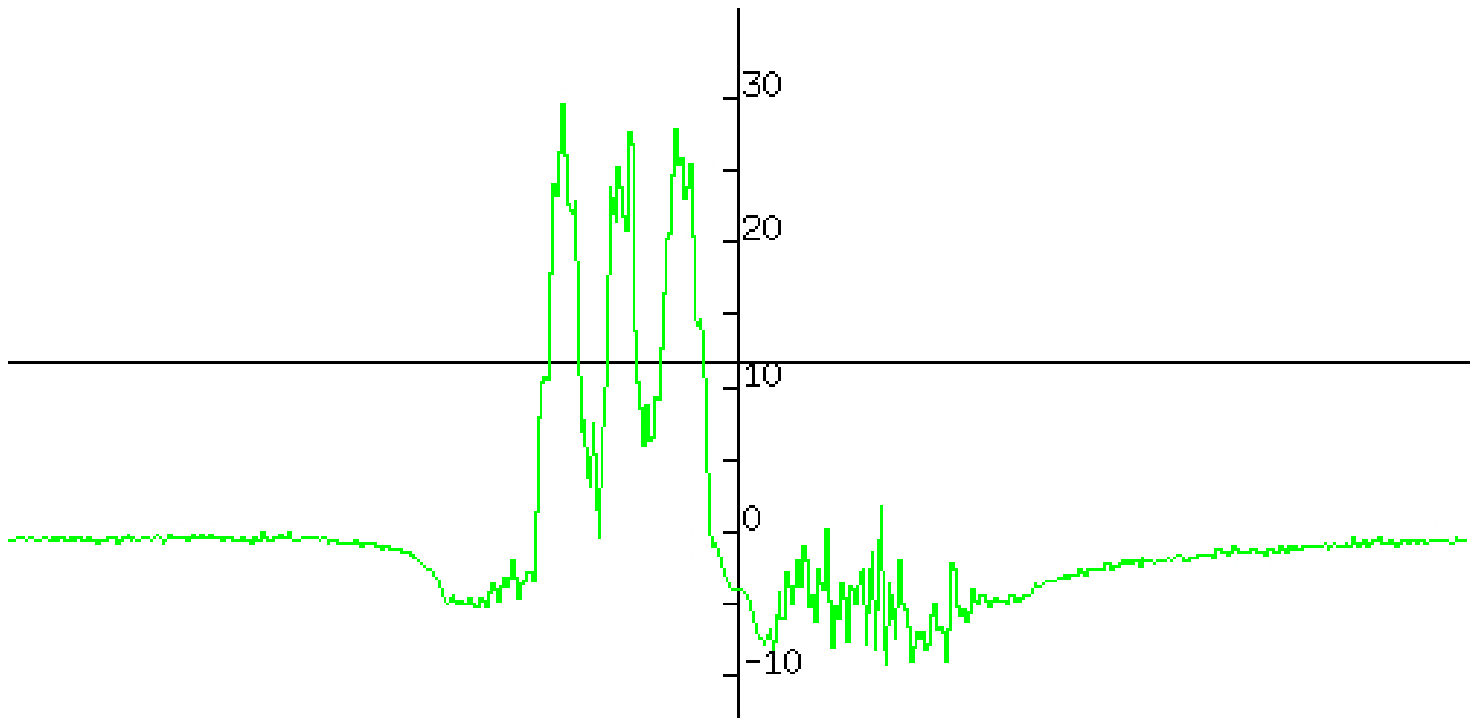}{Two antennas with distance 1, 10\% noise, spatial domain\label{2Antennen1Typ1}}{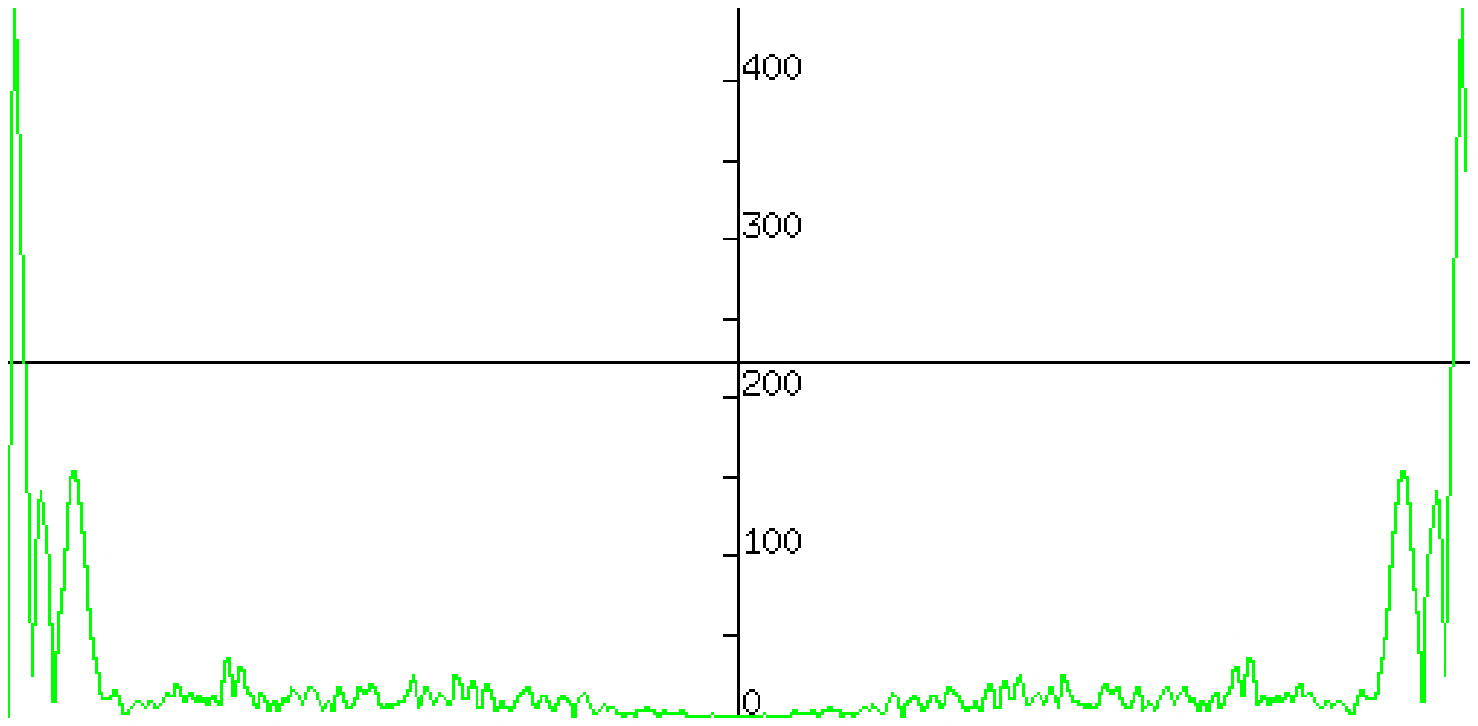}{Two antennas with distance 1, 10\% noise, frequency domain. Left: Logarithmic scale\label{2Antennen1Typ1FS}}
\asbildz{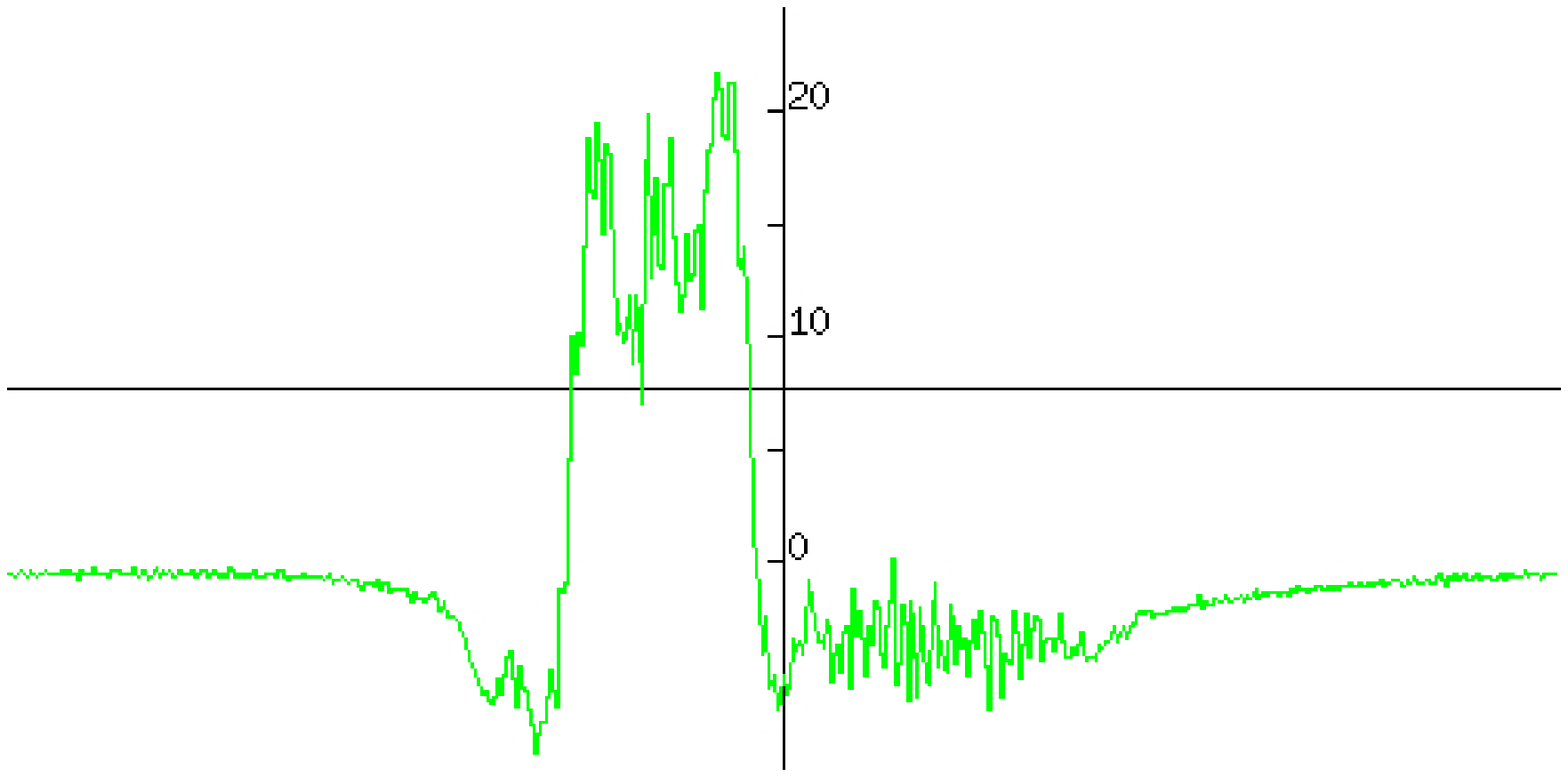}{Two antennas with distance 9, 10\% noise, spatial domain\label{2Antennen9Typ1}}{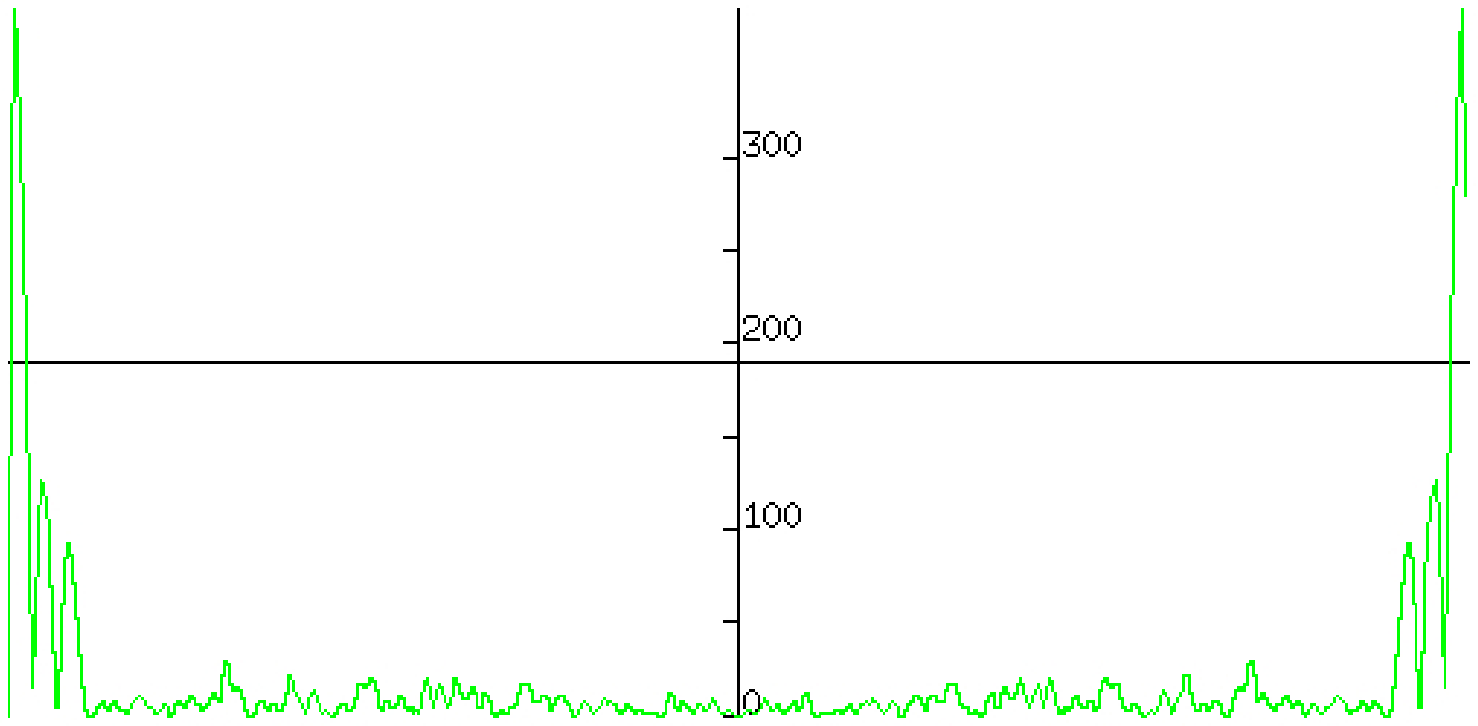}{Two antennas with distance 9, 10\% noise, frequency domain. Left: Logarithmic scale\label{2Antennen9Typ1FS}}

Figures \ref{2Antennen1Typ1} and \ref{2Antennen1Typ1FS} display the reconstruction and its Fourier transform for two antennas with a distance of one. The reconstruction of the phantom is satisfactory as is the suppression of the mirror image. In comparison with figures \ref{2Antennen1Typ2} and \ref{2Antennen1Typ2FS}, which show the corresponding images with the same number and placement of the antennas that are reconstructed from computed symmetric images, the error in the low frequencies for small antenna distances is negligible, however figure \ref{2Antennen1Typ1} is affected by circular artifacts that arise from the noise in the inversion of the spherical Radon transform. The small error in the low frequencies in figure \ref{2Antennen1Typ1} comes from a suppression of the low frequencies due to limited data in the inversion of the spherical Radon transform. Figures \ref{2Antennen9Typ1} and \ref{2Antennen9Typ1FS} show the reconstruction for two antennas with a distance of nine and the corresponding image in the frequency domain. The phantom in figure \ref{2Antennen9Typ1} is very blurry - an effect of the many missing frequencies that can be observed in figure \ref{2Antennen9Typ1FS}. Surprisingly however, the amplitude of the phantom is more exactly reconstructed in comparison to figure \ref{2Antennen1Typ1}. This points again to a more exact reconstruction of the low frequencies for large antenna distances, although the effect is not as striking as in subsection \ref{Typ2}. This is also reflected in the slightly smaller maxima in figure \ref{2Antennen9Typ1FS}.


\asbildz{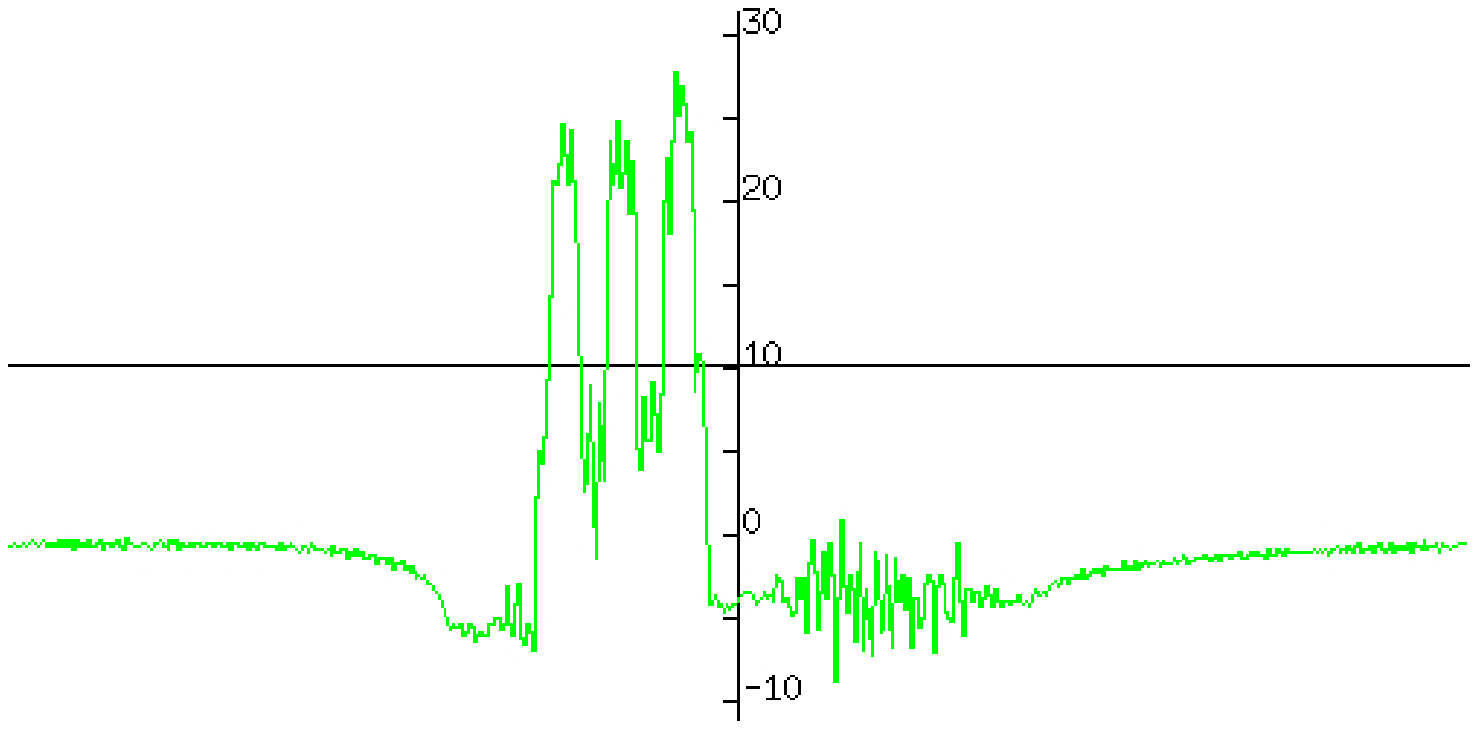}{Three antennas with positions 0, 2, and 5, 10\% noise, spatial domain\label{3Antennen25Typ1}}{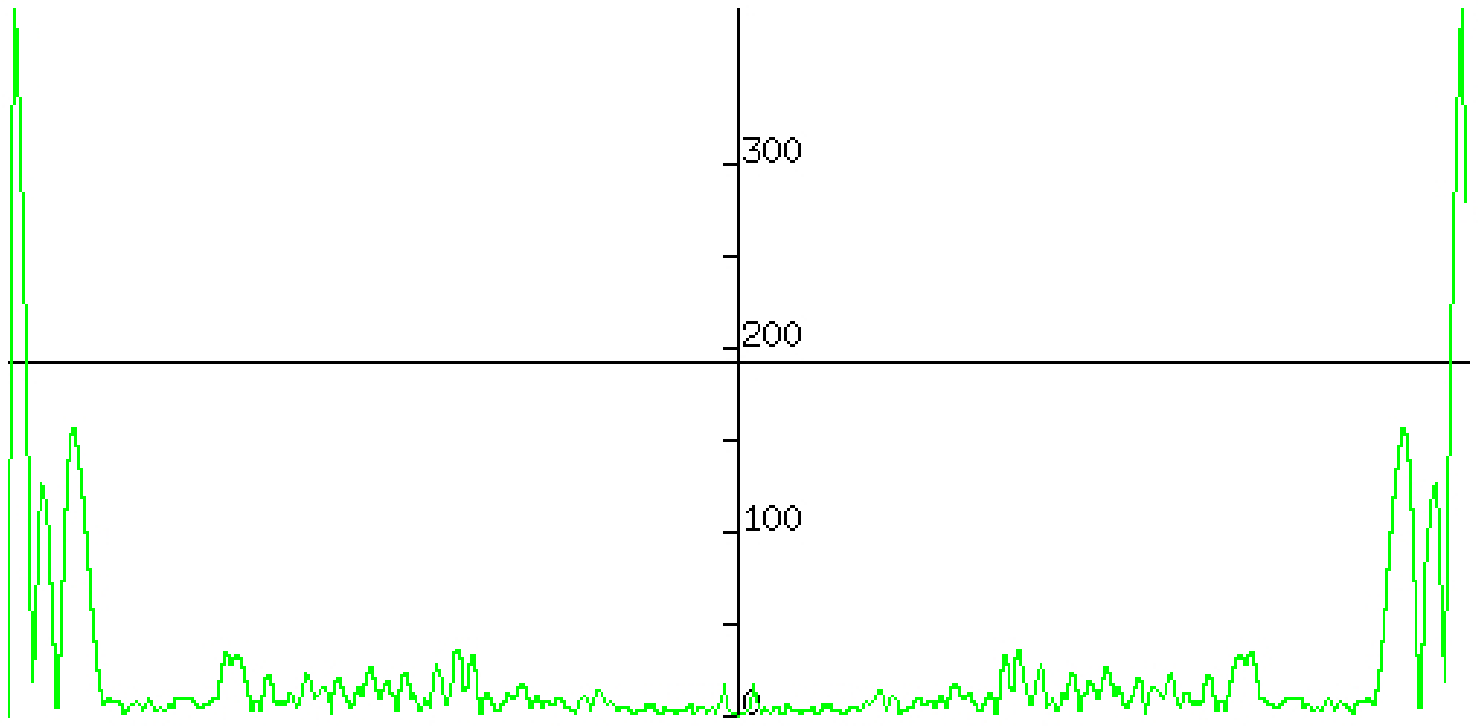}{Three antennas with positions 0, 2, and 5, 10\% noise, frequency domain. Left: Logarithmic scale\label{3Antennen25Typ1FS}}

Again, this emphasizes that it could be an interesting idea to combine more than two antennas as it could allow to gain the advantages of a small and a large antenna distance. The successes are demonstrated in figures \ref{3Antennen25Typ1} to \ref{4Antennen149Typ1FS}.

\asbildz{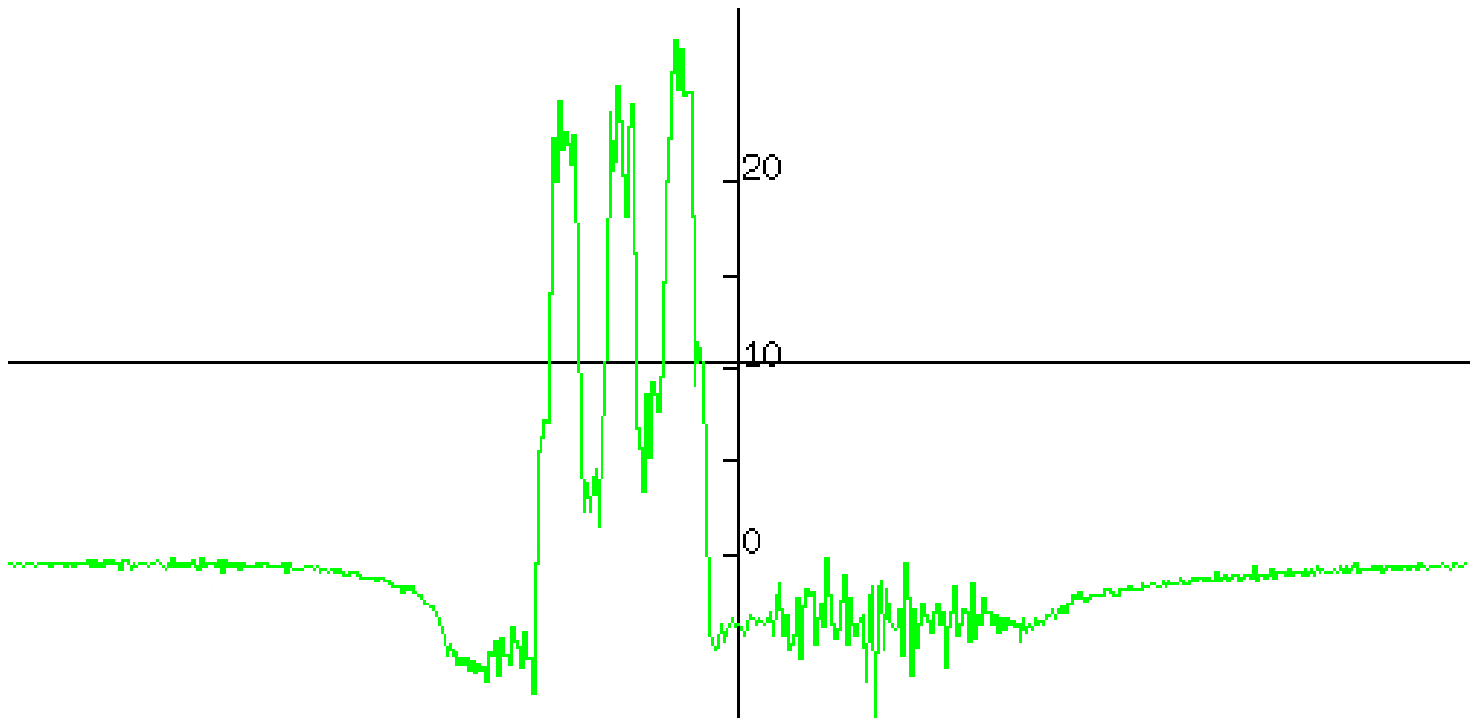}{Four antennas with positions 0, 1, 4, and 9, 10\% noise, spatial domain\label{4Antennen149Typ1}}{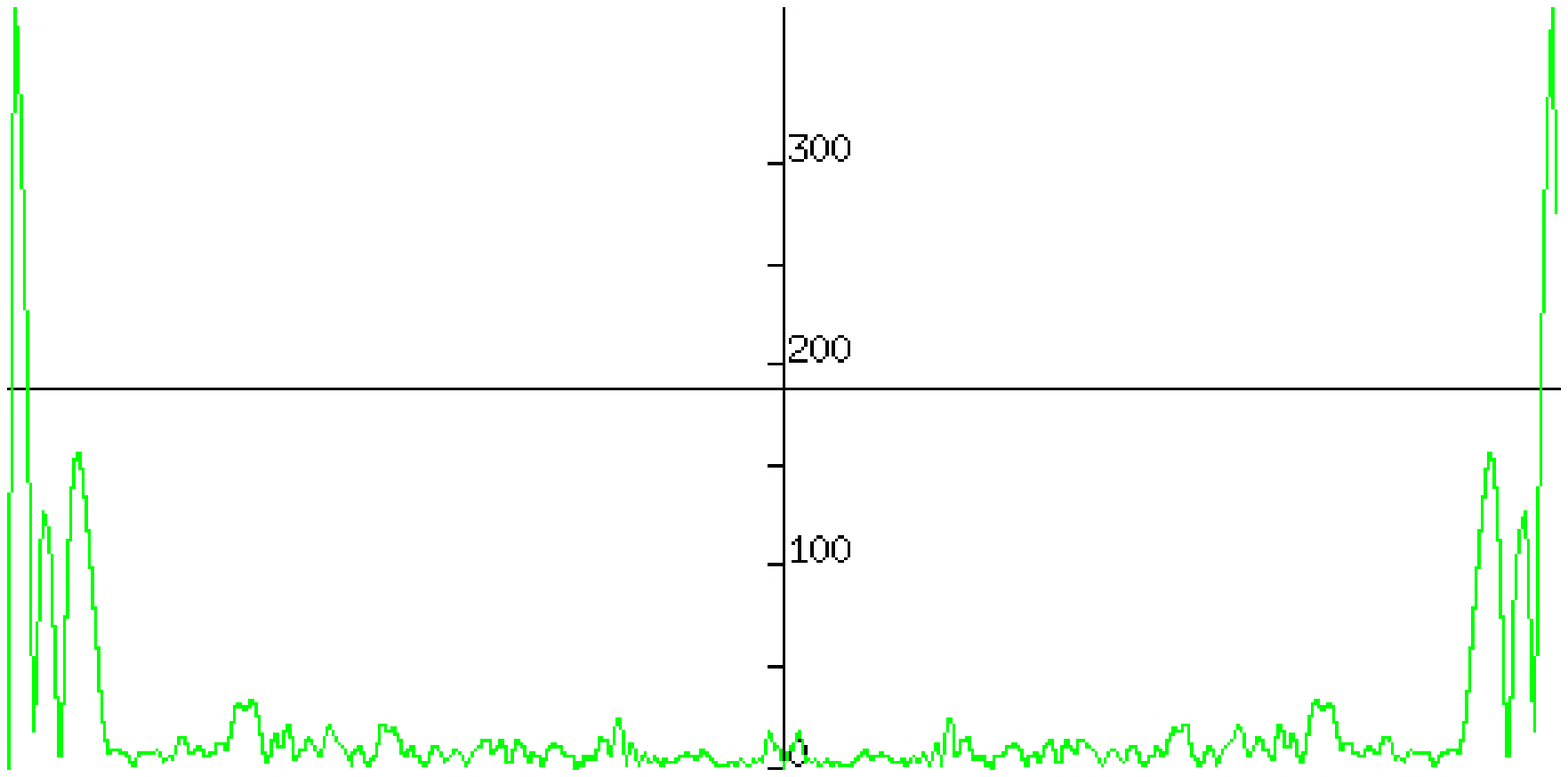}{Four antennas with positions 0, 1, 4, and 9, 10\% noise, frequency domain. Left: Logarithmic scale\label{4Antennen149Typ1FS}}

In figures \ref{3Antennen25Typ1} and \ref{3Antennen25Typ1FS} the results from three antennas with distances two and five can be seen. The reconstruction is superior to the preceding two, but the amplitude of the phantom is still not completely correct, as can be seen in the cross section of figure \ref{3Antennen25Typ1}. Again, with three antennas there are no frequencies missing, as seen in figure \ref{3Antennen25Typ1FS}, since no regularization is necessary. Figures \ref{4Antennen149Typ1} and \ref{4Antennen149Typ1FS} show the reconstruction and its Fourier transform for four antennas with distances one, four, and nine. The reconstruction improved over the preceding one, but still the amplitude of the phantom is not correct. This is an effect of the inversion of the spherical Radon transform with limited data. Nevertheless the suppression of the mirror image is good. As could be presumed, the low frequencies are reconstructed quite nicely and there are no frequencies missing, as can be seen in all four figures.
The reconstruction quality improves further with more antennas, as can be seen in comparison of figures \ref{2Antennen1Typ1}, \ref{2Antennen9Typ1}, \ref{3Antennen25Typ1}, and \ref{4Antennen149Typ1}. This is not only an effect of averaging out the noise, but is amplified by a weighted summation of the data, where the data is more heavily weighted the more noise resistant it is.


\subsection{Analysis of stability and the effects of noise\label{noise}}

Now an analysis of noise effects and stability is performed.\\
Figures \ref{Vergleich1Typ2} to \ref{30Vergleich1Typ1} show reconstructions from two antennas with varying amounts of noise. As in subsection \ref{symmetric}, figures \ref{Vergleich1Typ2} to \ref{30Vergleich1Typ2} demonstrate only the characteristics of the solution of the left-right ambiguity, whereas figures \ref{Vergleich1Typ1} to \ref{30Vergleich1Typ1} contain also the artifacts arising due to the inversion of the spherical Radon transform with limited data as in subsection \ref{spherical}.

\asbildd{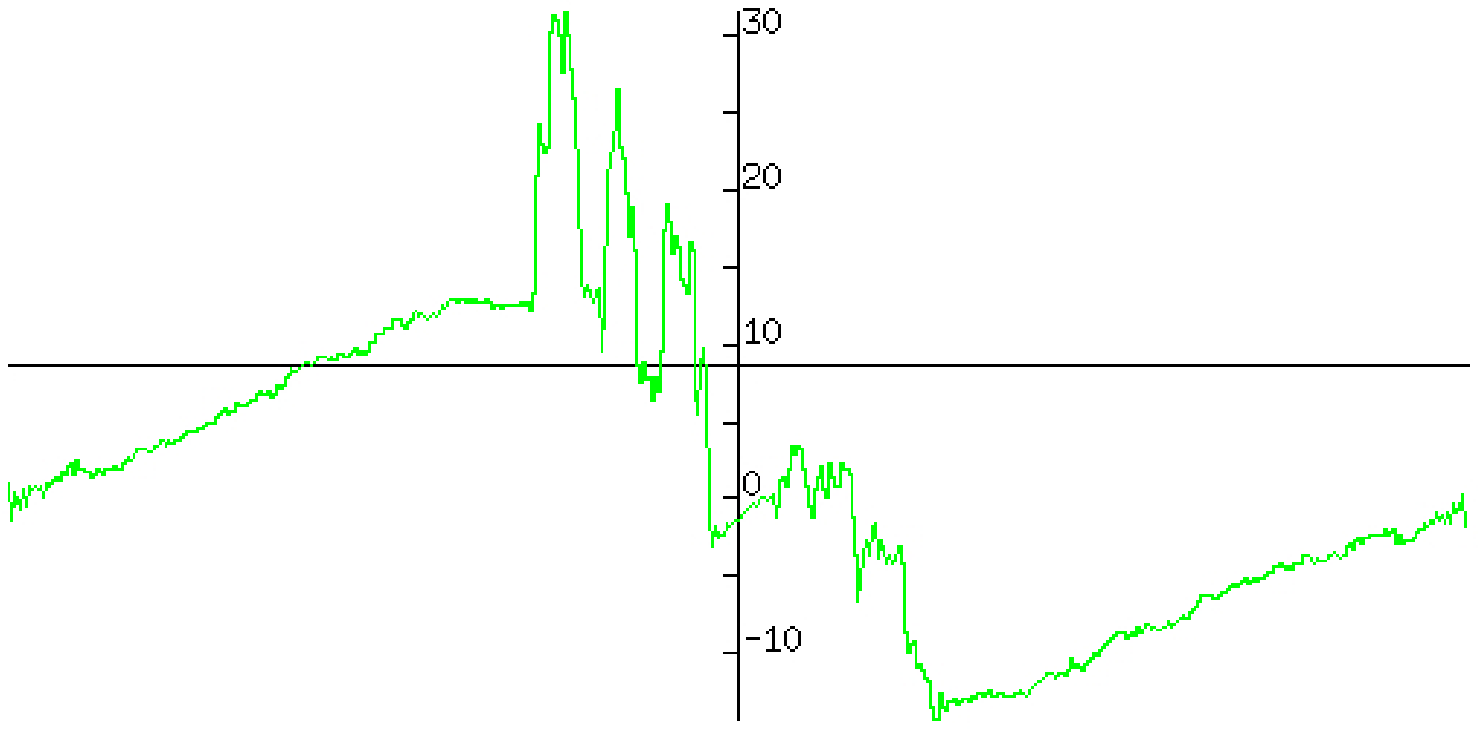}{Two antennas with positions 0 and 1, no spherical Radon transform, 10\% noise\label{Vergleich1Typ2}}{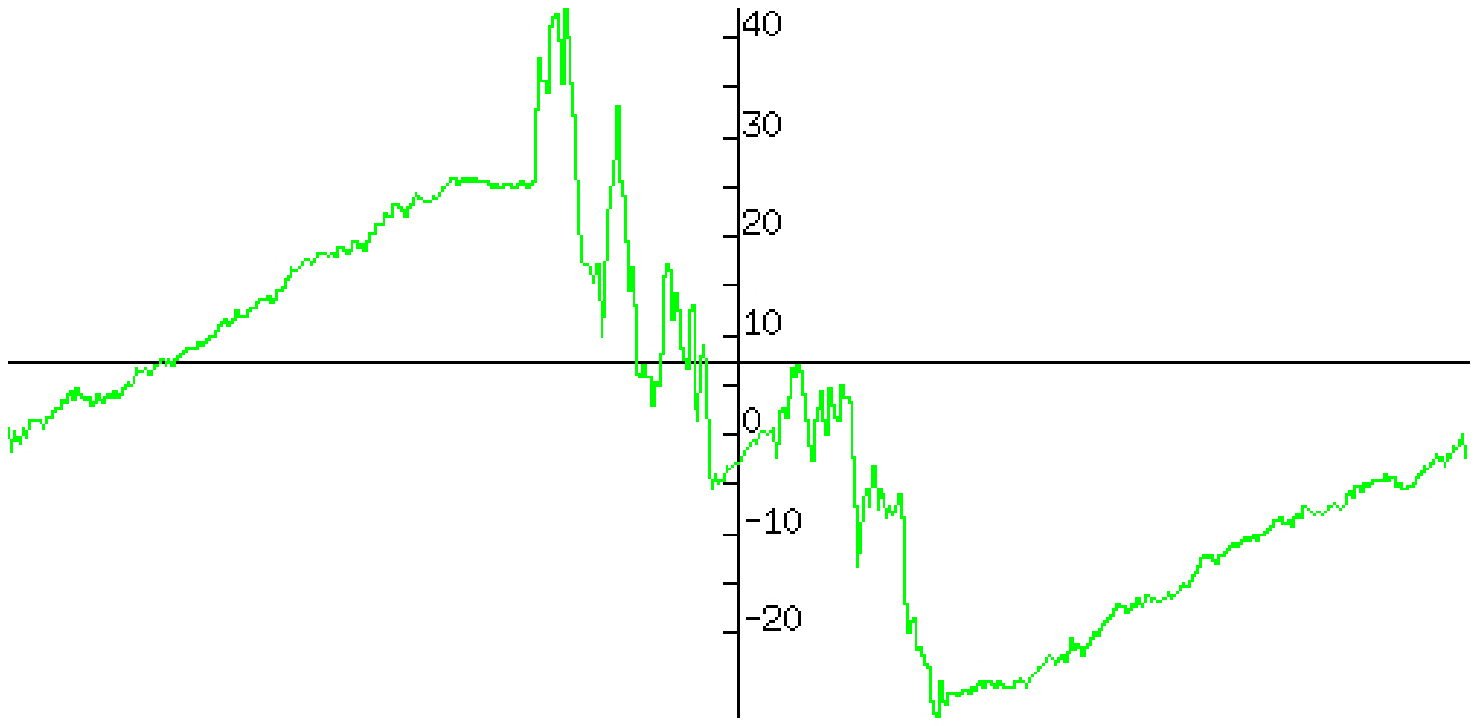}{Two antennas with positions 0 and 1, no spherical Radon transform, 20\% noise\label{20Vergleich1Typ2}}{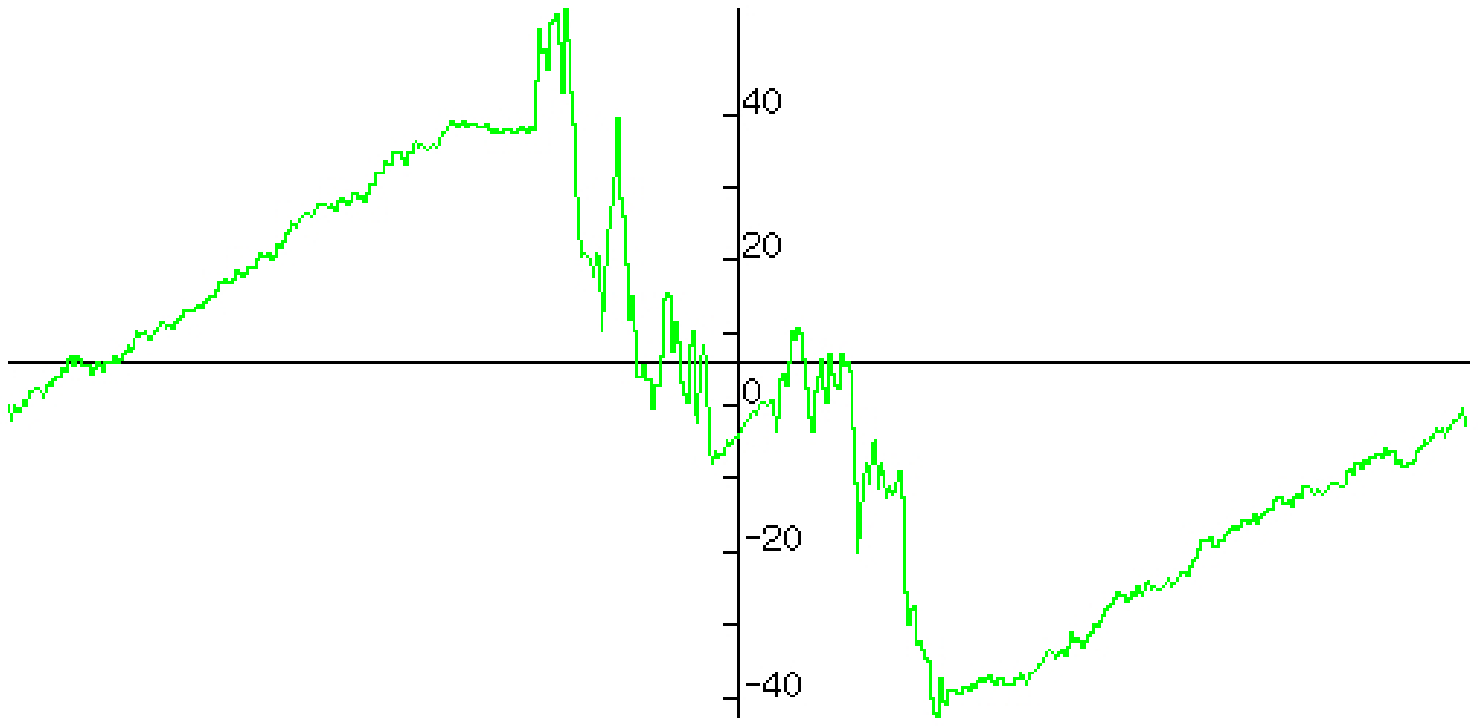}{Two antennas with positions 0 and 1, no spherical Radon transform, 30\% noise\label{30Vergleich1Typ2}}

Figures \ref{Vergleich1Typ2} to \ref{30Vergleich1Typ2} show the reconstructions from computed even images with respect to two antennas with a distance of one. In figure \ref{Vergleich1Typ2} the noise level is at 10\%, in figure \ref{20Vergleich1Typ2} at 20\%, and in figure \ref{30Vergleich1Typ2} at 30\%. In figure \ref{Vergleich1Typ2} the phantom is clearly visible despite the known error in the low frequencies. The phantom in figure \ref{20Vergleich1Typ2} is obscured more strongly due to this problem, and in figure \ref{30Vergleich1Typ2} it is very difficult to discern the phantom from the underlying error. Nevertheless it can be seen in the cross sections that the algorithm is stable since the cross sections look quite similar disregarding the scaling, and there are no signs of an overamplification of noise.

\asbildd{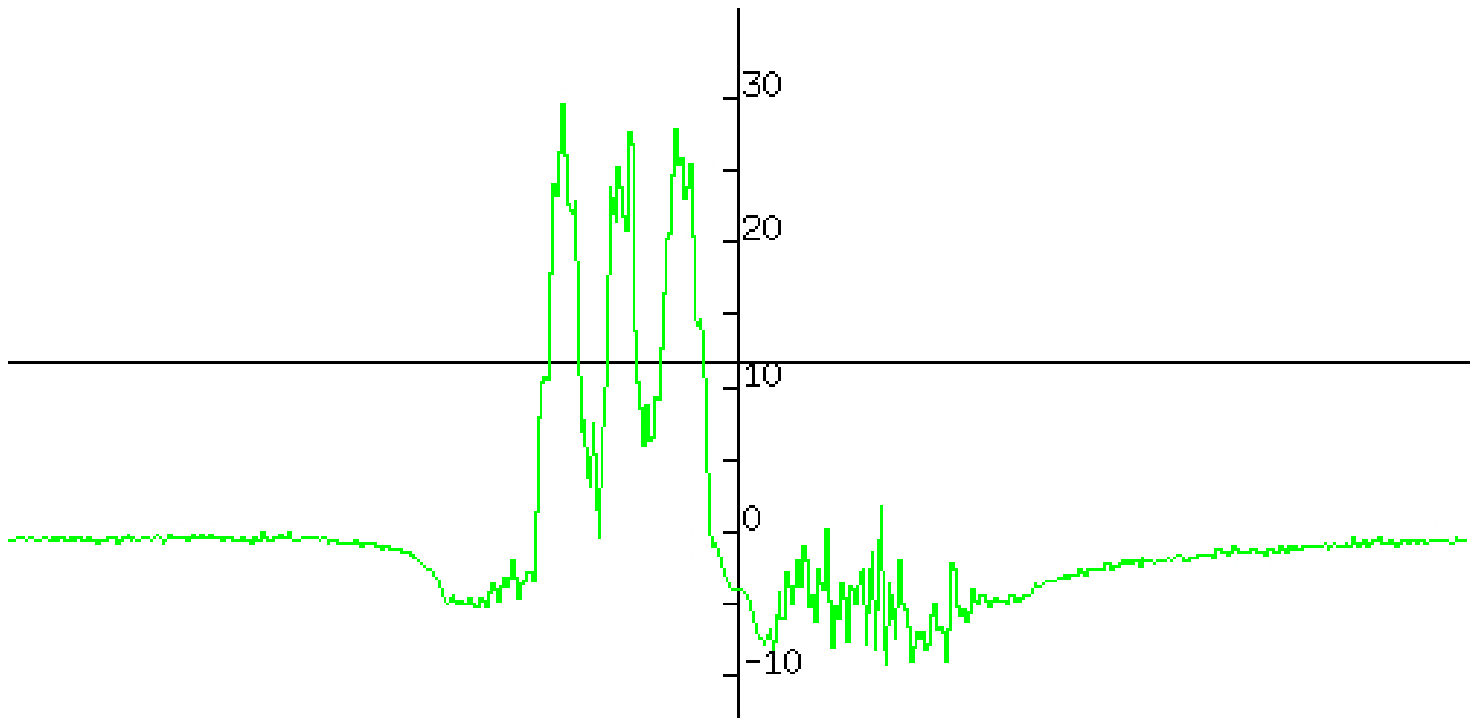}{Two antennas with positions 0 and 1, with spherical Radon transform, 10\% noise\label{Vergleich1Typ1}}{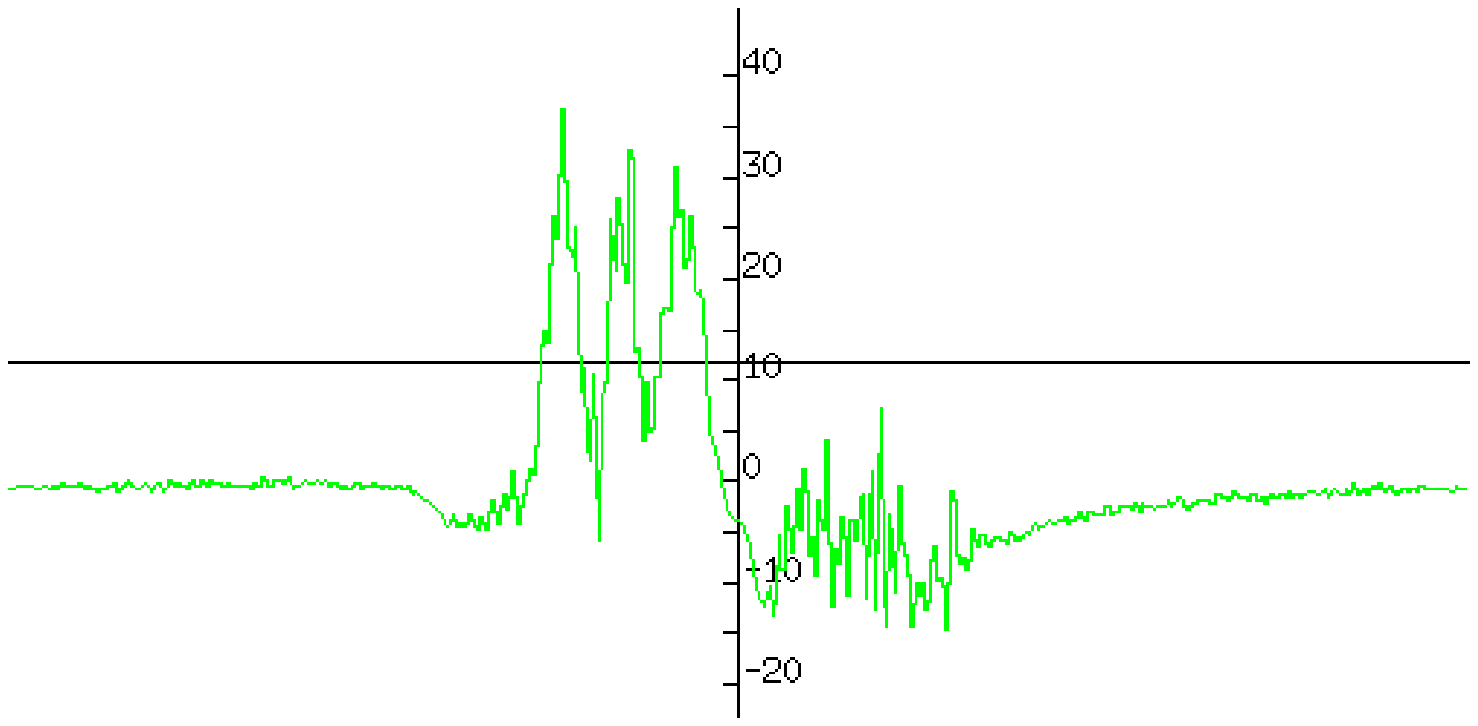}{Two antennas with positions 0 and 1, with spherical Radon transform, 20\% noise\label{20Vergleich1Typ1}}{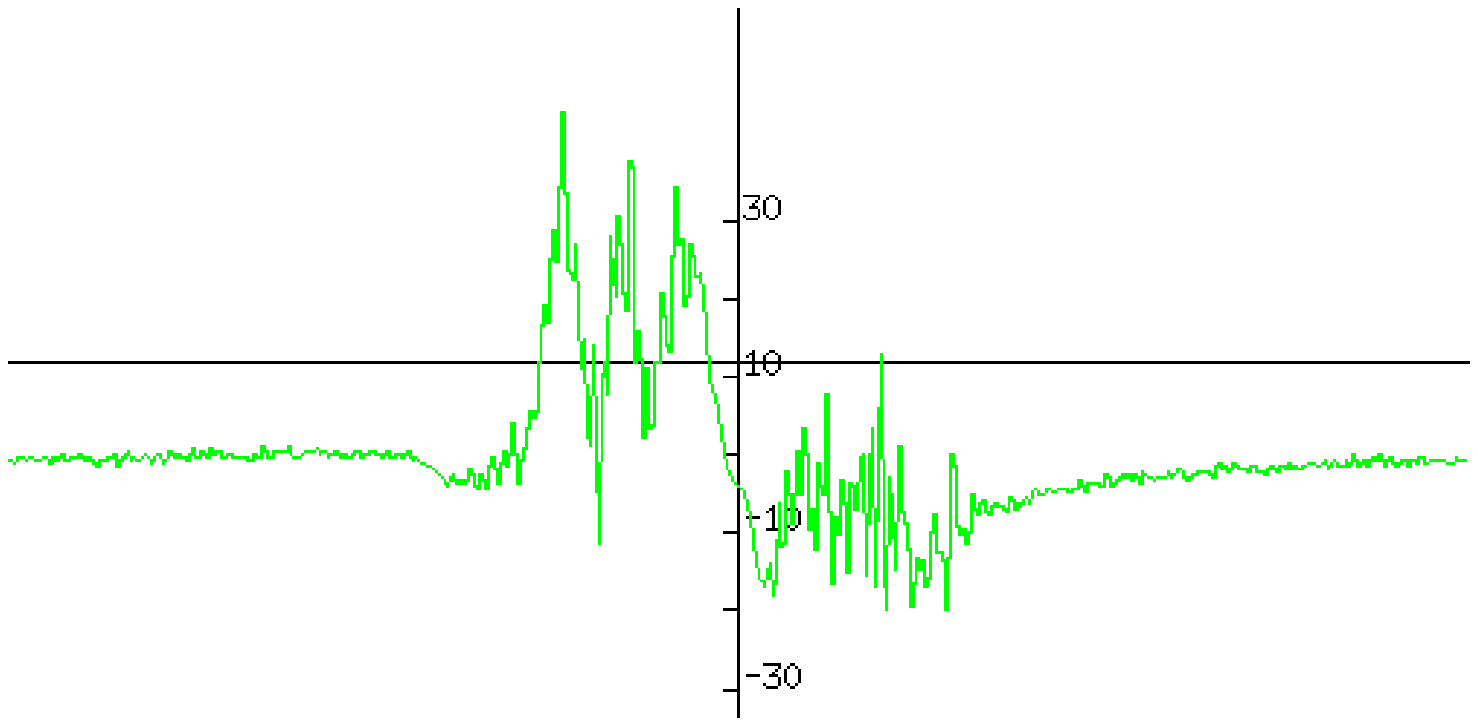}{Two antennas with positions 0 and 1, with spherical Radon transform, 30\% noise\label{30Vergleich1Typ1}}

In figures \ref{Vergleich1Typ1} to \ref{30Vergleich1Typ1} the reconstructions using the inversion of the spherical Radon transform for two antennas with a distance of one are depicted. In figure \ref{Vergleich1Typ1} the noise level is at 10\%, in figure \ref{20Vergleich1Typ1} at 20\%, and in figure \ref{30Vergleich1Typ1} at 30\%. The reconstruction quality in figure \ref{Vergleich1Typ1} is again acceptable, but declines for higher noise levels as in figure \ref{20Vergleich1Typ1}. Finally, for a noise level of 30\% the reconstruction quality is not very good, but the phantom is still discernible in contrast to figure \ref{30Vergleich1Typ2}. All three images do not display the right amplitude of the phantom, but as previously analyzed, this is an effect of the inversion with limited data. With increasing noise levels the reconstruction degrades, but no uncontrolled noise amplification can be observed.

Figures \ref{Vergleich1Typ2} to \ref{30Vergleich1Typ1} show that the reconstruction quality gets worse with a higher amount of noise, since the amplitude of the error gets larger. Nevertheless there are no signs of an overamplification of noise. The loss of reconstruction quality is within the expected order of magnitude, so this points to the stability of the algorithm.\\

Figures \ref{Vergleich13Typ2} to \ref{30Vergleich13Typ1} compare the results for calculations with three antennas with different noise severity. Figures \ref{Vergleich13Typ2} to \ref{30Vergleich13Typ2} are directly computed from symmetric images, whereas figures \ref{Vergleich13Typ1} to \ref{30Vergleich13Typ1} are gained from inversions of the spherical Radon transform.

\asbildd{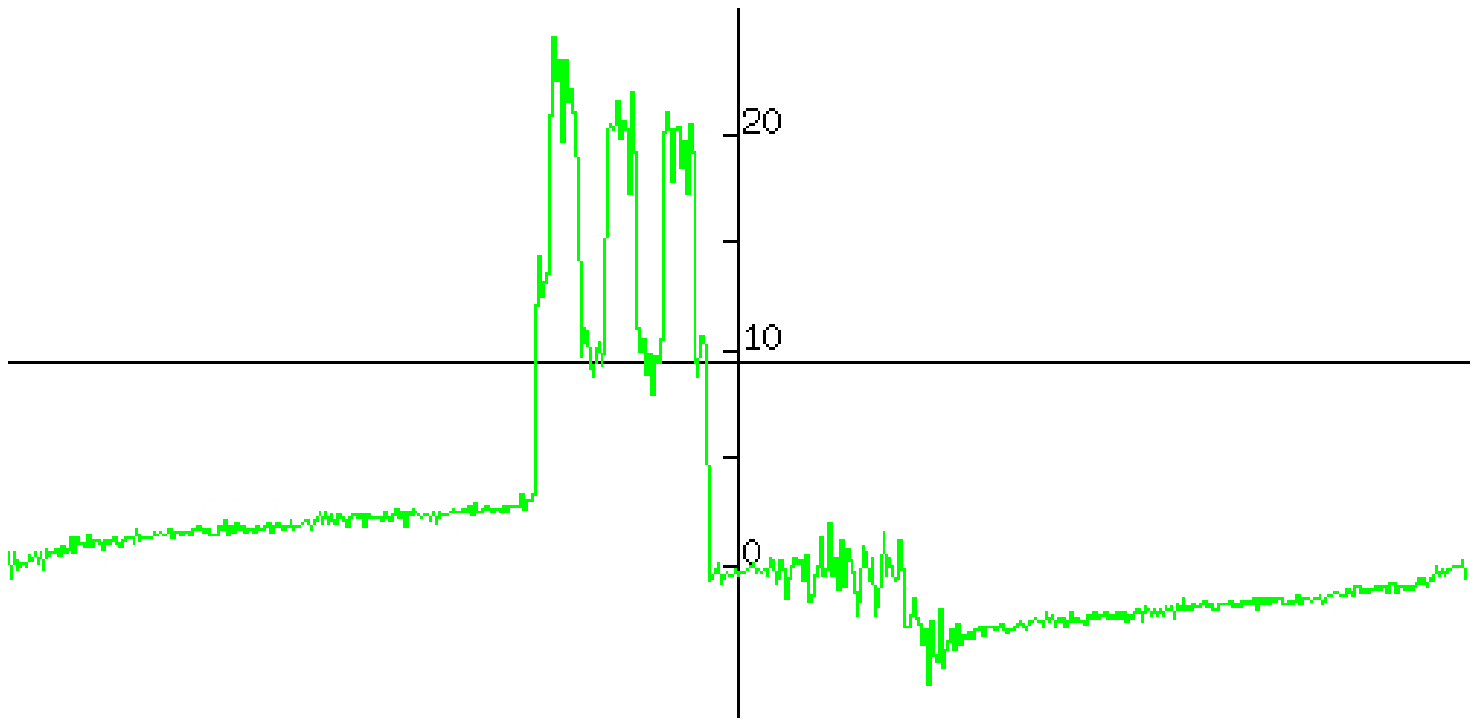}{Three antennas with positions 0, 1, and 3, no spherical Radon transform, 10\% noise\label{Vergleich13Typ2}}{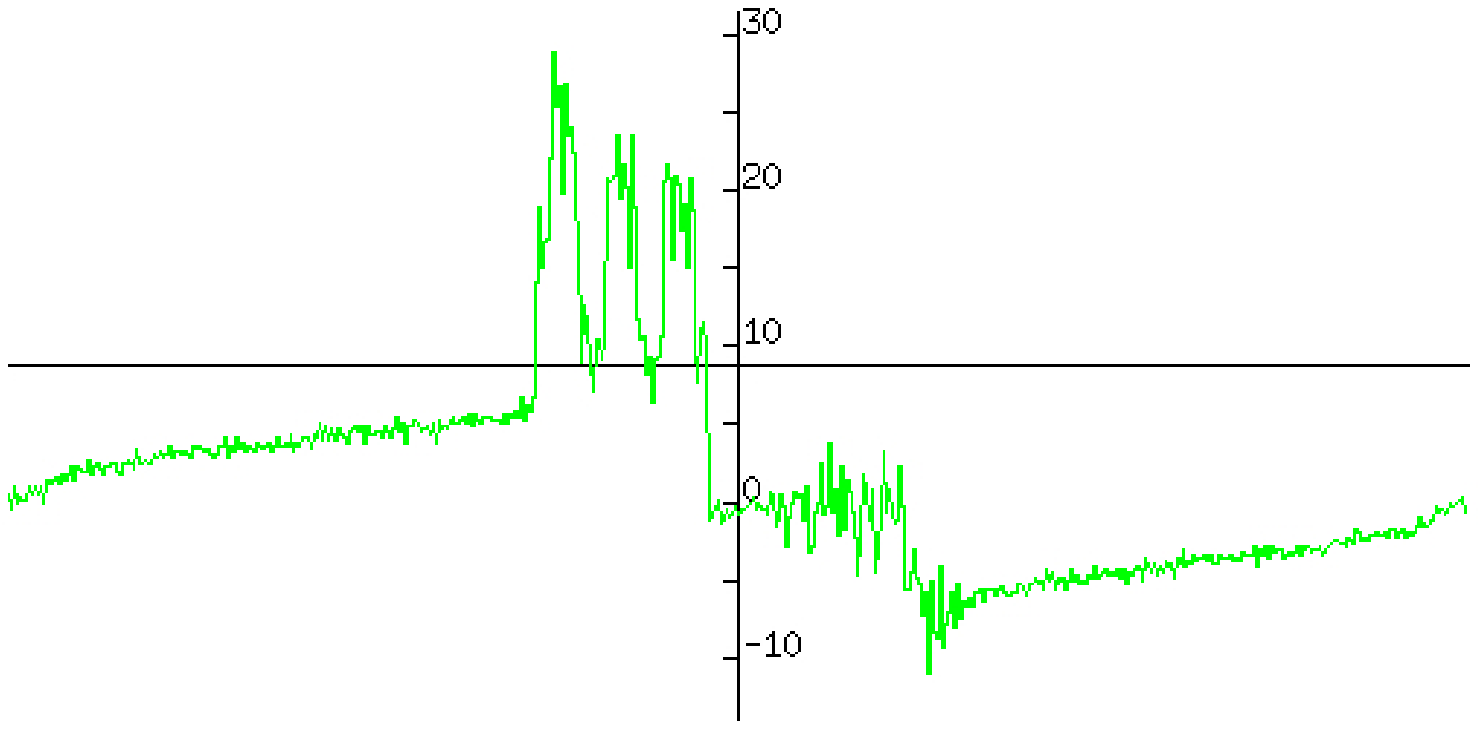}{Three antennas with positions 0, 1, and 3, no spherical Radon transform, 20\% noise\label{20Vergleich13Typ2}}{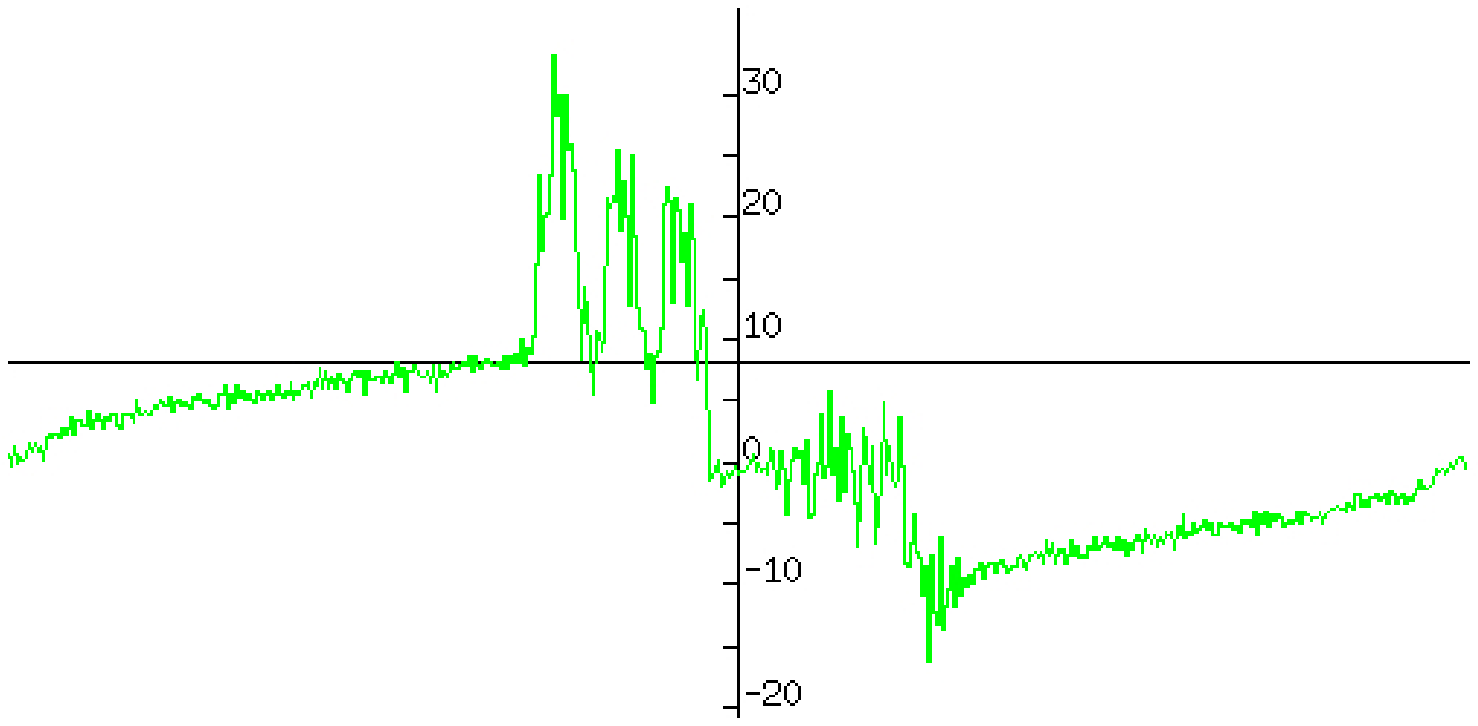}{Three antennas with positions 0, 1, and 3, no spherical Radon transform, 30\% noise\label{30Vergleich13Typ2}}

Figures \ref{Vergleich13Typ2} to \ref{30Vergleich13Typ2} show the reconstructions using computed symmetric images with regard to three antennas with distances one and three. The noise level in figure \ref{Vergleich13Typ2} is at 10\%, in figure \ref{20Vergleich13Typ2} at 20\%, and in figure \ref{30Vergleich13Typ2} at 30\%. The reconstruction of the phantom in figure \ref{Vergleich13Typ2} is satisfactory since the amplitude is close to the correct values. The suppression of the mirror image is also satisfying. However there is still a problem with low frequencies as can be seen in the slight slope in the cross section in figure \ref{Vergleich13Typ2}. The reconstruction quality in figure \ref{20Vergleich13Typ2} is worse. The phantom has a higher peak and the mirror image is more strongly visible. Additionally, the slope in the cross section of figure \ref{20Vergleich13Typ2} is steeper which points to a more severe problem with low frequencies. This tendency continues in figure \ref{30Vergleich13Typ2}. The amplitude of the phantom is again farther off, albeit not by much, and the mirror image is more pronounced. Also the slope in the cross section of figure \ref{30Vergleich13Typ2} steepened a bit.

\asbildd{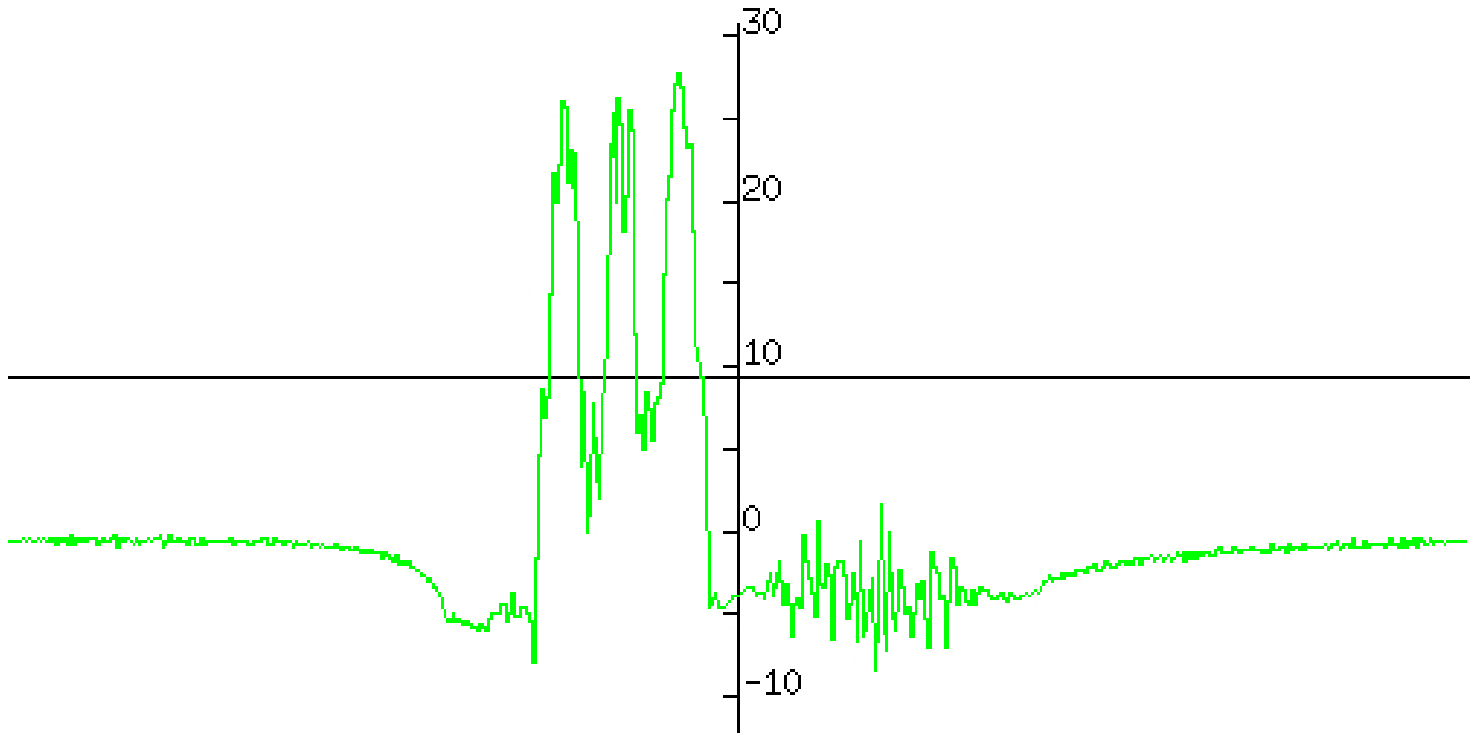}{Three antennas with positions 0, 1, and 3, with spherical Radon transform, 10\% noise\label{Vergleich13Typ1}}{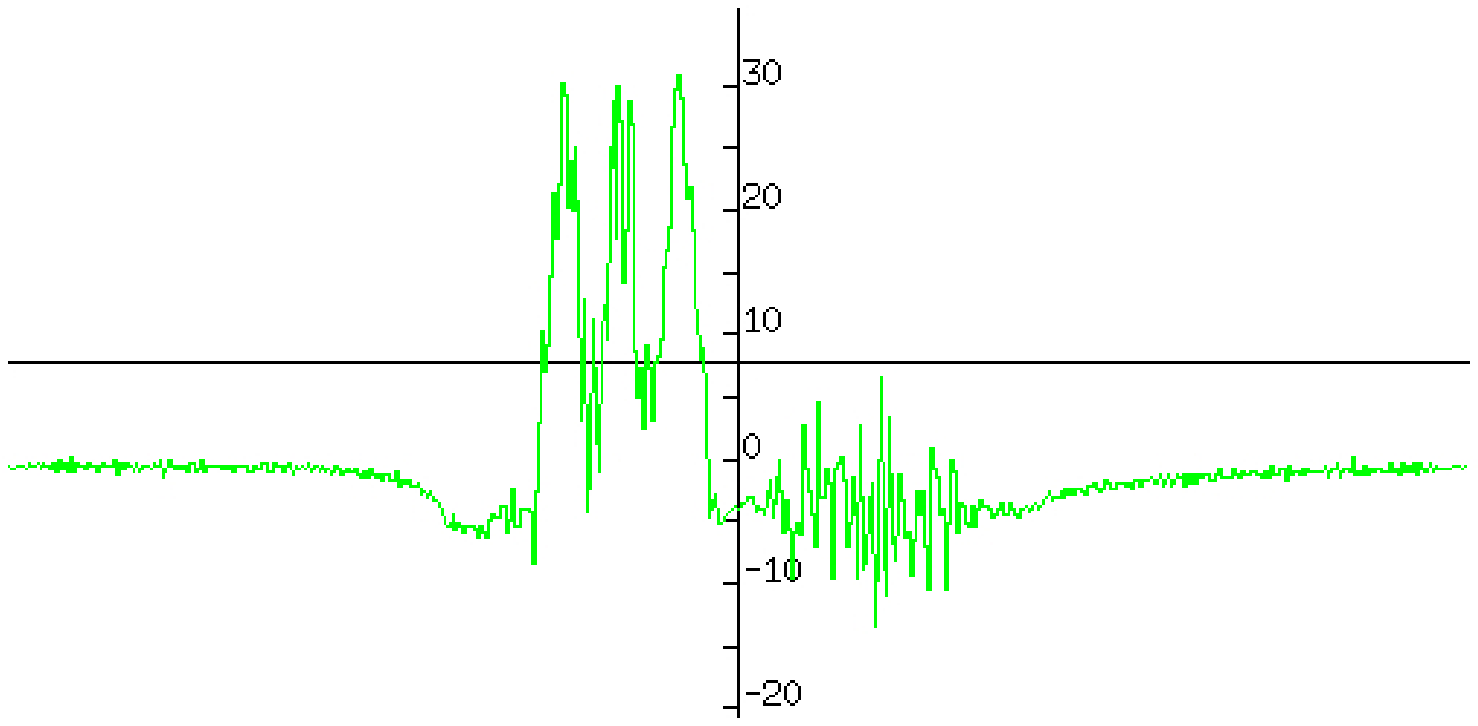}{Three antennas with positions 0, 1, and 3, with spherical Radon transform, 20\% noise\label{20Vergleich13Typ1}}{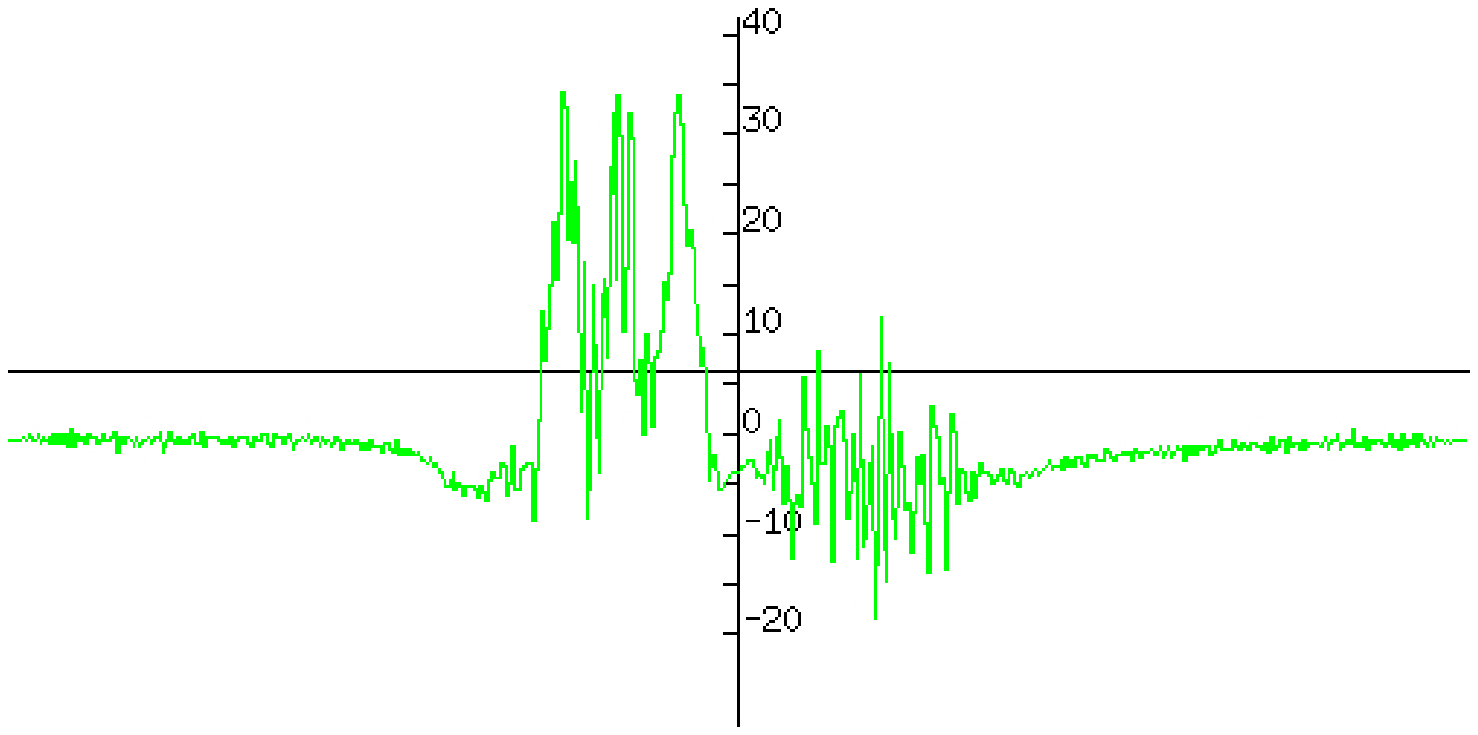}{Three antennas with positions 0, 1, and 3, with spherical Radon transform, 30\% noise\label{30Vergleich13Typ1}}

In figures \ref{Vergleich13Typ1} to \ref{30Vergleich13Typ1} the results for calculations with three antennas with distances of one and three using the inversion of the spherical Radon transform are displayed. In figure \ref{Vergleich13Typ1} the noise level is at 10\%, in figure \ref{20Vergleich13Typ1} at 20\%, and in figure \ref{30Vergleich13Typ1} at 30\%. The reconstruction in figure \ref{Vergleich13Typ1} is quite good. The amplitude of the phantom is not completely correct, but as already mentioned this is an effect of the inversion of the spherical Radon transform due to limited data. The suppression of the mirror image is also acceptable. In figure \ref{20Vergleich13Typ1} a clearly inferior reconstruction quality due to notable errors in  the amplitude of the phantom and less suppression of the mirror image can be observed. Figure \ref{30Vergleich13Typ1} is similar to figure \ref{20Vergleich13Typ1}, but overall the quality is a bit worse due to higher variances in the phantom's amplitude and the mirror image.

Figures \ref{Vergleich13Typ2} to \ref{30Vergleich13Typ1} show comparable results to figures \ref{Vergleich1Typ2} to \ref{30Vergleich1Typ1}, but the quality is clearly better. The images with $10\%$ noise are satisfactorily reconstructed.\\

In order to examine, whether the reconstruction quality - especially with respect to noise - can be enhanced with even more than three antennas, in the following the results of calculations with four antennas will be discussed.\\
Figures \ref{Vergleich138Typ2} to \ref{30Vergleich138Typ1} display reconstructions from four antennas with several noise amplitudes. Figures \ref{Vergleich138Typ2} to \ref{30Vergleich138Typ2} show only the effect of the solution of the left-right ambiguity whereas figures \ref{Vergleich138Typ1} to \ref{30Vergleich138Typ1} also include the artifacts from the inversion of the spherical Radon transform.

\asbildd{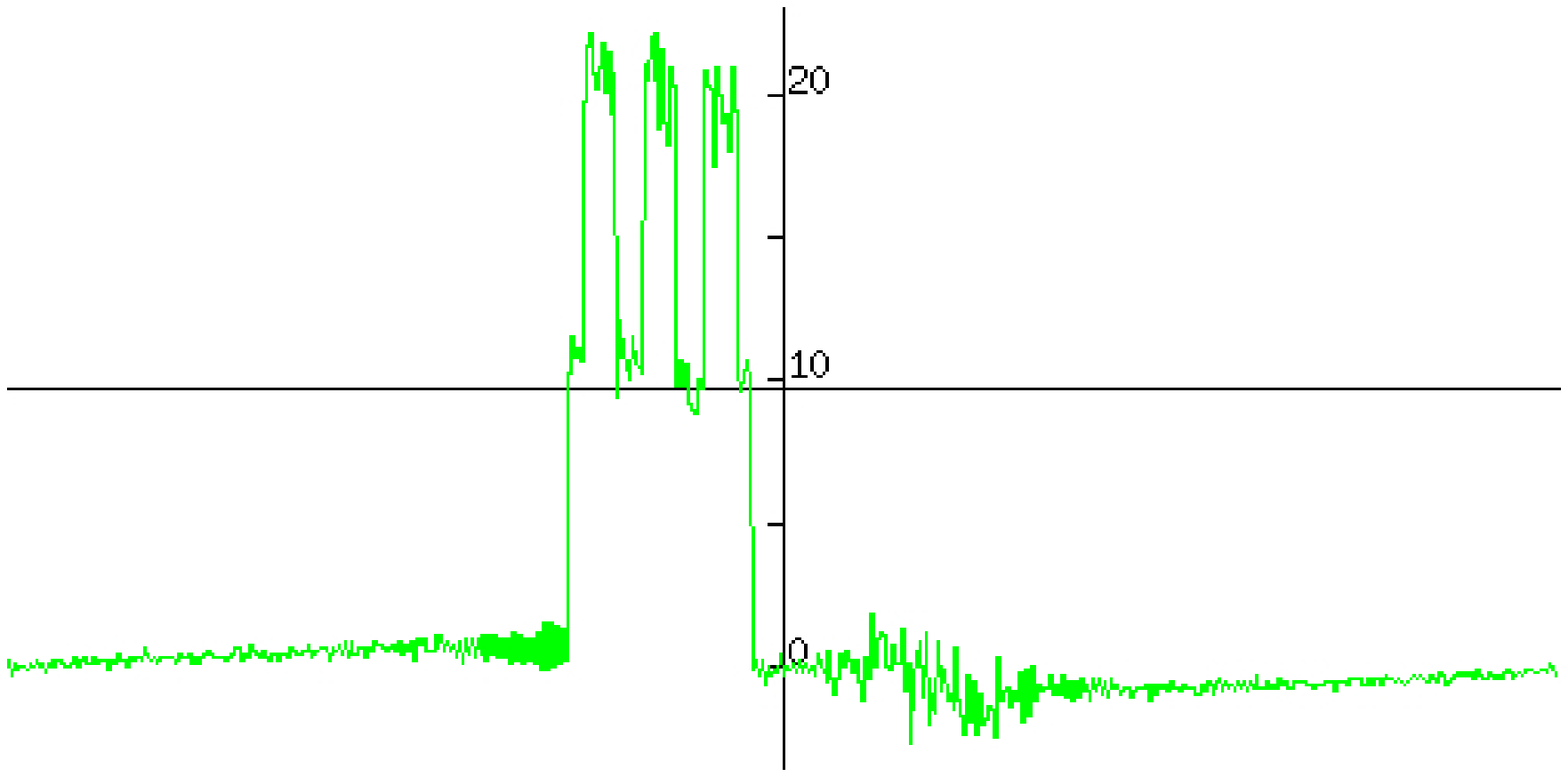}{Four antennas with positions 0, 1, 3, and 8, no spherical Radon transform, 10\% noise\label{Vergleich138Typ2}}{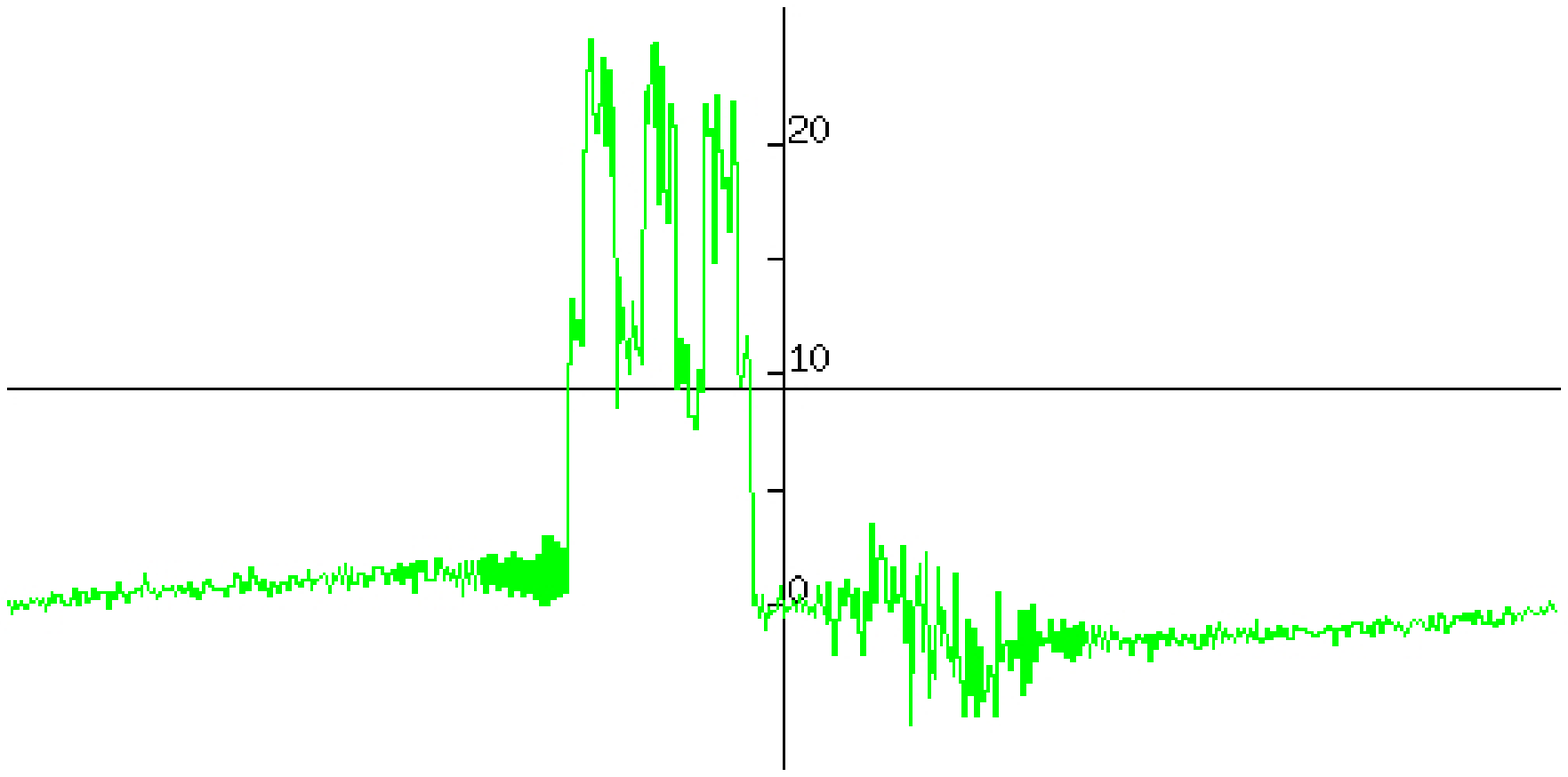}{Four antennas with positions 0, 1, 3, and 8, no spherical Radon transform, 20\% noise\label{20Vergleich138Typ2}}{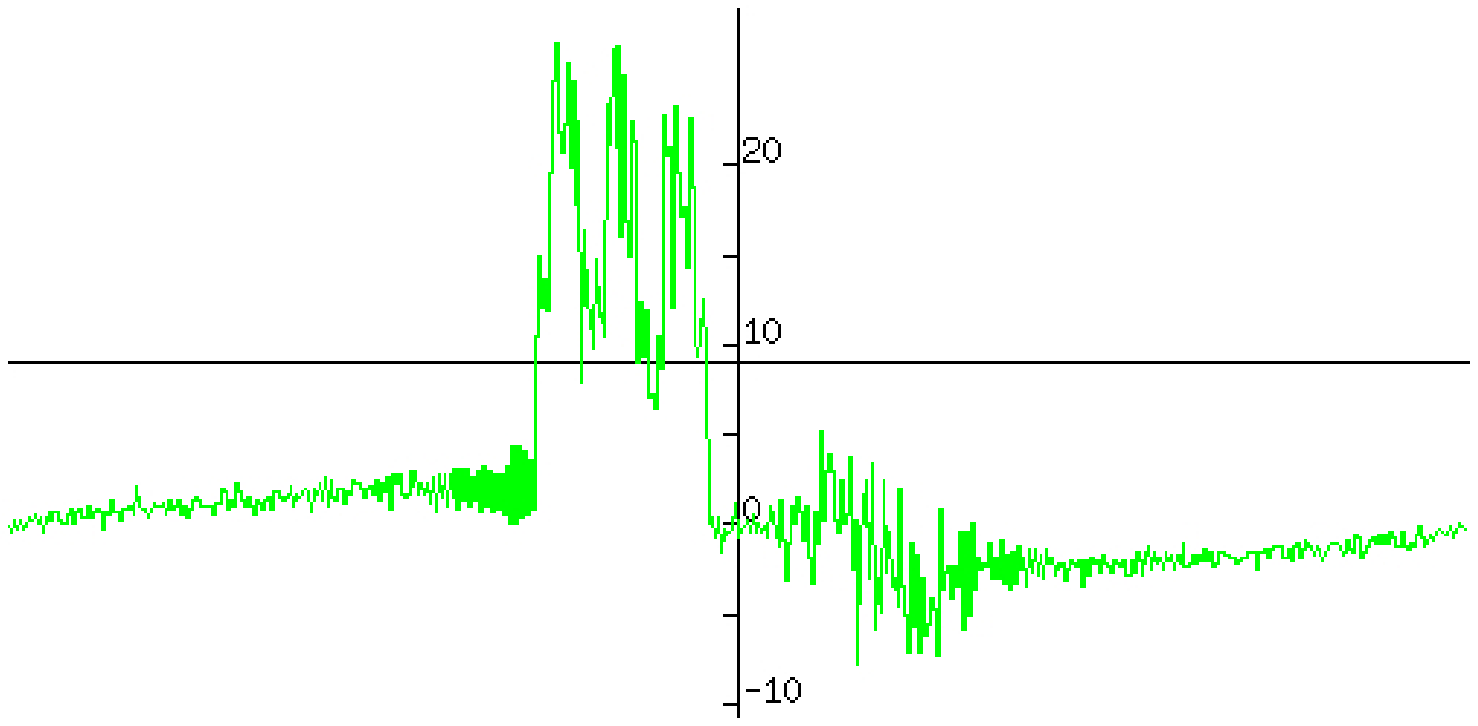}{Four antennas with positions 0, 1, 3, and 8, no spherical Radon transform, 30\% noise\label{30Vergleich138Typ2}}

Figures \ref{Vergleich138Typ2} to \ref{30Vergleich138Typ2} display reconstructions from computed even images for four antennas with distances one, three, and eight. The noise level in figure \ref{Vergleich138Typ2} is at a level of 10\%, in figure \ref{20Vergleich138Typ2} at 20\%, and in figure \ref{30Vergleich138Typ2} at 30\%. The reconstruction in figure \ref{Vergleich138Typ2} is very good. The amplitude of the phantom is recovered very well and the mirror image is quite small. Only a very slight slope can be detected in the cross section. The quality in figure \ref{20Vergleich138Typ2} declines, but nevertheless the amplitude is almost on target and the suppression of the mirror image is still acceptable. Noticable however is the increased slope in comparison to figure \ref{Vergleich138Typ2}. Figure \ref{30Vergleich138Typ2} shows an additional decline in quality, but not as severe. The amplitude deviates farther, and the slope is steeper. Most annoying however is the more strongly pronounced mirror image.

\asbildd{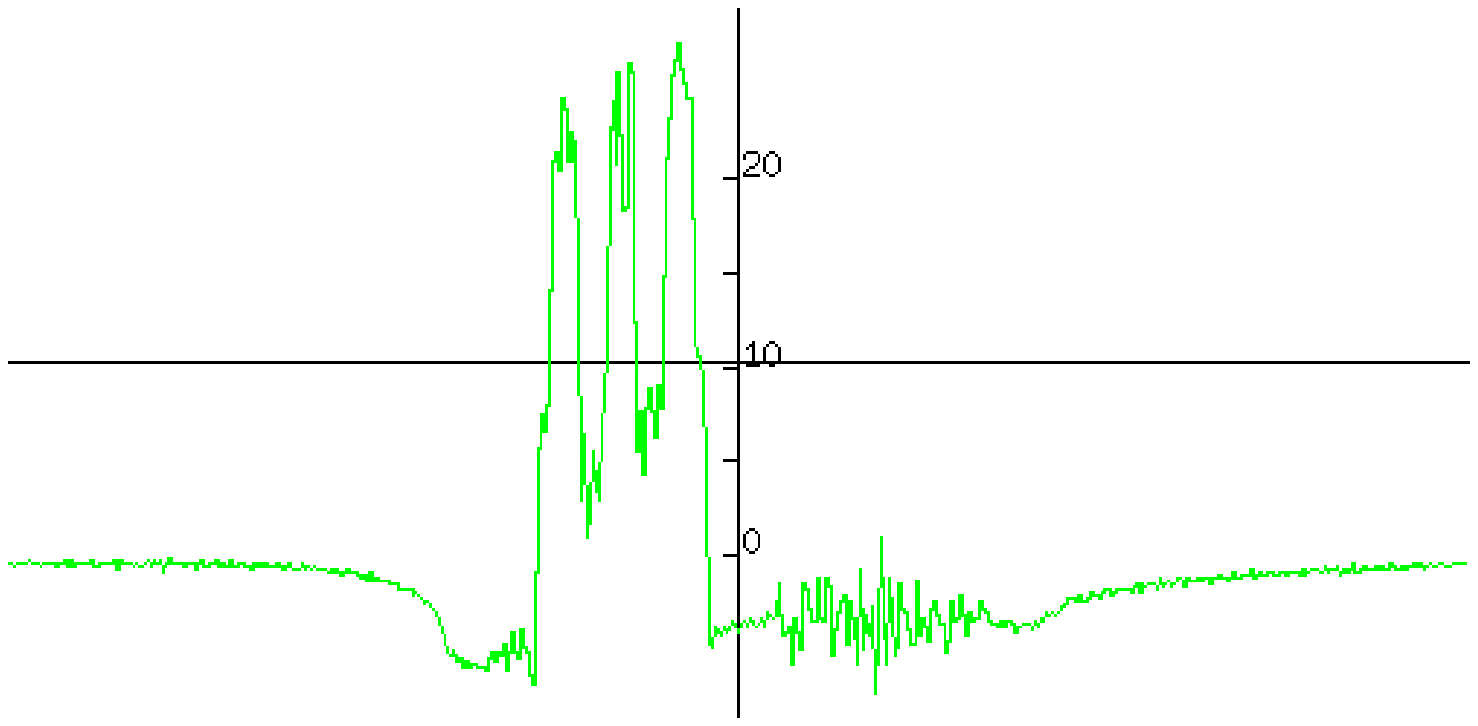}{Four antennas with positions 0, 1, 3, and 8, with spherical Radon transform, 10\% noise\label{Vergleich138Typ1}}{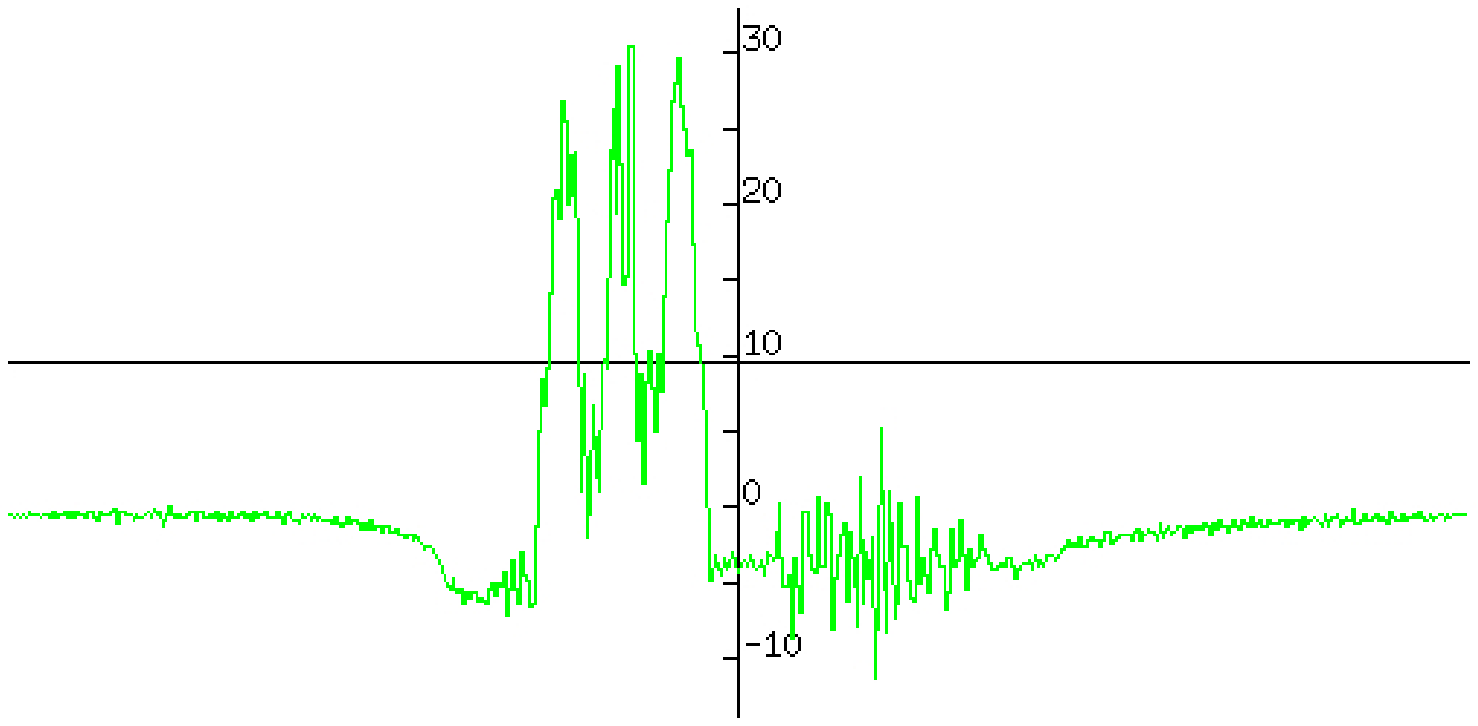}{Four antennas with positions 0, 1, 3, and 8, with spherical Radon transform, 20\% noise\label{20Vergleich138Typ1}}{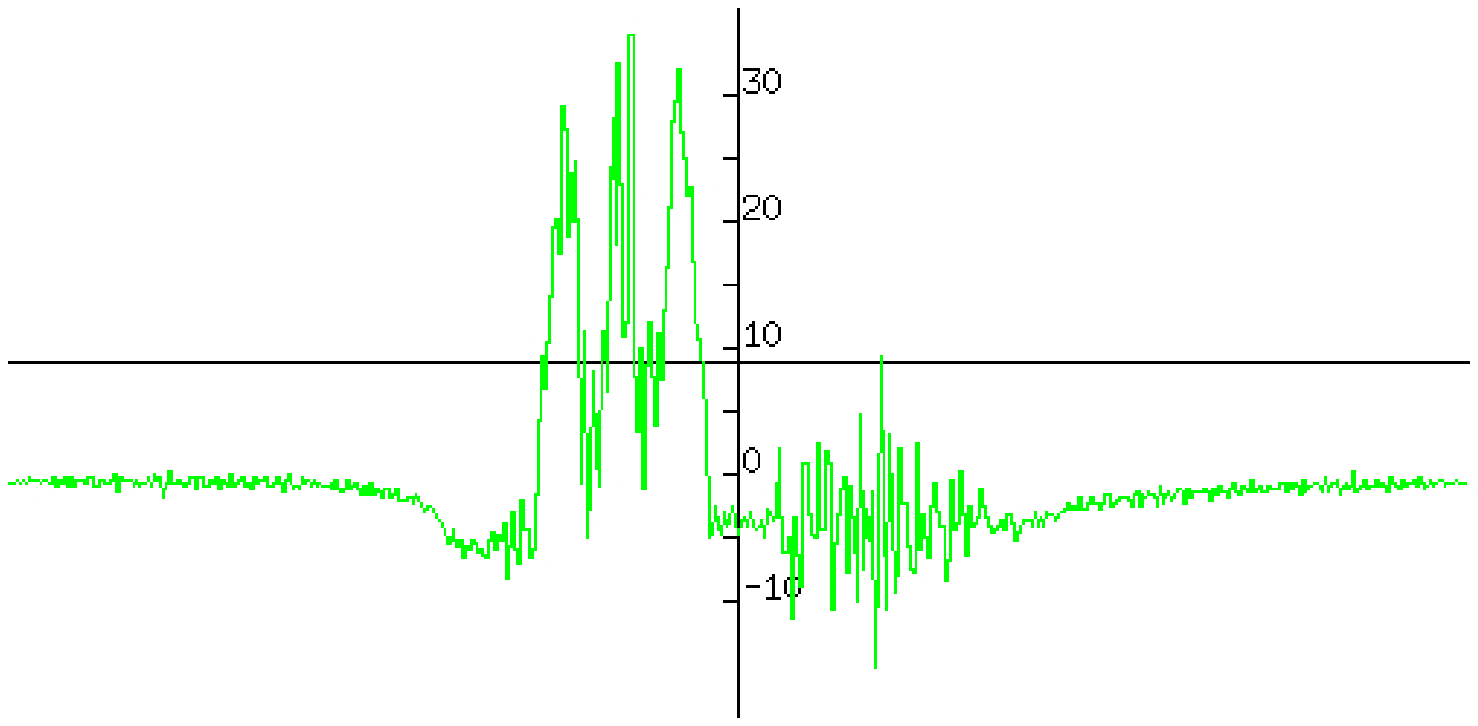}{Four antennas with positions 0, 1, 3, and 8, with spherical Radon transform, 30\% noise\label{30Vergleich138Typ1}}

Figures \ref{Vergleich138Typ1} to \ref{30Vergleich138Typ1} delineate the results from the inverted spherical Radon transforms using four antennas with distances one, three, and eight. In figure \ref{Vergleich138Typ1} the noise level is at 10\%, in figure \ref{20Vergleich138Typ1} at 20\%, and in figure \ref{30Vergleich138Typ1} at 30\%. Figure \ref{Vergleich138Typ1} looks quite good. The suppression of the mirror image is adequate and the amplitude of the reconstructed phantom is almost only plagued by the error due to the inversion of the spherical Radon transform with limited data. In figure \ref{20Vergleich138Typ1} the suppression is not as good and the amplitude of the phantom deviates farther. Finally figure \ref{30Vergleich138Typ1} displays  a slightly larger deviation in the object's amplitude and a more noticable mirror image.

Figures \ref{Vergleich138Typ2} to \ref{30Vergleich138Typ1} show a sizeable improvement with regard to figures \ref{Vergleich13Typ2} to \ref{30Vergleich13Typ1}. Now, even with $20\%$ noise the reconstructions seem quite useful and the images with $10\%$ noise look rather good.\\


But this effect can be enhanced even more.\\
Figures \ref{Vergleich13819Typ2} to \ref{30Vergleich13819Typ1} are reconstructions from five antennas with a maximum distance of 19 and varying amounts of noise. With regard to currently achievable resolutions and the wingspan of the used airplanes, this distance corresponds to an appropriate maximum antenna distance. Figures \ref{Vergleich13819Typ2} to \ref{30Vergleich13819Typ2} are computed from symmetric images. Figures \ref{Vergleich13819Typ1} to \ref{30Vergleich13819Typ1} are the results of a complete reconstruction process, including the inversion of the spherical Radon transform.

\asbildd{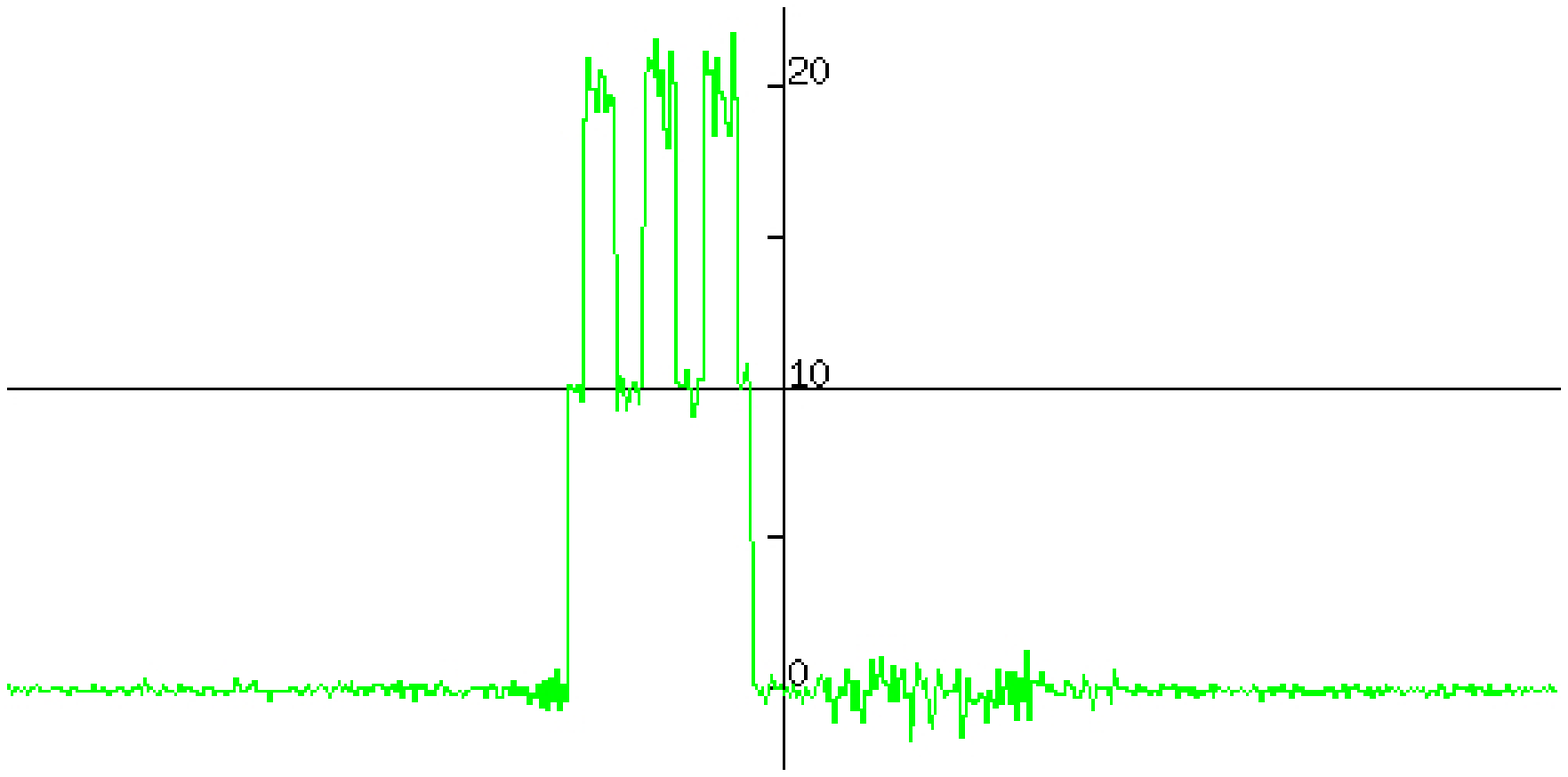}{Five antennas with positions 0, 1, 3, 8, and 19, no spherical Radon transform, 10\% noise\label{Vergleich13819Typ2}}{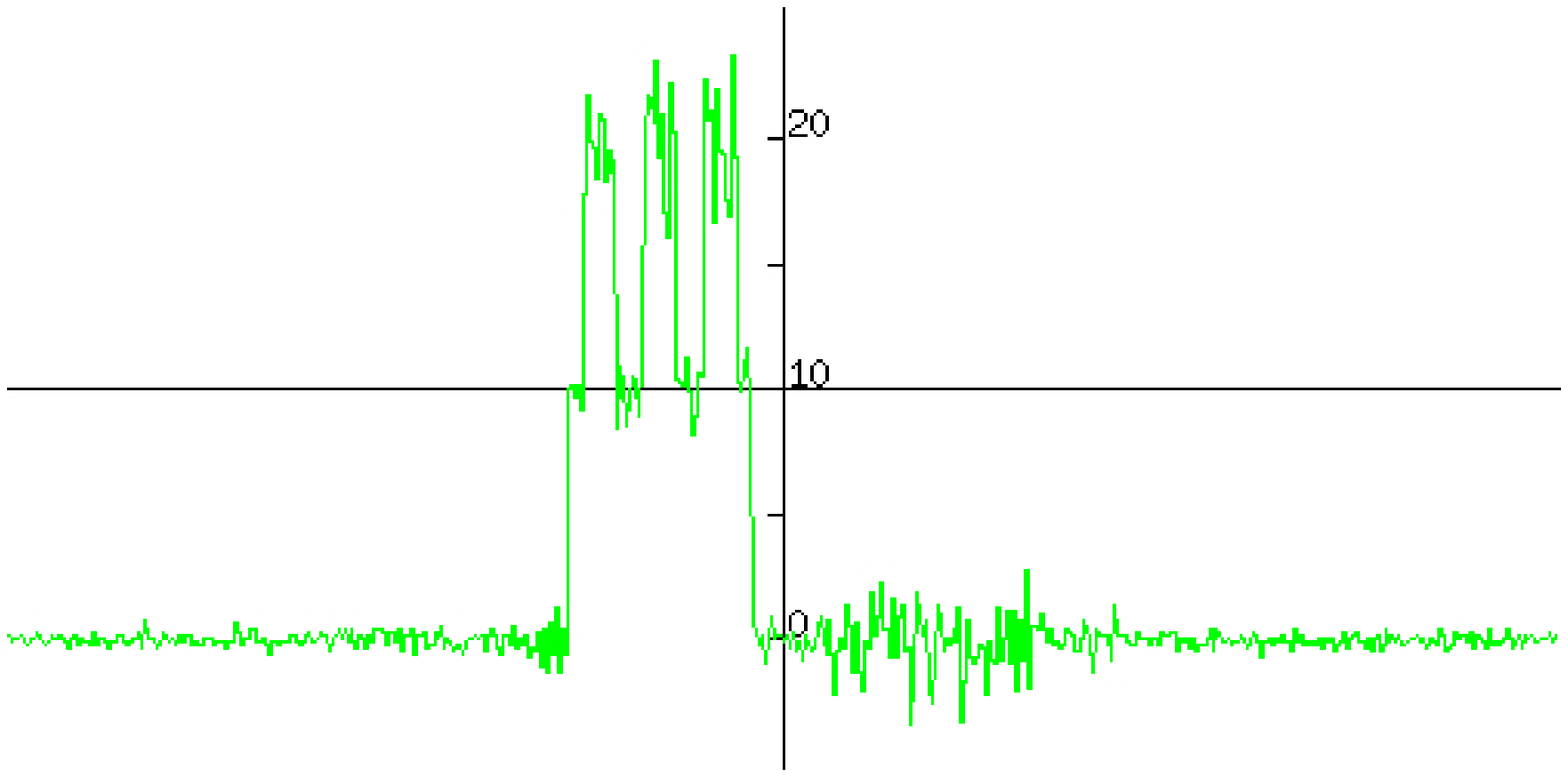}{Five antennas with positions 0, 1, 3, 8, and 19, no spherical Radon transform, 20\% noise\label{20Vergleich13819Typ2}}{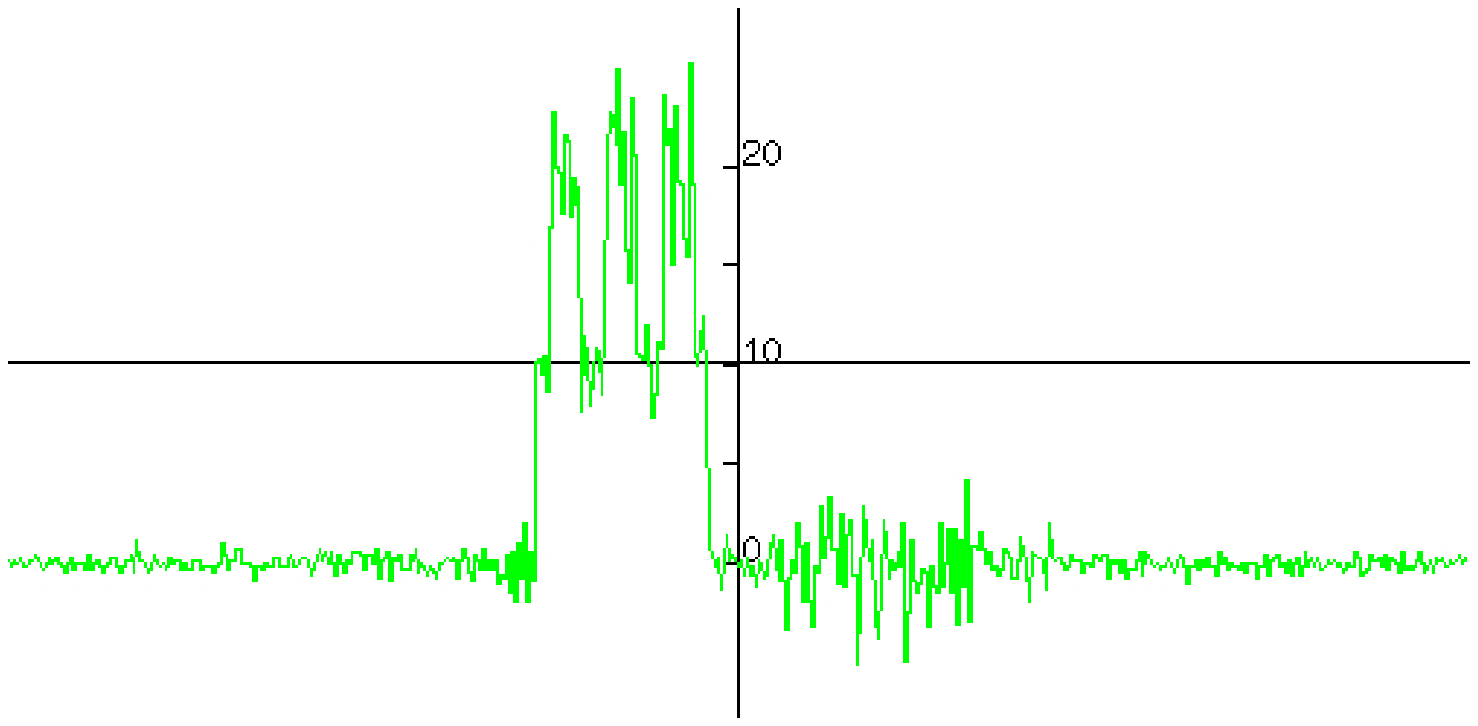}{Five antennas with positions 0, 1, 3, 8, and 19, no spherical Radon transform, 30\% noise\label{30Vergleich13819Typ2}}

Figures \ref{Vergleich13819Typ2} to \ref{30Vergleich13819Typ2} show reconstructions from computed even images with respect to five antennas with distances one, three, eight, and 19. The noise level in figure \ref{Vergleich13819Typ2} is at 10\%, in figure \ref{20Vergleich13819Typ2} at 20\%, and in figure \ref{20Vergleich13819Typ2} at 30\%. Figure \ref{Vergleich13819Typ2} exhibits an almost perfect reconstruction. The object's amplitude is recovered  very accurately, and no slope is visible. The most notable feature however is that the mirror image is almost undiscernible from the noise floor. Also in figure \ref{20Vergleich13819Typ2} a very good reconstruction can be seen. The amplitude of the phantom is quite good, and again there is no slope. The mirror image is visible, but does not stand out excessively. The reconstruction quality visible in figure \ref{30Vergleich13819Typ2} is still remarkable. The amplitude of the object is recovered closely and there is no slope. The suppression of the mirror image is also very good.

\asbildd{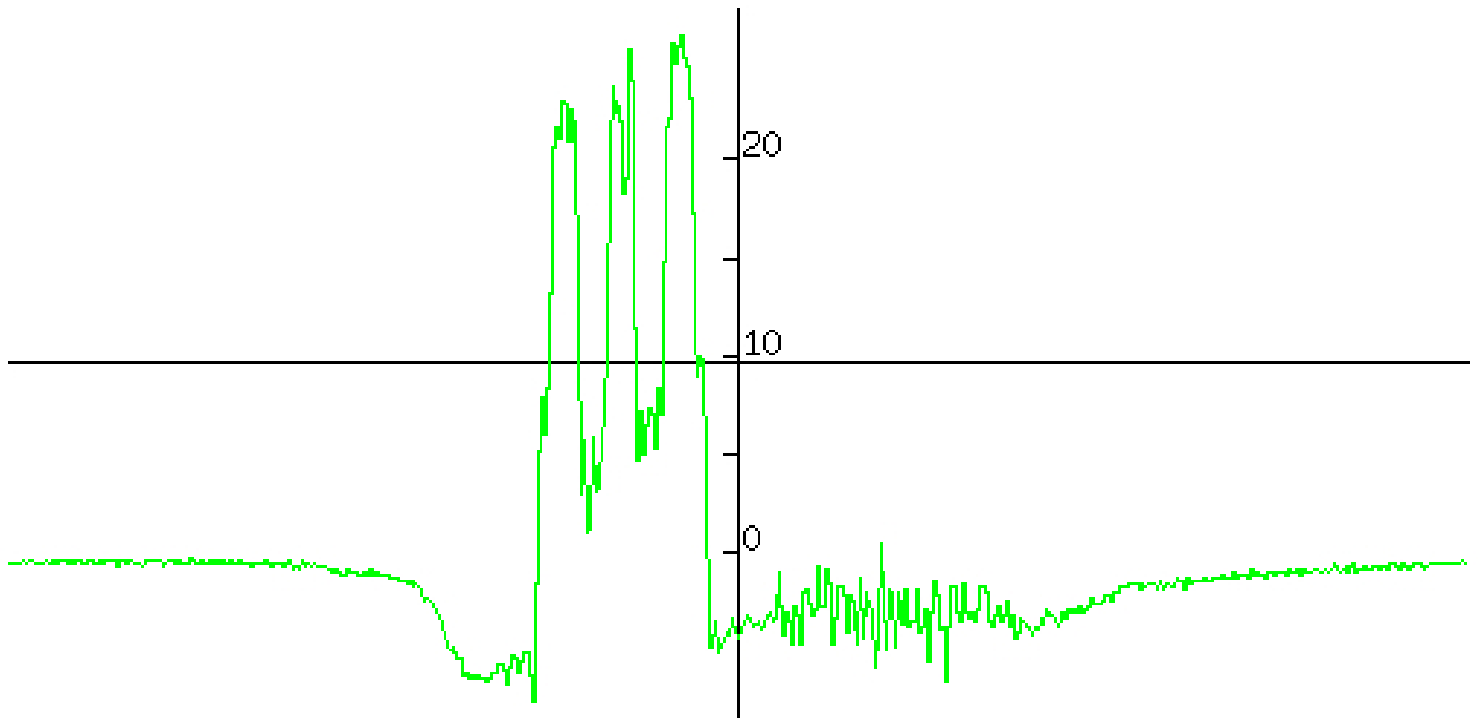}{Five antennas with positions 0, 1, 3, 8, and 19, with spherical Radon transform, 10\% noise\label{Vergleich13819Typ1}}{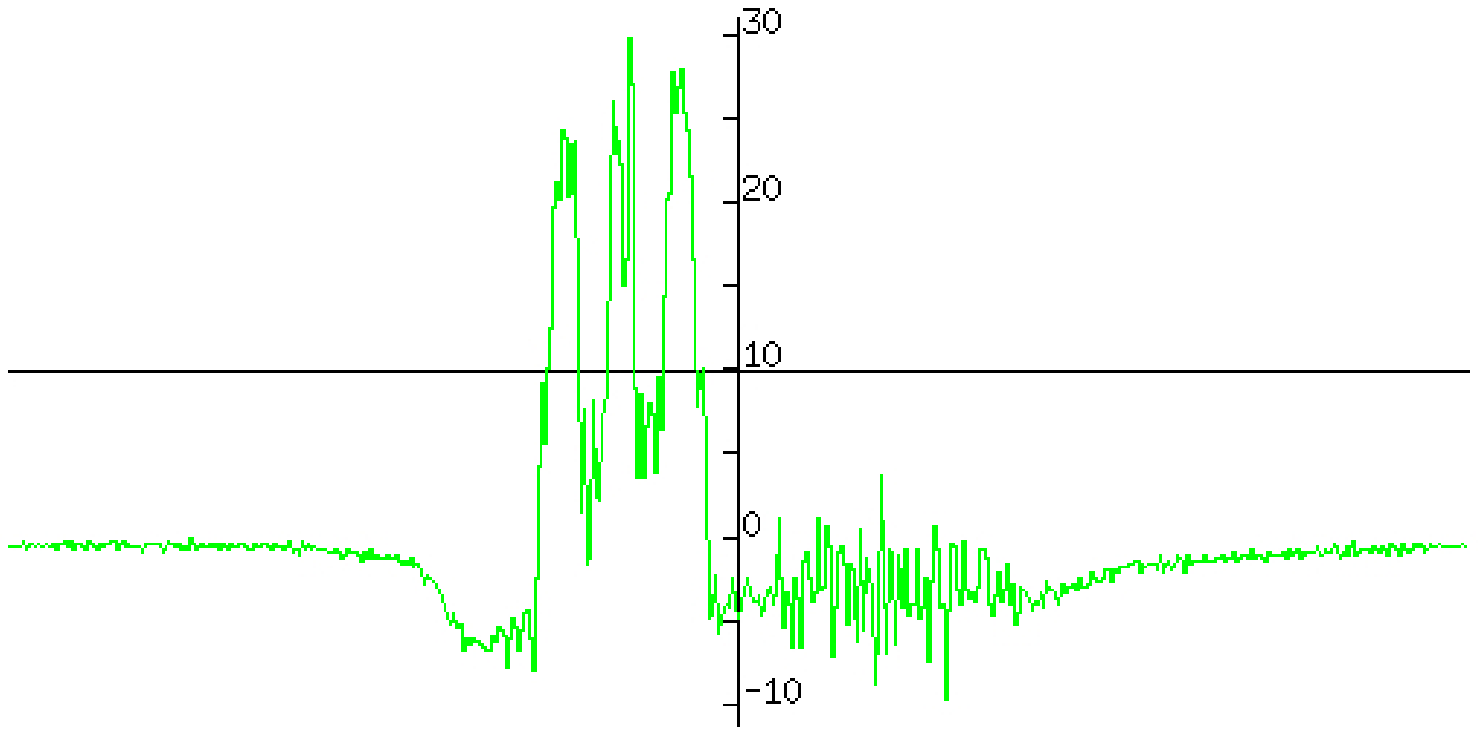}{Five antennas with positions 0, 1, 3, 8, and 19, with spherical Radon transform, 20\% noise\label{20Vergleich13819Typ1}}{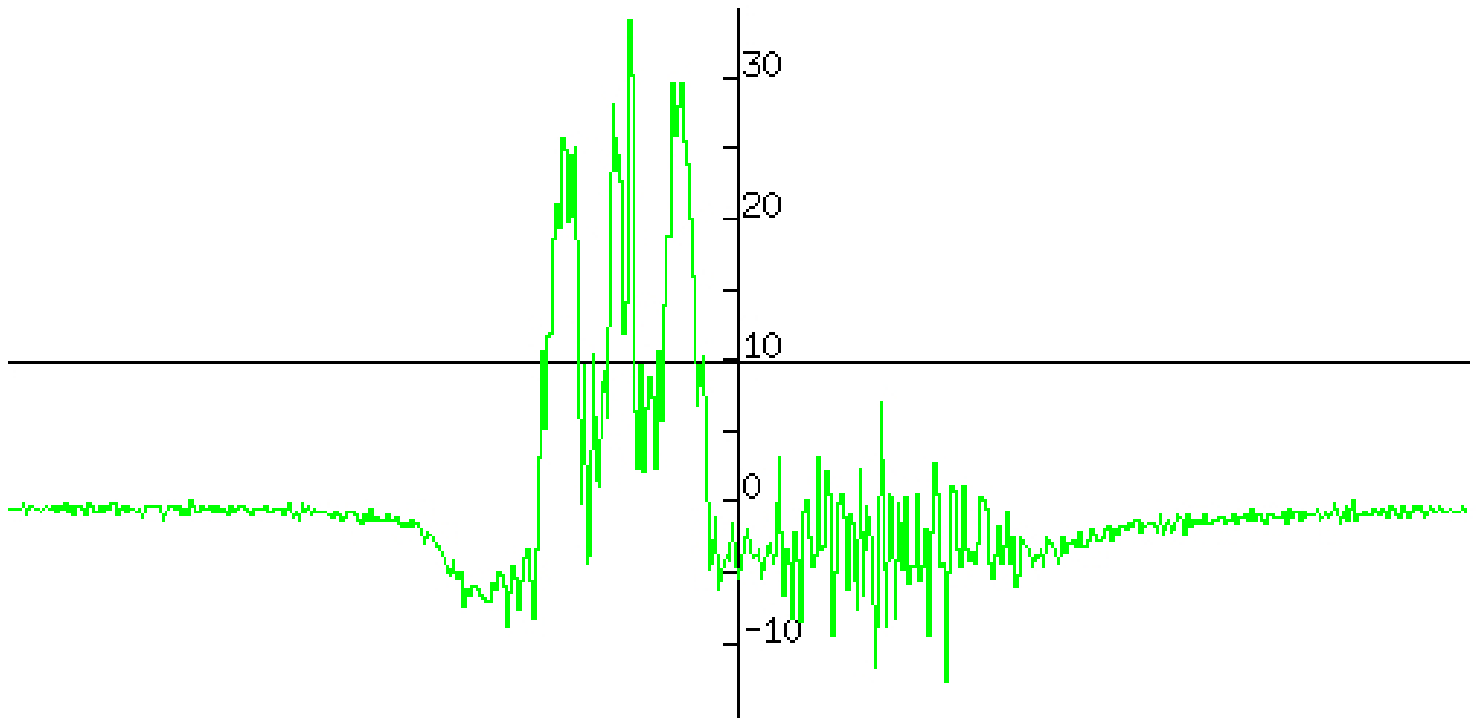}{Five antennas with positions 0, 1, 3, 8, and 19, with spherical Radon transform, 30\% noise\label{30Vergleich13819Typ1}}

Figures \ref{Vergleich13819Typ1} to \ref{30Vergleich13819Typ1} delineate reconstructions using inversions of the spherical Radon transform from five antennas with distances one, three, eight, and 19. In figure \ref{Vergleich13819Typ1} the noise level is at 10\%, in figure \ref{20Vergleich13819Typ1} at 20\%, and in figure \ref{30Vergleich13819Typ1} at 30\%. The reconstruction depicted in figure \ref{Vergleich13819Typ1} is very good. The mirror image almost drowns in the noise induced by the inversion of the spherical Radon transform with limited data. Unfortunately there is still a gradient in the amplitude of the phantom. However, as can be seen in comparison with figure \ref{Vergleich13819Typ2} this is obviously an effect of the limited data inversion, too. Figure \ref{20Vergleich13819Typ1} is still plagued by the problems of the limited data inversion, but apart from that the reconstruction is quite good. The suppression of the mirror image is still quite good. However the phantom's amplitude varies too much. Figure \ref{30Vergleich13819Typ1} displays a stronger deviance in the object's amplitude. The mirror image is also more pronounced. Therefore the reconstruction quality is a bit worse.

Figures \ref{Vergleich13819Typ2} and \ref{Vergleich13819Typ1} show only very scant signs of noise. But also figures \ref{20Vergleich13819Typ2} and \ref{20Vergleich13819Typ1} are acceptable and figures \ref{30Vergleich13819Typ2} and \ref{30Vergleich13819Typ1} show a large improvement in comparison with reconstructions that use less antennas. Figures \ref{Vergleich13819Typ2} to \ref{30Vergleich13819Typ1} show a noticable improvement in comparison with figures \ref{Vergleich138Typ2} to \ref{30Vergleich138Typ1}. The most important improvement is the absence of a slope in figures \ref{Vergleich13819Typ2} to \ref{30Vergleich13819Typ1}, but also the better suppression of the mirror image is obvious in all figures.

The results of subsection \ref{noise} show that the algorithms in theorem \ref{antisymformel} are stable. Additionally, they show that it is advantageous to employ more than two antennas and that better results can be expected for an increasing number of antennas. So if an airplane would be used to its maximal capability, then even with a large interference by noise very good images should be possible and the results should be very stable. This culminates in the insight that even for very adverse noise levels it is possible to receive good images if the whole wingspan of an airplane can be used.




\chapter{Conclusions and outlook}

In this thesis several topics connected with SAR are addressed. Most of them will help to increase image quality, but one raises questions that deserve further study.

In chapter \ref{Diplomkorrektur} an important error in a common approach to invert the spherical Radon transform was analyzed and a solution was proposed. This allowed for a promising new approach to alleviate the problem of limited data. The results of numerical simulations are very promising. To improve the results even further, it should be studied whether and what kind of decay in the approximatively continued data enhances the reconstruction quality.

Chapter \ref{Fehler} addressed the problem of limited data more thoroughly. The findings of this chapter contribute to an understanding of the artifacts that frequently appear in the numerical inversion of the spherical Radon transform. Unfortunately, it turned out that the artifacts caused by limited data are not restricted to high frequencies, as is the case in computerized tomography. Since this analysis should also hold for other models of SAR, e.g. models using the wave equation, the implications of this finding should be examined with great care. As can be seen in the numerical simulations of chapter \ref{Diplomkorrektur}, the applied regularization has a strong influence on the artifacts that appear in the reconstruction. Therefore the physical meaning of regularizations should be studied. Hopefully, this results in a reconstruction formula that minimizes artifacts and that can be relied upon, because it matches the physical conditions.

Chapter \ref{Fast} presented another new way to invert the spherical Radon transform that uses the ideas of chapter \ref{Fehler} to reduce the complexity of the problem with a projection onto orthogonal functions. With the insight gained from chapter \ref{Fehler} it might be possible to select some of these orthogonal functions to obtain good reconstructions. If only a handful of these function would be accounted for, this could also yield a fast way to invert the spherical Radon transform.

Chapter \ref{Antisym} dealt with the problem of the left-right ambiguity that causes objects that lie on one side of the flight track to appear in the reconstruction of both sides. With the postprocessing formula derived in this chapter it is possible to use reconstructions from at least two antennas to solve this problem. The numerical simulations are very promising. It can be concluded that the reconstruction quality improves with a higher number of antennas. It is however important to note that this is not only the usual effect of averaging over more data to reduce noise. This effect is amplified by a weighted summation of the data that emphasizes data less affected by noise. Therefore if an ample number of antennas would be mounted on an airplane, it should be possible to reconstruct images that show only very scant errors despite heavy noise.\\
Because of the jitter of the airplane a modification of the formulas in this chapter should be studied. It should be possible to obtain a similar formula that uses the data instead of the reconstructed images as input. However it would be necessary to mirror the data so it becomes an even function that has a meaning for negative radii. Then the left-right ambiguity in the data could be treated in a similar way as in chapter \ref{Antisym}, and for each side of the flight track a separate inversion could be performed. With this modification the fluctuation in airplane position and heading would probably no longer pose a problem because it is very small for one send/receive cycle. An additional advantage would be the fact, that the inversion of the spherical Radon transform would have to be performed only twice regardless of the number of antennas. This should speed up the reconstruction considerably.

\onecolumn

\end{sloppypar}

\end{document}